\documentclass{article}
\usepackage[english]{babel}
\usepackage[utf8]{inputenc}
\usepackage[T1]{fontenc}
\usepackage{graphicx}
\usepackage[outputdir=build, cache=false, draft=true]{minted}
\usepackage{stmaryrd}
\usepackage{bm} 
\usepackage{mathrsfs}
\usepackage{mathtools}
\usepackage{units}
\usepackage{tensor}
\usepackage{accents}
\usepackage{xspace}
\usepackage{enumitem}
\usepackage[obeyFinal,textsize=footnotesize,textwidth=55pt,
backgroundcolor=white,linecolor=black,bordercolor=black]{todonotes}
\usepackage{dsfont}
\usepackage{hyperref}
\usepackage{bookmark}
\usepackage{cleveref}
\usepackage{amssymb}
\usepackage{amsmath}
\usepackage{helvet}
\usepackage{mathrsfs}
\usepackage{authblk}
\usepackage[numbers]{natbib}     
\usepackage{pgf, tikz, pgfplots}
\usepackage[top=2.5cm, bottom=2.5cm, left=3cm, right=3cm]{geometry}
\usepackage{algpseudocode}
\usepackage{booktabs}   
\pgfplotsset{compat=1.15}
\usetikzlibrary{shapes.geometric, arrows,backgrounds,fit,positioning}

\newcommand{\R}{\mathbb R}

\newcommand{\T}{\mathcal T}    
\newcommand{\E}{\mathcal E}    

\DeclareMathOperator{\e}{e}    
\DeclareMathOperator*{\argmin}{arg\,min}   

\DeclarePairedDelimiter{\norm}{\lVert}{\rVert}  
\DeclarePairedDelimiter{\snorm}{\lvert}{\rvert}  
\newcommand{\dy}{\, \mathrm{d}y}    
\renewcommand{\phi}{\varphi}

\renewcommand{\epsilon}{\varepsilon}

\newcommand{\Pcal}{\mathcal P}    

\newcommand{\scal}[2]{\left(#1,#2\right)}    

\renewcommand{\tilde}[1]{\widetilde{#1}}

\newcommand{\comment}[1]{}  
\newcommand{\bw}{\mathrm{bw}}

\tikzset{block/.style={rectangle, draw, text width=4.2em, text centered, minimum height=3em, rounded corners=5pt},
title/.style={font=\bfseries},
line/.style={draw, -latex'},
int/.style={draw, minimum size=2em, font=\upshape}}

\title{An a posteriori error estimator for the spectral fractional power of the Laplacian}

 \author[1]{Raphaël Bulle}

 \author[4,5,6]{Olga Barrera}

 \author[,1]{Stéphane P. A. Bordas\thanks{Corresponding author: \texttt{stephane.bordas@uni.lu}}}

 \author[2,3,7]{Franz Chouly}

 \author[1]{Jack S. Hale}
 
 \affil[1]{Institute of Computational Engineering, University of Luxembourg, 6
 Avenue de la Fonte, 4362 Esch-sur-Alzette, Luxembourg.}
 \affil[2]{Université Bourgogne Franche-Comté, Institut de Mathématiques de
 Bourgogne, 21078 Dijon, France.}
 \affil[3]{Center for Mathematical Modeling and Department of Mathematical
 Engineering, University of Chile and IRL 2807 – CNRS, Santiago, Chile.}
\affil[4]{School of Engineering Computing and Mathematics, Oxford Brookes
University, Oxford, UK.}
\affil[5]{Department of Engineering Science, University of Oxford, Oxford, UK.}
\affil[6]{Department of Medical Research,China Medical University Hospital,
China Medical University, Taichung, Taiwan.}
\affil[7]{Departamento de Ingenería Matématica, CI²MA, Universidad de
Concepción, Casilla 160-C, Concepción, Chile.}

\begin{document}
\maketitle
\begin{center}
    \footnotesize
    \textbf{Abstract}\\
    \vspace{10pt}
    \begin{minipage}{0.8\textwidth}
        We develop a novel a posteriori error estimator for the $L^2$ error
        committed by the finite element discretization of the solution of the
        fractional Laplacian.
        Our a posteriori error estimator takes advantage of the semi--discretization
        scheme using rational approximations which allow to reformulate the
        fractional problem into a family of non--fractional parametric problems.
        The estimator involves applying the implicit Bank--Weiser error estimation
        strategy to each parametric non--fractional problem and reconstructing the
        fractional error through the same rational approximation used to compute the
        solution to the original fractional problem.
        In addition we propose an algorithm to adapt both the finite element mesh and the rational scheme in order to balance the discretization errors.
        We provide several numerical examples in both two and three-dimensions
        demonstrating the effectivity of our estimator for varying fractional powers
        and its ability to drive an adaptive mesh refinement strategy.
    \end{minipage}
\end{center}

{\footnotesize
\noindent \textbf{Published journal article:} \url{https://doi.org/10.1016/j.cma.2023.115943}\\
}

\section{Introduction}

Fractional partial differential equations (FPDEs) have gained in popularity during the last two decades and are now
applied in a wide range of fields\ \cite{Lischke2018} such as anomalous
diffusion\ \cite{Bonito2018a,Carlson2019,Defterli2015,Duo2019,Nochetto2015},
electromagnetism and geophysical electromagnetism\ \cite{Bonito2020a,Weiss2019}, phase fluids
\cite{Ainsworth2017,Akagi2016,Du2019}, porous
media\ \cite{Akagi2016,Caputo2000,Barrera2020,Bulle2021}, quasi-geostrophic
flows \cite{Bonito2021a} and spatial statistics\
\cite{Bolin2017,Lindgren2011a}.

The main interest in fractional models lies in their ability to reproduce
non-local behavior with a relatively small number of parameters \cite{Atangana2018,Caputo1967}.
While this non-locality can be interesting from a modeling perspective, it also constitutes an ongoing challenge for numerical methods since applying standard approaches naturally leads to large dense linear systems that are
computationally intractable.

In the last decade various numerical methods have been derived in order to
circumvent the main issues associated with the application of standard numerical
methods to FPDEs.
The two main ones being the non--locality leading to dense
linear systems and, for some particular definitions of the fractional operator,
the evaluation of singular integrals \cite{Ainsworth2017a,Ainsworth2018}.

We focus on discretization schemes based on finite element methods,
other methods can be found e.g.\ in
\cite{Aceto2019,Higham2013,Podlubny2009}.
Among the methods addressing the above numerical issues, we can cite: methods to
efficiently solve eigenvalue problems \cite{Carlson2019}, multigrid methods for
performing efficient dense matrix--vector products\
\cite{Ainsworth2017a,Ainsworth2018}, hybrid finite element--spectral schemes\
\cite{Ainsworth2018a}, Dirichlet-to-Neumann maps (such as the
Caffarelli--Silvestre extension)\
\cite{Antil2015,Caffarelli2007,Chen2015,Faustmann2020,Nochetto2015,Stinga2010},
semigroups methods\ \cite{Cusimano2020,Cusimano2018,Stinga2019}, rational
approximation methods\ \cite{Aceto2017,Harizanov2020,Hofreither2020},
Dunford--Taylor integrals\
\cite{Bolin2017,Bonito2017,Bonito2016,Bonito2019,Bonito2013,Borthagaray2019a,Gavrilyuk2004,Hofreither2020}
(which can be considered as particular examples of rational approximation
methods) and reduced basis methods\ \cite{Danczul2019,Danczul2020,Dinh2021}.

Although we focus exclusively on the spectral definition of the fractional
Laplacian, there is no unique definition of the fractional power of the
Laplacian operator.
The three most frequently found definitions of the fractional Laplacian are: the
integral fractional Laplacian, defined from the principal value of a singular
integral over the whole space $\R^d$\
\cite{Ainsworth2017a,Ainsworth2018,Bonito2019a,Caffarelli2007,DElia2013}, the
regional fractional Laplacian, defined by the same singular integral but over a
bounded domain only\ \cite{Chen2018,Duo2019,Fall2020,Mou2015} and the spectral
fractional Laplacian, defined from the spectrum of the standard Laplacian over a
bounded domain\
\cite{Ainsworth2018a,Antil2018,Banjai2019,Cusimano2020,Harizanov2020,Meidner2018}.
The different definitions are equivalent in the entire space $\R^d$, but this is
no longer the case on a bounded domain\
\cite{Bonito2018a,Duo2019,Kwasnicki2017,Lischke2018}. 
These definitions lead to significantly different mathematical problems
associated with infinitesimal generators of different stochastic processes
\cite{Lischke2018,Duo2019}.

Efficient methods for solving fractional problems typically rely on a
combination of different discretization methods.
For example, \cite{Bonito2013,Harizanov2016} which are also the foundation of this work,
combine a rational sum representation of the spectral fractional Laplacian with a standard finite element method in space.
Both the quadrature scheme and the finite element method induce discretization
errors.
In order to achieve a solution to a given accuracy while avoiding wasted
computational time, these errors need to be balanced.

A priori error estimation has been tackled for some definitions of the
fractional Laplacian, such as the integral
Laplacian\
\cite{Acosta2017,Ainsworth2017a,Ainsworth2018a,Bonito2018a,Bonito2013,Gimperlein2019}
and the spectral fractional Laplacian\
\cite{Antil2015,Antil2018,Banjai2019,Bonito2013,Meidner2018,Nochetto2015}.
Unlike the standard Laplacian equation, solutions to the fractional Laplacian
problems often exhibit strong boundary layers even for smooth data,
particularly when the fractional power is low \cite{Grubb2015}.
These singularities lead to computational difficulties and have to be taken into
account using, for example a priori geometric mesh refinement towards the
boundary of the domain
\cite{Acosta2017,Banjai2020,Bonito2019a,Borthagaray2019a,Faustmann2020}, or
partition of unity enrichment \cite{Bordas2016}.
We emphasize that \cite{Bonito2013} contains already an a priori error analysis
in the $L^2$ norm for a combined rational sum finite element method that we use in
this work.

A posteriori error estimation has also been considered in the literature on
fractional equations.
A simple residual based estimator is proposed for the integral fractional
Laplacian in\ \cite{Ainsworth2017a}.
A similar idea is used in the context of non-local variational inequalities
in\ \cite{Gimperlein2019,Nochetto2010}.
Gradient-recovery based a posteriori error estimation has been developed in
the context of fractional differential equations in\ \cite{Zhao2017}.
In\ \cite{Bonito2018a,Chen2015} the authors present another estimator, based on
the solution to local problems on cylindrical stars, for the integral fractional
Laplacian discretized using the Caffarelli--Silvestre extension.
A weighted residual estimator is derived in\ \cite{Faustmann2019} in the same
context.

\section{Contribution}
The main contribution of this work is the derivation of a novel a posteriori
error estimator for the combined rational finite element approximation of the
spectral fractional Laplacian. It is a natural a posteriori counterpart to the
a priori results developed in \cite{Bonito2013}. To our knowledge, our work is
the first to propose an a posteriori error estimator for the spectral
fractional Laplacian discretized using rational approximation techniques (as
opposed to the Caffarelli-Silvestre extension).

Our work starts with rational approximation--based discretization methods.
We are particularly interested in two of these methods: a method based on the quadrature rule for the Dunford--Taylor integral proposed in the seminal work\ \cite{Bonito2013} ---we will refer to this method as "the BP method" in the following--- and a method based on Best Uniform Rational Approximations, ---"the BURA method"--- applied for the first time to fractional partial differential equations in\ \cite{Harizanov2016}.
These methods decompose the original fractional problem into a set of independent parametric
non--fractional problems.
From this point we develop an associated set of independent non--fractional a
posteriori error estimation problems.
We compute the Bank--Weiser hierarchical estimators \cite{Bank1985} of the error
between each non--fractional parametric problem solution and its finite element
discretization, then the fractional problem discretization error is estimated by
the sum of the parametric contributions via the rational approximation.
In addition, we propose an algorithm to adaptively refine the rational scheme.
When two discretization methods are combined (here a rational scheme and a finite element method), the respective discretization errors must be balanced in order to prevent computational resources waste.
The algorithm we propose is based on a computable a priori estimator for the rational scheme error and rational schemes are chosen on--the--fly to balance the rational and finite element error estimators.

Our method leads to a fully local and parallelizable solution technique
for the spectral fractional Laplacian with computable $L^2$ error.
Our method is valid for any finite element degree (however, for the sake of
brevity we do not show results with higher degree finite elements) and for one, two and three dimensional problems\ \cite{Bulle2021a}. 

We implement our method in DOLFINx \cite{Habera2020}, the new problem solving
environment of the FEniCS Project\ \cite{Alnaes2015}.
A simple demonstration implementation is included in the supplementary
material.
We show numerical results demonstrating that the estimator can correctly
reproduce the a priori convergence rates derived in\ \cite{Bonito2013}.
Our newly developed error estimator is then used to steer an adaptive mesh
refinement algorithm, resulting in improved convergence rates for small
fractional powers and strong boundary layers.

\section{Motivation}
Given a fractional power $s$ in $(0,1)$ and a rational approximation $\mathcal
Q_s(\lambda)$ of the function $\lambda^{-s}$, it is possible to construct a
semi-discrete approximation $u_{\mathcal Q_s}$ of the solution $u$ to a fractional Laplace
equation as a weighted sum of solutions $(u_l)_l$ to non--fractional parametric
problems.
Then, a fully discrete approximation of $u$ is obtained by discretizing the
parametric solutions $(u_l)_l$ using a finite element method.

An a posteriori error estimator is then computed as the weighted sum of the
Bank--Weiser estimators of the error between each $u_l$ and its finite element
discretization.
As we will see in the following, the resulting numerical scheme is simple and
its implementation in code is straightforward.
Furthermore it maintains the appealing parallel nature of
rational approximation schemes \cite{Bonito2013,Gavrilyuk2004,Harizanov2020}.

We remark on why we have chosen to use the Bank--Weiser type error estimator, as
opposed to one of the many other error estimation strategies, e.g. explicit
residual, equilibrated fluxes, or recovery-type estimators (see
\cite{Ainsworth2000,Carstensen2010} and references therein).
In the case of fractional powers of the Laplacian operator, the resulting set of
parametric problems consists of singularly--perturbed reaction--diffusion
equations.
It has been proven in\ \cite{Verfurth1998} that the Bank--Weiser estimator is
robust with respect to the coefficients appearing in these parametric problems
when the error is measured in the natural norm.
To our knowledge, no such robustness, which our numerical experiments do
indicate, has been established for the $L^2$-norm for the Bank--Weiser
estimator.
Nevertheless, our numerical experiments indicate that this does appear to be the case.
Moreover, the Bank--Weiser estimator can be straightforwardly applied to
higher-order finite element methods and higher-dimension problems.
In addition, its computational stencil is highly local which is particularly
appealing for three-dimensional problems see e.g.\ \cite{Bulle2021a}.
Finally, our choice of the Bank--Weiser estimator is also justified in \cref{subsubsec:heuristics}.

In this work we focus on error estimation in the $L^2$ norm, the estimation of
the error in the `natural' fractional norm is the topic of ongoing work.
For simplicity, we only consider fractional powers of the Laplacian with
homogeneous Dirichlet boundary conditions.

\section{Problem statement}
For any subset $\omega$ of $\overline \Omega$ we denote $L^2(\omega)$ the space
of square integrable functions on $\omega$ and $\scal{\cdot}{\cdot}_{\omega}$
its usual inner product.
Let $H^1(\omega)$ be the Sobolev space of functions with first order weak
derivatives in $L^2(\omega)$.
The space $H^1(\omega)$ is endowed with the usual inner product $\scal{\nabla
\cdot}{\nabla \cdot}_{L^2(\omega)} + \scal{\cdot}{\cdot}_{L^2(\omega)}$.
We will omit the dependence in $\omega$ in the subscripts when $\omega= \Omega$.
We will make use of the notation $\partial v/\partial n := \nabla v
\cdot n$ for the normal derivative of a smooth enough function $v$.
We denote $H^1_0(\Omega)$ the subspace of functions in $H^1(\Omega)$ with a
zero trace on $\Gamma$.

We consider the family of eigenfunctions $\{\psi_i\}_{i=1}^{\infty} \subset
H^1_0(\Omega)$ of the standard Laplacian operator with uniform zero Dirichlet
boundary condition on $\Omega$ as well as the corresponding family of
eigenvalues $\{\lambda_i\}_{i=1}^{\infty}$ such that
\begin{equation}
    \scal{\nabla \psi_i}{\nabla v} = \lambda_i \scal{\psi_i}{v},\quad \forall v \in H^1_0(\Omega),\ \forall i \in [\![1,+\infty[\![.
\end{equation}
We assume the Laplacian eigenvalues are sorted in increasing order and we assume $\lambda_0 \in \R$
is a lower bound of the spectrum
\begin{equation}
    \lambda_0 \leqslant \lambda_1 \leqslant \cdots \leqslant \lambda_i \leqslant
    \lambda_{i+1} \leqslant \cdots
\end{equation}
The family $\{\psi_i\}_{i=1}^{\infty}$ is an orthonormal basis of $L^2(\Omega)$ \cite{Allaire2007}.
For $s$ in $(0,1)$ we introduce the spectral fractional Sobolev space $\mathbb
H^s$ and its natural norm
\begin{equation}
    \mathbb{H}^s := \left\{ v \in L^2(\Omega),\ \sum_{i=1}^{\infty} \lambda_i^s
    \scal{v}{\psi_i}^2 < \infty \right\},\qquad 
    \norm{v}_{\mathbb{H}^s}^2 := \sum_{i=1}^{\infty} \lambda_i^s
\scal{v}{\psi_i}^2.
    \label{eq:spectral_sobolev}
\end{equation}
Especially, for $0 \leqslant s \leqslant 1$ we have $H^1_0(\Omega) = \mathbb
H^1(\Omega) \subseteq \mathbb H^s(\Omega) \subseteq L^2(\Omega)
=: \mathbb H^0(\Omega)$ and the norm $\norm{\cdot}_{\mathbb H^s}$ coincide with
$\norm{\cdot}_{L^2}$ when $s=0$ and with $\snorm{\cdot}_{H^1}$ when $s=1$.
In the following, for a function $v \in L^2(\Omega)$ we will denote $v_i := \scal{v}{\psi_i}$ for all $i \geqslant 1$.

\subsection{The spectral fractional Laplacian}
Let $s$ be a real number in $(0,1)$ and $f$ be a given function in
$L^2(\Omega)$.
We consider the following fractional Laplacian problem: we look for a
function $u$ such that
\begin{equation}
    \mathcal L^s u = f,
    \label{eq:fractional_strong_form}
\end{equation}
where $\mathcal L$ stands for the Laplacian operator in $\Omega$, with homogeneous Dirichlet boundary conditions on $\Gamma$.
Let us consider $\{\lambda_i,\psi_i\}_{i=1}^{+\infty} \subset \R^{+*} \times L^2(\Omega)$, the spectrum of $\mathcal L$, defined by the following generalized eigenvalue problem
\begin{equation}
    \scal{\nabla \psi_i}{\nabla v} = \lambda_i \scal{\psi_i}{v},\quad \forall v \in H^1_0(\Omega),\ \forall i \in \llbracket 1,+\infty \llbracket.
\end{equation}
The solution $u$ of \cref{eq:fractional_strong_form} is defined using the
spectrum of the standard Laplacian \cite{Antil2018}
\begin{equation}
    \scal{u}{v} := \scal{\mathcal L^{-s} f}{v} = \sum_{i=1}^{\infty} \lambda_i^{-s} f_i \scal{\psi_i}{v},\quad \forall v \in H^1_0(\Omega).
    \label{eq:def_u}
\end{equation}
If we notice that
\begin{equation}
    u_i = \lambda_i^{-s} f_i,\qquad \forall i\geqslant
    1,
    \label{eq:coef_u_f}
\end{equation}
then, for $f$ in $L^2(\Omega)$ we can show that
\begin{equation}
    \norm{u}_{\mathbb H^{2s}} = \norm{f}_{L^2}.
\end{equation}

\subsection{Rational approximation}\label{subsec:rational_approximation}
Our method relies on rational approximations of the real function $\lambda
\mapsto \lambda^{-s}$ for $\lambda \geqslant
\lambda_0$ where $\lambda_0 > 0$ is a fixed lower bound and for $s$ in $(0, 1)$.
We present here the two examples we are interested in: the BP method introduced in \cite{Bonito2013} and the BURA method presented in \cite{Harizanov2016}.

Both these methods are rational approximation methods; they aim at  approximating the function $\lambda \mapsto \lambda^{-s}$ by a rational function of the form
\begin{equation}
    \mathcal Q_s(\lambda) := C_1(s,N) + C_2(s,N) \sum_{l=-M}^N \frac{a_l(s,N)}{c_l(s,N) + b_l(s,N) \lambda},
    \label{eq:general_rational_scheme}
\end{equation}
where $M, N, C_1, C_2, (a_l)_{l=1}^N, (b_l)_{l=1}^N$ and $(c_l)_{l=1}^N$ are properly chosen parameters.
Many rational approximation schemes have been proposed in the literature, see e.g.\ \cite{Aceto2017,Aceto2018,Gavrilyuk2003,Harizanov2020,Vabishchevich2019}.

Before we present the two particular schemes we used in our numerical experiments, we want to highlight again that the error estimation scheme developed later can be derived in the same manner regardless of the choice of the rational approximation, as long as it leads to a set of well--posed non--fractional parametric problems.

\subsubsection{The BP method}
The BP method is based on the following expression derived from the Balakrishnan's formula \cite{Balakrishnan1959}
\begin{equation}
    \lambda^{-s} = \frac{2\sin(\pi s)}{\pi}
    \int_{-\infty}^{+\infty} \e^{2s y} \left(1 + \e^{2y}
    \lambda\right)^{-1}\ \dy.
    \label{eq:balakrishnan}
\end{equation}
Then, the rational approximation is obtained from \cref{eq:balakrishnan} by discretizing the integral on the right-hand side with a
trapezoidal quadrature rule,
\begin{equation}
    \mathcal Q_s(\lambda) := \frac{2 \kappa \sin(\pi s)}{\pi
    } \sum_{l=-M(\kappa)}^{N(\kappa)} \e^{2sl\kappa} \left(1+ \e^{2l\kappa} \lambda\right)^{-1},
    \label{eq:bp_rational_approximation}
\end{equation}
where $\kappa > 0$ is the fineness parameter and
\begin{equation}
    M(\kappa) := \left \lceil \frac{\pi^2}{4s\kappa^2}\right \rceil,\quad
    \text{and}\quad N(\kappa) := \left \lceil \frac{\pi^2}{4(1-s)\kappa^2} \right
        \rceil,
\end{equation}
where $\lceil \cdot \rceil$ is the ceiling function.
Thus,
\begin{equation}
    \lim_{\kappa \to 0} \mathcal Q_s(\lambda) = \lambda^{-s},\quad \forall (\lambda, s) \in [\lambda_0, +\infty)\times(0,1).
\end{equation}

The convergence of the BP method has been studied in \cite{Bonito2013,Bonito2019}.
Especially, it has been proved that $\mathcal Q_s$ converges uniformly to $\lambda \mapsto \lambda^{-s}$ at an exponential rate as $\kappa$ tends to zero.
The method has been applied to the discretization of various types of PDEs in e.g. \cite{Bonito2013,Bonito2017a,Bonito2016,Bonito2017}.

\subsubsection{The BURA method}

The starting point of the BURA method is the approximation of the function $g(\lambda)=\lambda^s$ for $\lambda \in [0, \lambda_0^{-1}]$.
In this study we used the novel algorithm from \cite{Hofreither2021} to compute the residuals and poles of the BURA of $g(\lambda)$.
From this algorithm we obtain the following rational approximation
\begin{equation}
    \mathcal R_s(\lambda) := \mathcal R_s(0) + \sum_{l=1}^N \frac{r_l}{p_l} + \sum_{l=1}^N \frac{r_l}{\lambda - p_l},
    \label{eq:bura_lambda_s}
\end{equation}
where $N$ is the degree of the rational function, $(r_l)_{l=1}^N \subset \R^{+*}$ are its residuals and $(p_l)_{l=1}^N \subset \R^{-*}$ its poles.
Especially, if $\mathbb P_N$ is the space of polynomial functions of degree $N$ on $(0,\lambda_0^{-1}]$ and
\begin{equation}
    \mathbb Q_N = \left\{ \frac{p}{q},\ p,q \in \mathbb P_N, q\neq 0\right\},
\end{equation}
is the space of rational functions of degree $N$ over $(0,\lambda_0^{-1}]$, then
\begin{equation}
    \mathcal R_s \simeq \argmin_{\tilde r \in \mathbb Q_N} \norm{g-\tilde r}_{L^{\infty}(0,\lambda_0^{-1}]}.
\end{equation}
See \cite{Hofreither2021} for a detailed discussion.

Then, an approximation of $\lambda \mapsto \lambda^{-s}$ can be obtained using partial fraction decomposition
\begin{align}
    \mathcal Q_s(\lambda) := \mathcal R_s(\lambda^{-1}) &= \mathcal R_s(0) + \sum_{l=1}^N \frac{r_l}{p_l} + \sum_{l=1}^N \frac{r_l}{\lambda^{-1}-p_l}\notag\\
    & = \mathcal R_s(0) + \sum_{l=1}^N \frac{r_lp_l + r_l (\lambda^{-1}-p_l)}{p_l(\lambda^{-1} - p_l)}\notag\\
    & = \mathcal R_s(0) + \sum_{l=1}^N \frac{r_l}{p_l(1-p_l \lambda)}.
    \label{eq:bura_rational_approximation}
\end{align}

Like for the BP method, the convergence of the BURA method is of exponential rate, as proved in \cite{Stahl2003}.
However, numerical evidence in the literature show that the BURA method is more efficient than the BP method (see e.g. \cite{Hofreither2020}).
In particular, it leads to a lowest number of parametric solves for a given accuracy.
For more detailed discussions on BURA methods, see e.g. \cite{Aceto2017,Harizanov2020,Hofreither2020,Hofreither2021}.

\section{Discretization}\label{sec:discretization}
We combine the rational approximation \cref{eq:general_rational_scheme} with a
finite element method to derive a fully discrete approximation of the solution
$u$ to \cref{eq:fractional_strong_form}.

\subsection{Rational semi-discrete approximation}

The general scheme used for the semi--discrete approximation will be replaced by the BP method \cref{eq:bp_rational_approximation} and the BURA method \cref{eq:bura_rational_approximation} in the numerical experiments of \cref{sec:numerical_results}.

We define a semi--discrete approximation of the solution $u$, defined in \cref{eq:def_u} by replacing $\lambda_i^{-s}$ by $\mathcal Q_s(\lambda_i)$ in \cref{eq:def_u}.
Thus,
\begin{equation}
    \scal{u_{\mathcal Q_s}}{v} := \sum_{i=1}^{+\infty} \left(C_1 + C_2 \sum_{l=-M}^N \frac{a_l}{c_l + b_l \lambda_i} \right) f_i \scal{\psi_i}{v},\quad \forall v \in H^1_0(\Omega).
    \label{eq:semi-discrete_u}
\end{equation}
If we interchange the two sums in \cref{eq:semi-discrete_u}, we obtain
\begin{equation}
    \scal{u_{\mathcal Q_s}}{v} = C_1 \scal{f}{v} + C_2 \sum_{l=-M}^N a_l \scal{u_l}{v},\quad \forall v \in H^1_0(\Omega),
    \label{eq:rational_discretization}
\end{equation}
where the functions $\{u_l\}_{l=-M}^N$ are solutions to the parametric problems:
for each $l$ in $\llbracket-M,N\rrbracket$, find $u_l$ in $H^1_0$ such that
\begin{equation}
    c_l \scal{u_l}{v} + b_l \scal{\nabla u_l}{\nabla v} =
    \scal{f}{v},\quad \forall v \in H^1_0(\Omega).
    \label{eq:subproblem_u}
\end{equation}
Note, \cref{eq:rational_discretization} can be written as
\begin{equation}
    u_{\mathcal Q_s} = C_1 f + C_2 \sum_{l=-M}^N a_l u_l.
\end{equation}
Since the family of eigenfunctions of the Laplacian $\{\psi_i\}_{i=1}^{+\infty} \subset H^1_0(\Omega)$ is a Hilbert basis of $L^2(\Omega)$, the space $H^1_0(\Omega)$ is dense in $L^2(\Omega)$ and the orthogonal projection of $f$ onto $H^1_0(\Omega)$ in \cref{eq:rational_discretization} is $f$ itself.

Using \cref{eq:rational_discretization}, we reduce the problem \cref{eq:fractional_strong_form} to a series of parametric non--fractional problems \cref{eq:subproblem_u} to solve.
As we will see in \cref{sec:numerical_results}, the number of parametric problems to solve can be relatively small, especially when the BURA method is used in the semi--discrete approximation.
Moreover, these non--fractional problems are fully independent from each other so their solve can be performed in parallel.

\subsection{Finite element discretization}
In order to get a fully discrete approximation of $u$, we use a finite element
method to discretize the parametric problems \cref{eq:subproblem_u}.
Although it is not mandatory, we use the same mesh and same finite element space
for all the parametric problems. We discuss this choice, and possible alternative
strategies, in \cref{subsec:fe_adaptive_refinement}.

Let $\T$ be a mesh on the domain $\Omega$, composed of cells $\T =
\left\lbrace T \right\rbrace$, facets $\E = \left\lbrace E \right\rbrace$
(we call \textit{facets} the edges in dimension two and the faces in
dimension three), and vertices.
The mesh $\T$ is supposed to be regular, in Ciarlet's sense: $h_T/\rho_T
\leqslant \gamma,\ \forall T \in \T$, where $h_T$ is the diameter of a cell
$T$, $\rho_T$ the diameter of its inscribed ball, and $\gamma$ is a positive
constant fixed once and for all.
The subset of facets that are not coincident with the boundary $\Gamma$
(called interior facets) is denoted $\E_I$.
Let $n^+$ and $n^-$ in $\R^d$ be the outward unit normals to a given edge as
seen by two cells $T^+$ and $T^-$ incident to a common edge $E$.
The space of polynomials of order $p$ on a cell $T$ is denoted $\Pcal_p(T)$
and the continuous Lagrange finite element space of order $p$ on the mesh
$\T$ is defined by
\begin{equation}
    V^p := \left\lbrace v_p \in H^1(\Omega), v_{p|T} \in \Pcal_p(T) \;
    \forall T \in \T \right\rbrace.
    \label{eq:FE_space}
\end{equation}
We denote $V^p_0$ the finite element space composed by functions of $V^p$
vanishing on the boundary $\Gamma$.
For a given index $l$, the finite element discretization of
\cref{eq:subproblem_u} reads: for each $l$ in $\llbracket -M, N\rrbracket$, find
$u_{l,p}$ in $V^p_0$ such that
\begin{equation}
    c_l \scal{u_{l,p}}{v_p} + b_l \scal{\nabla u_{l,p}}{\nabla
    v_p} = \scal{f}{v_p},\quad \forall v_p \in V^p_0.
    \label{eq:FE_subproblem}
\end{equation}
Then, combining the solutions to \cref{eq:FE_subproblem} as in \cref{eq:semi-discrete_u}
we can give a fully discrete approximation of the solution to
\cref{eq:fractional_strong_form} $u_{\mathcal Q_s,p} \in V^p_0$, such that
\begin{equation}
    \scal{u_{\mathcal Q_s,p}}{v_p} = C_1 \scal{f}{v_p} + C_2 \sum_{l=-M}^N a_l \scal{u_{l,p}}{v_p},\quad \forall v_p \in V^p_0.
    \label{eq:full_discretization}
\end{equation}
Note, \cref{eq:full_discretization} can be rewritten as
\begin{equation}
    u_{\mathcal Q_s, p} = C_1 f_{V^p_0} + C_2 \sum_{l=-M}^N a_l u_{l,p},
    \label{eq:full_discretization_strong}
\end{equation}
where $f_{V^p_0}$ is the $L^2$ projection of $f$ onto the finite element space $V^p_0$.
The computation of $u_{\mathcal Q_s,p}$ is summarized in the top part of
\cref{fig:diagram_discrete_solution}.
For a detailed discussion on the derivation of $u_{\mathcal Q_s, p}$ in terms of matrices, see \cite{Hofreither2020}.

\section{Error analysis}

This section aims at studying the total discretization error defined by
\begin{equation}
    e := \norm{u - u_{\mathcal Q_s,p}}_{L^2}.
    \label{eq:total_discretization_error}
\end{equation}
Since for any $s \in (0,1)$, the discrepancy $u - u_{\mathcal Q_s,p}$ belongs to
$\mathbb H^{2s}(\Omega) \subset L^2(\Omega)$, the error can be measured in the
$L^2$ norm for any value of the fractional power $s$.

Using the triangle inequality, $e$ is controlled from above by the sum of the rational discretization error $\norm{u - u_{\mathcal Q_s}}_{L^2}$  and the finite element discretization error $\norm{u_{\mathcal Q_s} - u_{\mathcal Q_s,p}}_{L^2}$ as follow
\begin{equation}
    e = \norm{u-u_{\mathcal Q_s, p}}_{L^2} = \norm{u - u_{\mathcal Q_s} + u_{\mathcal Q_s} - u_{\mathcal Q_s, p}}_{L^2} \leqslant \norm{u - u_{\mathcal Q_s}}_{L^2} + \norm{u_{\mathcal Q_s} - u_{\mathcal Q_s, p}}_{L^2}
    \label{eq:upper_triangle}
\end{equation}
In the following we describe estimators for each contribution however the main novelty of this study comes from the finite element discretization error estimation.

\subsection{Rational approximation error analysis}
\label{subsec:rational_estimation}
The estimation of the error induced by the approximation of $u$ by $u_{\mathcal Q_s}$ reduces to the estimation of the scalar rational approximation error.
First, we notice that, if we take $v = \psi_i$ in \cref{eq:semi-discrete_u}, the coefficients $u_{\mathcal Q_s,i}$ of $u_{\mathcal Q_s}$ in the basis $\{\psi_i\}_{i=1}^{+\infty}$ are given by
\begin{equation}
    u_{\mathcal Q_s,i} = \mathcal Q_s(\lambda_i) f_i,\qquad \forall i\geqslant 1.
    \label{eq:coefficient_semi_discrete_approximation}
\end{equation}

Then, using the expansion in the basis $\{\psi_i\}_{i=1}^{+\infty}$, \cref{eq:def_u} and \cref{eq:coefficient_semi_discrete_approximation} we have
\begin{align*}
    \norm{u - u_{\mathcal Q_s}}_{L^2}^2 &= \sum_{i=1}^{+\infty} (u_i - u_{\mathcal Q_s, i})^2\\
    &= \sum_{i=1}^{+\infty} \left(\lambda_i^{-s} - \mathcal Q_s(\lambda_i)\right)^2 f_i^2\\
    &\leqslant \max_{\lambda \geqslant \lambda_0}(\lambda^{-s} - \mathcal Q_s(\lambda))^2 \norm{f}_{L^2}^2.
\end{align*}
So,
\begin{equation}
    \norm{u - u_{\mathcal Q_s}}_{L^2} \leqslant \max_{\lambda \geqslant \lambda_0}(|\lambda^{-s} - \mathcal Q_s(\lambda)|) \norm{f}_{L^2}.
    \label{eq:upper_bound_rational_error}
\end{equation}
Thus, the following rational approximation error estimator
\begin{equation}
    \eta_{\mathcal Q_s} = \max_{\lambda_0 \leqslant \lambda \leqslant \lambda^+}(|\lambda^{-s} - \mathcal Q_s(\lambda)|) \norm{f}_{L^2},
    \label{eq:rational_error_estimator}
\end{equation}
for a large value of $\lambda^+$, is cheap to compute since it consists in the approximation of the maximum of a scalar function and the approximation of the $L^2$ norm of the data $f$.
In addition, we emphasize that this estimator does not require the discrete solution $u_{\mathcal Q_s, p}$ unlike the finite element error estimator we describe in the next section.
However, $\eta_{\mathcal Q_s}$ depends on the lower bound of the Laplacian spectrum $\lambda_0$ and thus can be optimized by taking $\lambda_0 = \lambda_1$ when $\lambda_1$ is known.
When it is not the case, precise guaranteed lower bounds for $\lambda_1$ could be obtained
following e.g.\ \cite{Cances2017,Carstensen2014}.

\subsection{Finite element discretization error analysis}
Our goal is to derive a quantity $\eta$, depending on $f$ and the finite element approximations $(u_{l,p})_{l=-M}^N$ such that
\begin{equation}
    \eta \simeq \norm{u_{\mathcal Q_s} - u_{\mathcal Q_s, p}}_{L^2}.
\end{equation}
Our method is based on the Bank--Weiser finite element error estimator introduced in \cite{Bank1985} and its implementation in the FEniCSx software described in \cite{Bulle2021a}.

\subsubsection{Heuristics}\label{subsubsec:heuristics}
Let us start with some heuristics motivating the derivation of our a posteriori error
estimator.
The main idea is to derive a function $e^{\bw}_{\mathcal Q_s,T}$ that locally
represents the discretization error in the solution to the fractional problem
$(u_{\mathcal Q_s} - u_{\mathcal Q_s,p})_{|T}$ on a cell $T$ of the mesh.
Combining \cref{eq:rational_discretization} and \cref{eq:full_discretization_strong}, on a cell $T$ of the mesh we have
\begin{equation}
    (u_{\mathcal Q_s} - u_{\mathcal Q_s,p})_{|T} = C_1 (f-f_{V^p_0})_{|T} + C_2 \sum_{l=-M}^N a_l (u_l - u_{l,p})_{|T}.
\end{equation}
We approximate the difference $f-f_{V^p_0}$ by $f_{\overline V} - f_{V^p_0}$ where $f_{\overline V}$ is the $L^2$ projection of $f$ onto a finer finite element space $\overline V$ (e.g. we can choose $\overline V := V^{p+1}$).
Note that this step is not necessary in the BP method since the projection of $f$ onto $V^p_0$ is not involved in the rational sum.

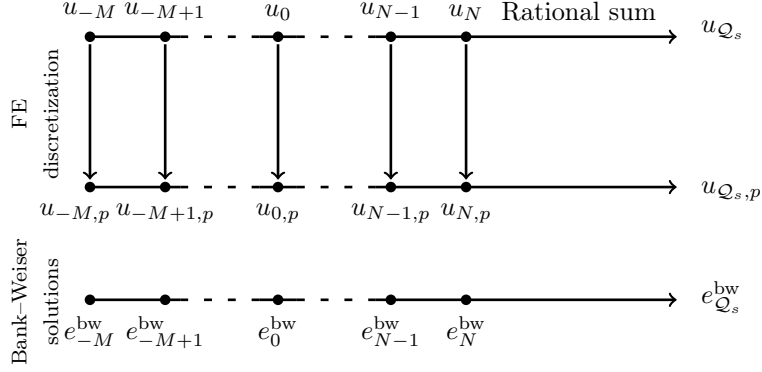
\begin{figure}
    \begin{center}
\begin{tikzpicture}
    \draw[thick, line width = 1pt, -] (-2.5, 0) -- (-1.3, 0);
    \draw[loosely dashed, line width = 1pt, -] (-1.3, 0) -- (-0.25, 0);
    \draw[thick, line width = 1pt, -] (-0.25, 0) -- (0.25, 0);
    \draw[loosely dashed, line width = 1pt, -] (0.25, 0) -- (1.3, 0);
    \draw[thick, line width = 1pt, ->] (1.3, 0) -- (5.3, 0);

    \coordinate(1) at (-2.5, 0);
        \fill[black] (1) circle (2pt);
        \node[above = 1mm of 1] {$u_{-M}$};
    \coordinate(2) at (-1.5, 0);
        \fill[black] (2) circle (2pt);
        \node[above = 1mm of 2] {$u_{-M+1}$};
    \coordinate(3) at (0, 0);
        \fill[black] (3) circle (2pt);
        \node[above = 1mm of 3] {$u_0$};
    \coordinate(4) at (1.5, 0);
        \fill[black] (4) circle (2pt);
        \node[above = 1mm of 4] {$u_{N-1}$};
    \coordinate(5) at (2.5, 0);
        \fill[black] (5) circle (2pt);
        \node[above = 1mm of 5] {$u_N$};
    \coordinate(6) at (4, 0);
        \node[above = 1mm of 6, text width=3cm, align=center] {Rational sum}; 
        \node[above right = -1.5mm and 15mm of 6] {$u_{\mathcal Q_s}$};

    \draw[thick, line width = 1pt, -] (-2.5, -2) -- (-1.3, -2);
    \draw[loosely dashed, line width = 1pt, -] (-1.3, -2) -- (-0.25, -2);
    \draw[thick, line width = 1pt, -] (-0.25, -2) -- (0.25, -2);
    \draw[loosely dashed, line width = 1pt, -] (0.25, -2) -- (1.3, -2);
    \draw[thick, line width = 1pt, ->] (1.3, -2) -- (5.3, -2);

    \coordinate(7) at (-2.5, -2);
        \fill[black] (7) circle (2pt);
        \node[below left = 1mm and -4mm of 7] {$u_{-M,p}$};
    \coordinate(8) at (-1.5, -2);
        \fill[black] (8) circle (2pt);
        \node[below = 1mm of 8] {$u_{-M+1,p}$};
    \coordinate(9) at (0, -2);
        \fill[black] (9) circle (2pt);
        \node[below = 1mm of 9] {$u_{0,p}$};
    \coordinate(10) at (1.5, -2);
        \fill[black] (10) circle (2pt);
        \node[below = 1mm of 10] {$u_{N-1,p}$};
    \coordinate(11) at (2.5, -2);
        \fill[black] (11) circle (2pt);
        \node[below = 1mm of 11] {$u_{N, p}$};

    \coordinate(13) at (4, -2);
    \node[above right = -3mm and 15mm of 13] {$u_{\mathcal Q_s,p}$};

    \draw[thick, line width=1pt, ->] (-2.5, -0.1) -- (-2.5, -1.9);
    \draw[thick, line width=1pt, ->] (-1.5, -0.1) -- (-1.5, -1.9);
    \draw[thick, line width=1pt, ->] (0, -0.1) -- (0, -1.9);
    \draw[thick, line width=1pt, ->] (1.5, -0.1) -- (1.5, -1.9);
    \draw[thick, line width=1pt, ->] (2.5, -0.1) -- (2.5, -1.9);

    \coordinate(12) at (-2.5, -1);
        \node[above left = 11mm and 3mm of 12, align=center, text width=2cm, rotate=90] {\footnotesize FE\\ discretization};
    
    \draw[thick, line width = 1pt, -] (-2.5, -3.5) -- (-1.3, -3.5);
    \draw[loosely dashed, line width = 1pt, -] (-1.3, -3.5) -- (-0.25, -3.5);
    \draw[thick, line width = 1pt, -] (-0.25, -3.5) -- (0.25, -3.5);
    \draw[loosely dashed, line width = 1pt, -] (0.25, -3.5) -- (1.3, -3.5);
    \draw[thick, line width = 1pt, ->] (1.3, -3.5) -- (5.3, -3.5);
    
    \coordinate(14) at (-2.5, -3.5);
        \fill[black] (14) circle (2pt);
        \node[below = 1mm of 14] {$e^{\bw}_{-M}$};
    \coordinate(15) at (-1.5, -3.5);
        \fill[black] (15) circle (2pt);
        \node[below = 1mm of 15] {$e^{\bw}_{-M+1}$};
    \coordinate(16) at (0, -3.5);
        \fill[black] (16) circle (2pt);
        \node[below = 1mm of 16] {$e^{\bw}_0$};
    \coordinate(17) at (1.5, -3.5);
        \fill[black] (17) circle (2pt);
        \node[below = 1mm of 17] {$e^{\bw}_{N-1}$};
    \coordinate(18) at (2.5, -3.5);
        \fill[black] (18) circle (2pt);
        \node[below = 1mm of 18] {$e^{\bw}_N$};

    \coordinate(19) at (4, -3.5);
        \node[above right = -3mm and 15mm of 19] {$e^{\bw}_{\mathcal Q_s}$};

    \coordinate(12) at (-2.5, -3.5);
        \node[above left = 11mm and 3mm of 12, align=center, text width=2cm,
        rotate=90] {\footnotesize Bank--Weiser\\ solutions};
\end{tikzpicture}
\end{center}
    \caption{Summary of the computation of the fractional solution approximation
    and of the fractional Bank--Weiser solution.}
    \label{fig:diagram_discrete_solution}
\end{figure}

To approximate the differences $(u_l - u_{l,p})_{l=-M}^N$ we can use the framework proposed by Bank and Weiser in \cite{Bank1985} to
derive solutions $e^{\bw}_{l,T}$ such that
\begin{equation}
    e^{\bw}_{l,T} \simeq (u_l - u_{l,p})_{|T},\quad \forall l \in \llbracket -M,
    N \rrbracket,\ \forall T \in \T.
\end{equation}
We obtain $e^{\bw}_{\mathcal Q_s,T}$ using the rational approximation sum
\begin{equation}
    e^{\bw}_{\mathcal Q_s,T} := C_1 (f_{\overline V} - f_{V^p_0})_{|T} + C_2 \sum_{l=-M}^N a_l e^{\bw}_{l,T} \simeq (u_{\mathcal Q_s} - u_{\mathcal Q_s,p})_{|T},\quad
    \forall T \in \T.
\end{equation}
Finally, we can estimate the $L^2$ error on the cell $T$ by taking the norm of the function
$e^{\bw}_{\mathcal Q_s,T}$
\begin{equation}
    \norm{e^{\bw}_{\mathcal Q_s,T}}_{L^2} \simeq \norm{(u_{\mathcal Q_s} -
    u_{\mathcal Q_s,p})_{|T}}_{L^2}.
\end{equation}
The heuristic of the approximation of the parametric errors $(u_l - u_{l,p})_{l=-M}^N$ is summarized in \cref{fig:diagram_discrete_solution}.

We would like to emphasize that the Bank--Weiser estimator is not the only
possible choice.
In fact, the Bank--Weiser estimator could be replaced with another estimator
based on the solves of local problems, such as e.g.\ the one used in
\cite{Nochetto2015}.

\subsubsection{A posteriori error estimation}\label{subsubsec:FE_estimation}
Let us now derive our a posteriori error estimation method more precisely.
As mentioned in the last subsection, this estimator is based on a hierarchical
estimator computed from the solves of local Neumann problems on the cells and
introduced for the first time by Bank and Weiser in \cite{Bank1985}.

Let $T$ be a cell of the mesh.
We make use of the following local finite element spaces
\begin{equation}
    V^p_T := \left \{ v_{p,T} \in \mathcal P_p(T),\ v_{p,T} = 0\ \text{in}\
        (\Omega \setminus \overline T) \cup (\overline T \cap \partial
    \Omega)\right \}.
\end{equation}
Let us now consider two non-negative integers $p_+$ and $p_-$ such that $p_+
> p_- \geqslant 0$ and $\mathcal L_T: V^{p_+}_T \longrightarrow V^{p_-}_T$
the local Lagrange interpolation operator.
We introduce the \textit{local Bank--Weiser space}, defined by
\begin{equation}
    V^{\bw}_T := \ker(\mathcal L_T) = \left \{ v_{p_+,T} \in V^{p_+}_T,\
    \mathcal L_T(v_{p_+,T}) =0 \right \}.
    \label{eq:local_bw_space}
\end{equation}
The local parametric Bank--Weiser problem associated to the parametric problems
\cref{eq:subproblem_u} and \cref{eq:FE_subproblem} reads
\begin{equation}
    \int_T e^{\bw}_{l,T}v^{\bw}_T + \e^{2l\kappa} \int_T \nabla e^{\bw}_{l,T} \cdot \nabla
    v^{\bw}_T = \int_T r_{l,T} v^{\bw}_T - \frac{1}{2} \sum_{E \in
    \partial T} \int_E J_{l,E} v^{\bw}_T,\ \forall v^{\bw}_T \in V^{\bw}_T
\label{eq:local_bw_equation}
\end{equation}
where $r_{l,T}$ and $J_{l,T}$ are defined as follow:
\begin{equation}
    r_{l,T} := f_{|T} - {u_{l,p}}_{|T} + \e^{2l\kappa} \Delta {u_{l,p}}_{|T},\quad
    \text{and}\quad J_{l,T} := \e^{2l\kappa} \left(\frac{\partial
    u_{l,p|T_+}}{\partial n} - \frac{\partial u_{l,p|T_-}}{\partial n} \right),
\end{equation}
where $T_+$ and $T_-$ are the cells sharing the edge $E$ such that the normal $n$ is outward $T_+$.
The solution $e^{\bw}_{l,T}$ in $V^{\bw}_T$ is the local parametric
Bank--Weiser solution.
More details about the computation and implementation of the Bank--Weiser
solutions can be found in\ \cite{Bank1985,Bulle2021a}.

Then, we derive the local \textit{fractional} Bank--Weiser solution by summing
the local parametric Bank--Weiser solutions into the rational approximation sum
\begin{equation}
    e_{\mathcal Q_s,T}^{\bw} := C_1 (f_{\overline V}-f_{V^p_0})_{|T} + C_2 \sum_{l=-M}^N a_l e_{l,T}^{\bw}.
    \label{eq:local_fractional_bank-weiser_solution}
\end{equation}
The local fractional Bank--Weiser estimator is then defined as the $L^2$ norm of
this local solution
\begin{equation}
    \eta^{\bw}_{\mathcal Q_s,T} := \norm{e_{\mathcal Q_s,T}^{\bw}}_{L^2(T)}.
    \label{eq:local_L2_bank-weiser}
\end{equation}
The global fractional Bank--Weiser estimator is then defined by
\begin{equation}
    {\eta^{\bw}_{\mathcal Q_s}}^2 := \sum_{T \in \T} {\eta^{\bw}_{\mathcal Q_s,T}}^2.
    \label{eq:global_L2_bank-weiser}
\end{equation}

\section{Adaptive refinement}

One of the main applications of a posteriori error estimation is to drive adaptive mesh refinement algorithms.
When the error is unevenly spread across the mesh, refining uniformly is a
waste of computational resources leading to suboptimal convergence rates in the number of degrees of freedom.
This problem is compounded for computationally expensive problems like fractional problems.
Moreover, it is known that fractional problems often show a boundary layer behavior, the discretization error is consequently large in a localized region near the boundary \cite{Acosta2017a,Borthagaray2020,Vabishchevich2019}.
This problem has been tackled using graded meshes that are refined near the
boundary based on a priori or a posteriori considerations
\cite{Bonito2018a,Chen2015,Gimperlein2019,Meidner2018}.
As expected, the use of graded meshes improves the convergence of the methods.

Adaptive refinement algorithms are based on the loop
\begin{equation*}
    \cdots \longrightarrow \mathrm{Solve} \longrightarrow \mathrm{Estimate}
    \longrightarrow \mathrm{Mark} \longrightarrow \mathrm{Refine}
    \longrightarrow \cdots
\end{equation*}
In this work we are concerned with developments in the modules \textit{solve}, \textit{estimate} and, when an adaptive rational scheme is used, the module \textit{refine}.

In \cref{subsec:fe_adaptive_refinement}, we focus on the finite element mesh adaptive refinement, choosing a fixed rational scheme fine enough for the rational approximation error to be negligible.
We are using the Dörfler algorithm \cite{dorfler_convergent_1996} for the \textit{mark} module and the Plaza--Carey algorithm\ \cite{plaza_local_2000} for the \textit{refine} module.

In \cref{subsec:fe_ra_adaptive_refinement}, we allow the rational scheme to vary from one refinement step to another, in order to balance the discretization errors.
Thus, the \textit{refine} module is composed of the Plaza--Carey algorithm for the finite element mesh and an algorithm in charge of picking the right rational scheme in order for the rational and finite element approximation errors to be balanced at each refinement step.

Rational approximation methods have the advantage of being fully parallelizable due to the independence of the parametric problems from each other.
Similarly, the local a posteriori error estimation method we have presented
earlier is also parallelizable since the computation of the local Bank--Weiser solutions on the cells are independent from each other.
Our error estimation strategy combines these advantages and is fully
parallelizable both with respect to the parametric problems and local estimators computation.

\subsection{Finite element mesh adaptive refinement}\label{subsec:fe_adaptive_refinement}

An example of error estimation and mesh adaptive refinement algorithm based on our method is shown in \cref{fig:fe_refinement_algorithm}.
In this algorithm, we focus on the finite element error approximation.
Thus, the rational scheme is fixed and chosen so that the rational approximation error is negligible.
In this context, we assume that the total discretization error satisfies
\begin{equation}
    \norm{u-u_{\mathcal Q_s,p}}_{L^2} \simeq \norm{u_{\mathcal Q_s} - u_{\mathcal Q_s,p}}_{L^2}.
    \label{eq:negligible_rational_error}
\end{equation}

\begin{figure}
\begin{algorithmic}
    \State Choose a tolerance $\epsilon >0$, an initial mesh $\T_{n=0}$ and an initial rational scheme $\mathcal Q_{s, n=0}$
    \State Generate $\mathcal Q_{s, n=0}$ coefficients
    \State Initialize the total estimator $\eta^{\bw}_{\mathcal Q_s} = \epsilon +1$
    \While{$\eta^{\bw}_{\mathcal Q_s} > \epsilon$}
        \State Initialize the local Bank--Weiser solutions $\{e_{\mathcal Q_s, T}^{\mathrm{bw}}\}_T$ and the solution $u_{\mathcal Q_s,p}$ to zero
        \For{$l \in \llbracket-M,N\rrbracket$}
            \State Solve \cref{eq:FE_subproblem} on $\mathcal T_n$ to obtain $u_{l,p}$
            \State Add $a_l u_{l,p}$ to $u_{\mathcal Q_s,p}$
            \For{$T \in \mathcal T_n$}
                \State Solve \cref{eq:local_bw_equation} to obtain $e_{l,T}^{\mathrm{bw}}$
                \State Add $a_l e_{l,T}^{\mathrm{bw}}$ to $e_{\mathcal Q_s,T}^{\mathrm{bw}}$
            \EndFor
        \EndFor
        \State Multiply $u_{\mathcal Q_s, p}$ and $e_{\mathcal Q_s, T}^{\mathrm{bw}}$ by $C_2$
        \State Compute $f_{V^p}$ the $L^2$ projection of $f$ onto $V^p$ and add $C_1 f_{V^p}$ to $u_{\mathcal Q_s,p}$
        \State Compute $f_{V^{p+1}}$ the $L^2$ projection of $f$ onto $V^{p+1}$ and add $C_1 (f_{V^{p+1}} - f_{V^p})_{|T}$ to $e_{\mathcal Q_s,T}^{\mathrm{bw}}$
        \State Compute $\eta_{\mathcal Q_s, T}^{\mathrm{bw}} := \norm{e_{\mathcal Q_s, T}^{\mathrm{bw}}}_{L^2(T)}$ for all $T \in \mathcal T_n$ and $\eta_{\mathcal Q_s}^{\mathrm{bw}} :=\sqrt{\sum_T {\eta_{\mathcal Q_s,T}^{\mathrm{bw}}}^2}$
        \If{$\eta^{\bw}_{\mathcal Q_s} > \epsilon$}
            \State Return $u_{\mathcal Q_s,p}$
        \Else
            \State Mark the mesh $\mathcal T_n$ using $\{\eta_{\mathcal Q_s,T}^{\mathrm{bw}}\}_T$
            \State Refine the mesh $\mathcal T_n$ to obtain $\mathcal T_{n+1}$
        \EndIf
    \EndWhile
\end{algorithmic}
\caption{Finite element error estimation and mesh adaptive refinement algorithm outline in pseudo--code.}
\label{fig:fe_refinement_algorithm}
\end{figure}

The algorithm presented in \cref{fig:fe_refinement_algorithm} is based on three loops: one \textbf{While} loop and two \textbf{For} loops.
The \textbf{While} loop is due to the mesh adaptive refinement procedure and can not be parallelized.
However, the two \textbf{For} loops are fully parallelizable and this parallelization can be highly advantageous for large three-dimensional problems.

Note that there is no guarantee that the mesh we obtain at the end of the main
\textbf{While} loop in \cref{fig:fe_refinement_algorithm} is optimal for all the parametric
problems.
For some of the parametric solutions without boundary layers the mesh is
certainly over-refined.
An alternative approach could be to compute the $L^2$ norms of the parametric
Bank--Weiser solutions $e^{\bw}_{l,T}$ in order to derive parametric
Bank--Weiser estimators and refine the meshes independently for each parametric
problem.
This would require the storage of a possibly different mesh for each parametric
problem at each iteration.
More importantly, this would mean summing parametric finite element solutions
coming from different and possibly non-nested meshes.
Properly addressing this question is beyond the scope of this study.
Nonetheless, we give some hints in the numerical results \cref{subsec:2d-sines}.

\subsection{Rational scheme and finite element mesh adaptive refinement}\label{subsec:fe_ra_adaptive_refinement}

In this section, we introduce a method to adaptively refine the rational scheme in addition to the finite element mesh.
This method is partly inspired from the Continuation Multilevel Monte Carlo method, applied to stochastic PDEs, introduced in \cite{Collier2015}.
As in \cref{eq:upper_triangle}, we can use the second triangle inequality to obtain a lower bound on the total error
\begin{equation}
    \norm{u - u_{\mathcal Q_s,p}}_{L^2} \geqslant \big|\norm{u-u_{\mathcal Q_s}}_{L^2} - \norm{u_{\mathcal Q_s} - u_{\mathcal Q_s,p}}_{L^2}\big|.
    \label{eq:lower_triangle}
\end{equation}
Thus, the only way to reduce this lower bound to zero is to balance $\norm{u-u_{\mathcal Q_s}}_{L^2}$ and $\norm{u_{\mathcal Q_s} - u_{\mathcal Q_s, p}}_{L^2}$.
This also make sense from a more practical perspective, using a very fine rational approximation scheme is a waste of computational resources if a too coarse finite element scheme generates a large error, the inverse being also true.

The rational approximation error will be controlled via the rational estimator $\eta_{\mathcal Q_s}$, defined in \cref{eq:rational_error_estimator}, which will be used to choose the proper rational scheme at each refinement step.
According to the results from \cite{Bonito2013} and \cite{Stahl2003} both the BP and BURA rational schemes converge exponentially fast.
On the other hand, the finite element error usually shows a polynomial convergence rate.
Thus, in order to balance the rational approximation error with the finite element error, we must reduce the rational approximation convergence rate to match the finite element one.

At step $n$ we choose the rational approximation scheme of step $n+1$ such that the rational error estimator matches the Bank--Weiser estimator value.
To do so we need to estimate what will be the value of the Bank--Weiser estimator.
We assume that the logarithm of the Bank--Weiser estimator follows a linear trend (this is usually true in the asymptotic regime).
Thus, the next value of the estimator can be estimated by a simple linear regression on its values at step $n$ and $n-1$.

More precisely, let us denote $\rho_{\mathcal Q_s, n}^{\mathrm{bw}} = \ln(\eta_{\mathcal Q_s, n}^{\mathrm{bw}})$ the logarithm of the Bank--Weiser estimator and $d_n$ the logarithm of the number of degrees of freedom at the $n^{\mathrm{th}}$ refinement step.
We assume that
\begin{equation}
    \rho_{\mathcal Q_s, m}^{\mathrm{bw}} = \delta d_m + C,\quad \text{for } m = n+1, n, n-1,
    \label{eq:asymptotic_regime}
\end{equation}
where $\delta$ is the convergence slope and $C$ is a real constant.
Then, we have
\begin{equation}
    \rho_{\mathcal Q_s, n}^{\mathrm{bw}} - \rho_{\mathcal Q_s, n-1}^{\mathrm{bw}} = \delta (d_n - d_{n-1}),\quad \text{and}\quad \rho_{\mathcal Q_s, n+1}^{\mathrm{bw}} - \rho_{\mathcal Q_s, n}^{\mathrm{bw}} = \delta (d_{n+1} - d_n).
\end{equation}
In other words,
\begin{equation}
    \frac{\rho_{\mathcal Q_s, n+1}^{\mathrm{bw}} - \rho_{\mathcal Q_s, n}^{\mathrm{bw}}}{d_{n+1} - d_n} = \frac{\rho_{\mathcal Q_s, n}^{\mathrm{bw}} - \rho_{\mathcal Q_s, n-1}^{\mathrm{bw}}}{d_n - d_{n-1}} = \delta.
\end{equation}
Thus,
\begin{equation}
    \rho_{\mathcal Q_s, n+1}^{\mathrm{bw}} = (\rho_{\mathcal Q_s, n}^{\mathrm{bw}} - \rho_{\mathcal Q_s, n-1}^{\mathrm{bw}}) \frac{d_{n+1} - d_n}{d_n - d_{n-1}} + \rho_{\mathcal Q_s, n}^{\mathrm{bw}}.
\end{equation}
Now, if we denote $\rho_{\mathcal Q_s, n+1} = \ln(\eta_{\mathcal Q_s, n+1})$ the logarithm of the rational error estimator at the $n+1^{\mathrm{th}}$ step, our goal is to compute the coarsest rational scheme such that
\begin{equation}
    \rho_{\mathcal Q_s, n+1} < \rho_{\mathcal Q_s, n+1}^{\mathrm{bw}} = (\rho_{\mathcal Q_s, n}^{\mathrm{bw}} - \rho_{\mathcal Q_s, n-1}^{\mathrm{bw}}) \frac{d_{n+1} - d_n}{d_n - d_{n-1}} + \rho_{\mathcal Q_s, n}^{\mathrm{bw}}.
    \label{eq:upper_bound_ln_rational_estimator}
\end{equation}
In terms of estimators values \cref{eq:upper_bound_ln_rational_estimator} can be reformulated into
\begin{equation}
    \eta_{\mathcal Q_s, n+1} < \tilde \delta \eta_{\mathcal Q_s, n}^{\mathrm{bw}},
    \label{eq:upper_bound_rational_estimator}
\end{equation}
where
\begin{equation}
    \tilde \delta := \exp \left((\rho_{\mathcal Q_s, n}^{\mathrm{bw}} - \rho_{\mathcal Q_s, n-1}^{\mathrm{bw}}) \frac{d_{n+1} - d_n}{d_n - d_{n-1}}\right).
    \label{eq:delta_tilde}
\end{equation}
Notice that at the end of the $n^{\mathrm{th}}$ refinement step, all the quantities involved in the definition \cref{eq:delta_tilde} of $\tilde \delta$ are known.
Especially, $d_{n+1}$ can be computed once the $n^{\mathrm{th}}$ mesh is refined using $\eta_{\mathcal Q_s, n}^{\mathrm{bw}}$.

Since this method requires at least two levels of refinement we choose a coarse rational scheme and keep it for step $n=0$ and step $n=1$.
Despite the fact that \cref{eq:asymptotic_regime} is only true in the asymptotic regime, numerical experiments suggest that our method stabilizes (in the sense that $\eta_{\mathcal Q_s, n} \simeq \eta_{\mathcal Q_s, n}^{\mathrm{bw}}$) after a few refinement steps.

Finally, we stop the algorithm when the total estimator, defined by
\begin{equation}
    \eta = \eta_{\mathcal Q_s} + \eta_{\mathcal Q_s}^{\mathrm{bw}},
    \label{eq:total_estimator}
\end{equation}
reaches the prescribed tolerance (according to the bound \cref{eq:upper_triangle}).
An example of algorithm based on our method is shown in \cref{fig:fe_ra_refinement_algorithm}.
\begin{figure}
    \begin{algorithmic}
        \State Choose a tolerance $\epsilon >0$, an initial mesh $\T_{n=0}$, an initial rational scheme $\mathcal Q_{s, n=0}$
        \State Generate $\mathcal Q_{s, n=0}$ coefficients
        \State Compute $\eta_{\mathcal Q_s}$
        \State Initialize the total estimator $\eta = \epsilon +1$
        \While{$\eta > \epsilon$}
            \State Initialize the local Bank--Weiser solutions $\{e_{\mathcal Q_s, T}^{\mathrm{bw}}\}_T$ to zero
            \State Initialize the solution $u_{\mathcal Q_s,p}$ to zero
            \For{$l \in \llbracket-M,N\rrbracket$}
                \State Solve \cref{eq:FE_subproblem} on $\mathcal T_n$ to obtain $u_{l,p}$
                \State Add $a_l u_{l,p}$ to $u_{\mathcal Q_s,p}$
                \For{$T \in \mathcal T_n$}
                    \State Solve \cref{eq:local_bw_equation} to obtain $e_{l,T}^{\mathrm{bw}}$
                    \State Add $a_l e_{l,T}^{\mathrm{bw}}$ to $e_{\mathcal Q_s,T}^{\mathrm{bw}}$
                \EndFor
            \EndFor
            \State Multiply $u_{\mathcal Q_s, p}$ and $e_{\mathcal Q_s, T}^{\mathrm{bw}}$ by $C_2$
            \State Compute $f_{V^p}$ the $L^2$ projection of $f$ onto $V^p$ and add $C_1 f_{V^p}$ to $u_{\mathcal Q_s,p}$
            \State Compute $f_{V^{p+1}}$ the $L^2$ projection of $f$ onto $V^{p+1}$ and add $C_1 (f_{V^{p+1}} - f_{V^p})_{|T}$ to $e_{\mathcal Q_s,T}^{\mathrm{bw}}$
            \State Compute $\eta_{\mathcal Q_s, T}^{\mathrm{bw}} := \norm{e_{\mathcal Q_s, T}^{\mathrm{bw}}}_{L^2(T)}$ for all $T \in \mathcal T_n$, $\eta_{\mathcal Q_s}^{\mathrm{bw}} :=\sqrt{\sum_T {\eta_{\mathcal Q_s,T}^{\mathrm{bw}}}^2}$ and $\eta := \eta_{\mathcal Q_s}^{\mathrm{bw}} + \eta_{\mathcal Q_s}$
            \If{$\eta > \epsilon$}
                \State Return $u_{\mathcal Q_s,p}$
            \Else
                \State Mark the mesh $\mathcal T_n$ using $\{\eta_{\mathcal Q_s,T}^{\mathrm{bw}}\}_T$
                \State Refine the mesh $\mathcal T_n$ to obtain $\mathcal T_{n+1}$
                \If{$n>1$}
                    \State Compute $\tilde \delta$ using \cref{eq:delta_tilde}
                    \While{$\eta_{\mathcal Q_s, n+1} > \tilde \delta \eta_{\mathcal Q_s,n}^{\mathrm{bw}}$}
                        \State Refine the rational scheme $\mathcal Q_{s,n}$ to obtain $\mathcal Q_{s,n+1}$
                        \State Compute $\eta_{\mathcal Q_s, n+1}$
                    \EndWhile
                \EndIf
            \EndIf
        \EndWhile
    \end{algorithmic}
\caption{Total error estimation and rational scheme and mesh adaptive refinement algorithm outline in pseudo--code.}
\label{fig:fe_ra_refinement_algorithm}
\end{figure}

\section{Implementation}
We have implemented our method using the DOLFINx finite element
solver of the FEniCS Project\ \cite{Alnaes2015}.
Each parametric subproblem is submitted to a batch job queue.
A distinct MPI communicator is used for each job.
We use a standard first-order Lagrange finite element method and the resulting
linear system is solved using the conjugate gradient method.
The conjugate gradient method is preconditioned using
BoomerAMG from HYPRE\ \cite{falgout_hypre_2002} via the interface in PETSc\
\cite{balay_petsc_2016}.
To compute the Bank--Weiser error estimator for each subproblem we use the
methodology outlined in\ \cite{Bulle2021a} and implemented in the FEniCSx--EE
package\ \cite{bulle_fenics-ee_2019}.
For every subproblem the computed solution and error estimate is written to disk
in HDF5 format.
A final step, running on a single MPI communicator, reads the solutions and
error estimates for all subproblems, computes the quadrature sums using
\texttt{axpy} operations, defines the marked set of cells to be refined using
the Dörfler algorithm\ \cite{dorfler_convergent_1996}, and finally refines the
mesh using the Plaza--Carey algorithm\ \cite{plaza_local_2000}.

A more complex implementation using a single MPI communicator split into
sub-communicators would remove the necessity of reading and writing the
solution and error estimate for each subproblem to and from disk.
However, in practice the cost of computing the parametric solutions massively dominates
all other costs.

\section{Numerical results}\label{sec:numerical_results}

In \cref{subsec:2d-sines,subsec:2d-checkerboard,subsec:3d-sines,subsec:3d-checkerboard,subsec:3d-torus} we only consider finite element mesh adaptive refinement.
Thus, we choose the rational scheme $\mathcal Q_s$ in order to guarantee that the
rational approximation error is negligible (i.e. of the order of machine precision).
This can be achieved thanks to \cref{eq:upper_bound_rational_error}.
\Cref{subsec:ra_adaptive_results} is dedicated to the application of the algorithm combining mesh and rational schemes adaptive refinement, sketched in \cref{fig:fe_ra_refinement_algorithm}, to the two dimensional test cases.

\subsection{Two-dimensional product of sines test case}\label{subsec:2d-sines}

We solve \cref{eq:fractional_strong_form} on the square $\Omega = (0,\pi)^2$
with data $f(x,y)=\sin(x)\sin(y)$.
The analytical solution to this problem is given by $u(x,y) = 2^{-s} \sin(x)\sin(y)$.
Moreover, the analytical solutions to the parametric problems
\cref{eq:subproblem_u} are also known  $u_{l}(x,y) = (c_l + 2b_l)^{-1}
\sin(x)\sin(y)$.
The problem is solved on a hierarchy of structured (triangular) meshes.
For this test case the solution $u$ shows no boundary layer behavior, therefore
mesh adaptive refinement cannot improve the convergence rate.
Consequently, we only perform uniform mesh refinement.

The accuracy of the estimator is measured with the efficiency indices (defined as the ratios of the estimators over the exact errors) shown in \cref{tab:2d_sines_efficiency_fixed_ra_schemes}.
As we can see, the efficiency varies between $1$ and $2$ in most of the cases.
In the case of a fixed fine rational scheme the efficiency values are not influenced by the choice of the rational method, except when $s=0.1$.
This independence of the method is remarkable since the parametric problems coefficients can be very different from BP to BURA.

\begin{figure}
    \begin{center}
        \includegraphics[width=0.45\textwidth]{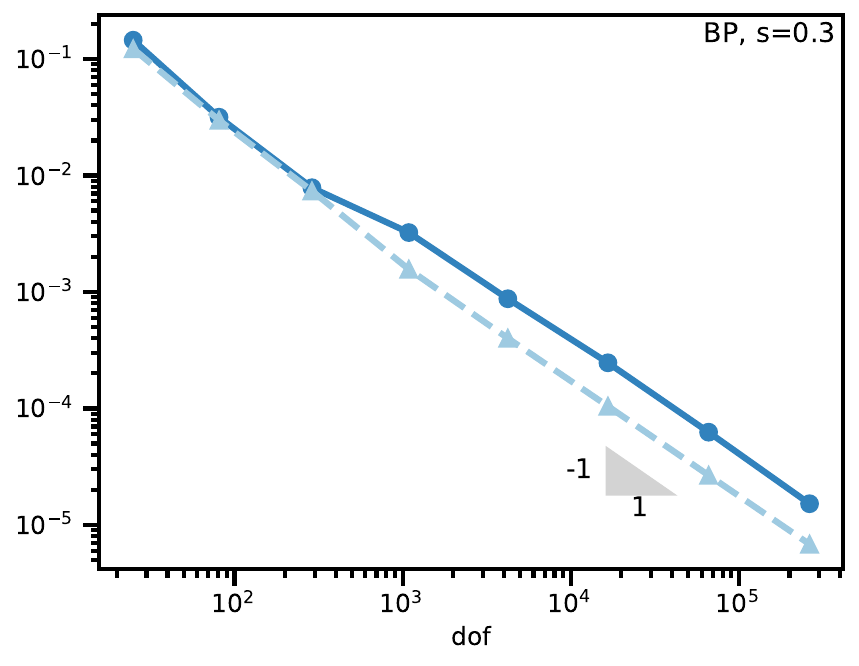}
        \includegraphics[width=0.45\textwidth]{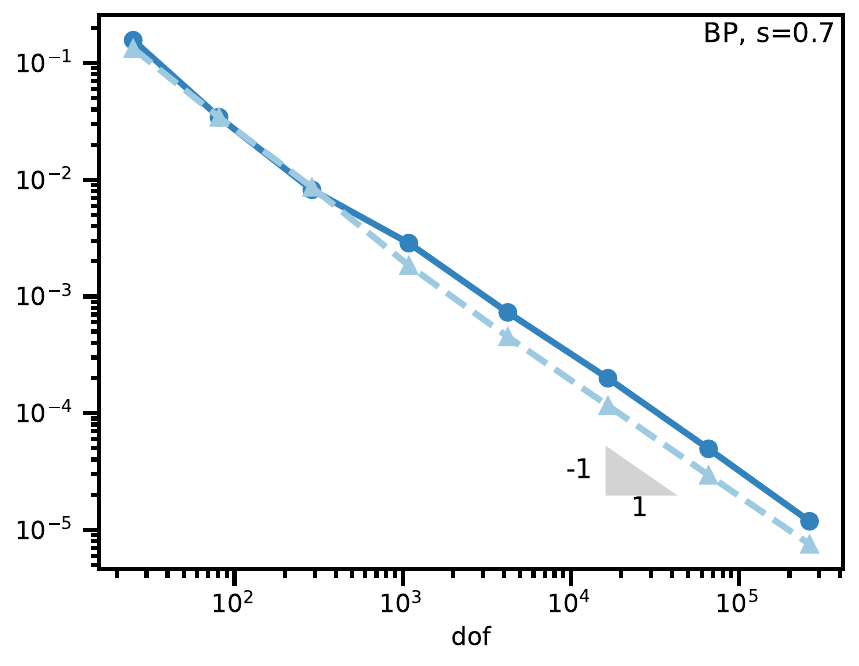}\\
        \includegraphics[width=0.45\textwidth]{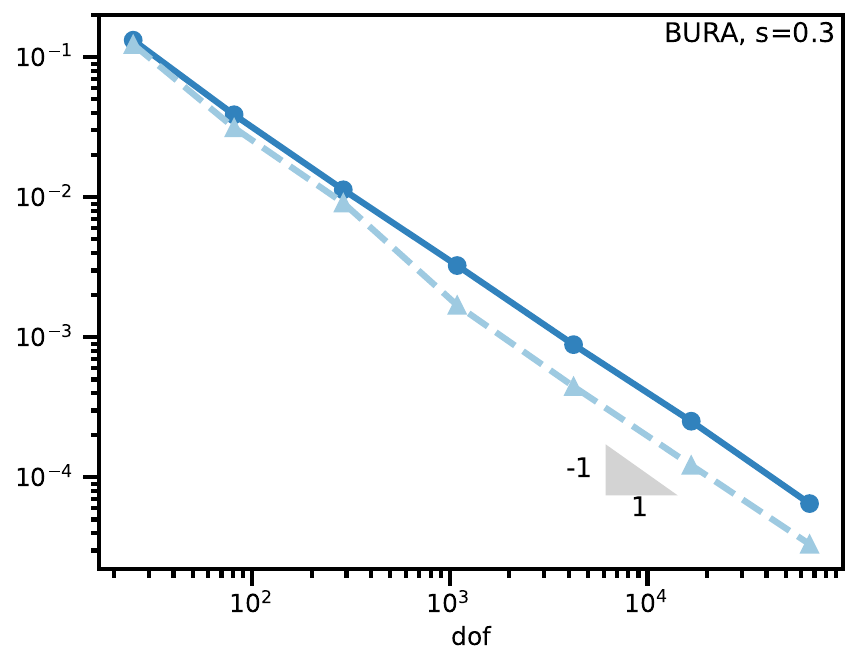}
        \includegraphics[width=0.45\textwidth]{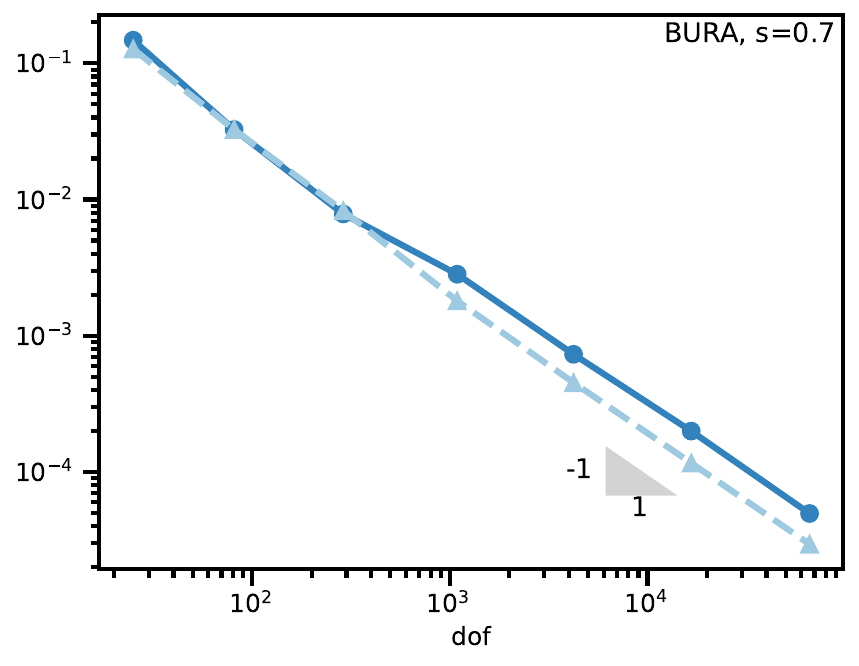}
    \end{center}
    \caption{\textbf{Two-dimensional product of sines test case:} the Bank--Weiser estimator $\eta^{\mathrm{bw}}_{\mathcal Q_s}$ in solid blue line is compared to the exact error in dashed light blue line for two different rational schemes and two different fractional powers.}
    \label{fig:2d_sines_convergence}
\end{figure}

Theorem 4.3 from \cite{Bonito2013} gives a convergence rate for the finite
element scheme discretization error when the BP method is used, depending on the elliptic regularity index $\alpha$ of the Laplacian over $\Omega$, on the fractional power and on the regularity index
$\delta$ of the data $f$.
Since $\Omega$ is convex the elliptic "pick-up" regularity index $\alpha$ can be
taken to be $1$ \cite{Bonito2018a} and since $f$ is infinitely smooth the
coefficient $\delta$ can be taken as large as wanted.
Consequently, Theorem 4.3 in \cite{Bonito2013} predicts a convergence rate of
$\mathrm{dof}^{-1}$ for this test case.
The convergence rates we measure in practice, shown in
\cref{tab:2d_sines_slopes_fixed_ra_schemes}, are mostly coherent with this prediction.
These rates are computed from a linear regression fit on the values obtained on the ten last meshes of the hierarchy (or on all the meshes when less than ten meshes are computed).
To our knowledge, there is no such convergence theorem available for the BURA method.
Convergence plots for the Bank--Weiser estimator and the exact finite element discretization error are shown in \cref{fig:2d_sines_convergence} for $s=0.3$ and $0.7$.

\begin{table}
    \begin{center}
        \begin{tabular}{llrrrrr}
            \toprule
               & Frac. power & 0.1   & 0.3   & 0.5   & 0.7   & 0.9   \\ \midrule
            BP & Estimator   & -0.92 & -0.93 & -0.95 & -0.96 & -0.97 \\
               & Exact error & -1.03 & -1.03 & -1.04 & -1.04 & -1.04 \\ 
               \midrule
            BURA & Estimator   & -0.79 & -0.93 & -0.95 & -0.96 & -0.97 \\
                 & Exact error & -0.83 & -1.04 & -1.05 & -1.06 & -1.05 \\
            \bottomrule
            \end{tabular}
    \end{center}
    \caption{\textbf{Two-dimensional product of sines test case:} convergence rates of
    the Bank--Weiser estimator and of the exact finite element error for various fractional powers and for fixed BP and BURA schemes.}
    \label{tab:2d_sines_slopes_fixed_ra_schemes}
\end{table}

\begin{table}
    \begin{center}
        \begin{tabular}{lrrrrr}
            \toprule
            Frac. power & 0.1  & 0.3  & 0.5  & 0.7  & 0.9  \\ \midrule
            BP          & 1.73 & 2.04 & 1.79 & 1.5  & 1.22 \\
            BURA        & 1.07 & 2.04 & 1.79 & 1.51 & 1.27 \\ 
            \bottomrule
        \end{tabular}
    \end{center}
    \caption{\textbf{Two-dimensional product of sines test case:} efficiency indices of the Bank--Weiser estimator for various fractional powers and for fixed BP and BURA schemes.}
    \label{tab:2d_sines_efficiency_fixed_ra_schemes}
\end{table}

\subsubsection{Parametric problems discretization error}
\label{subsubsec:parametric_problems}
\begin{figure}
    \begin{center}
        \includegraphics[width=0.32\textwidth]{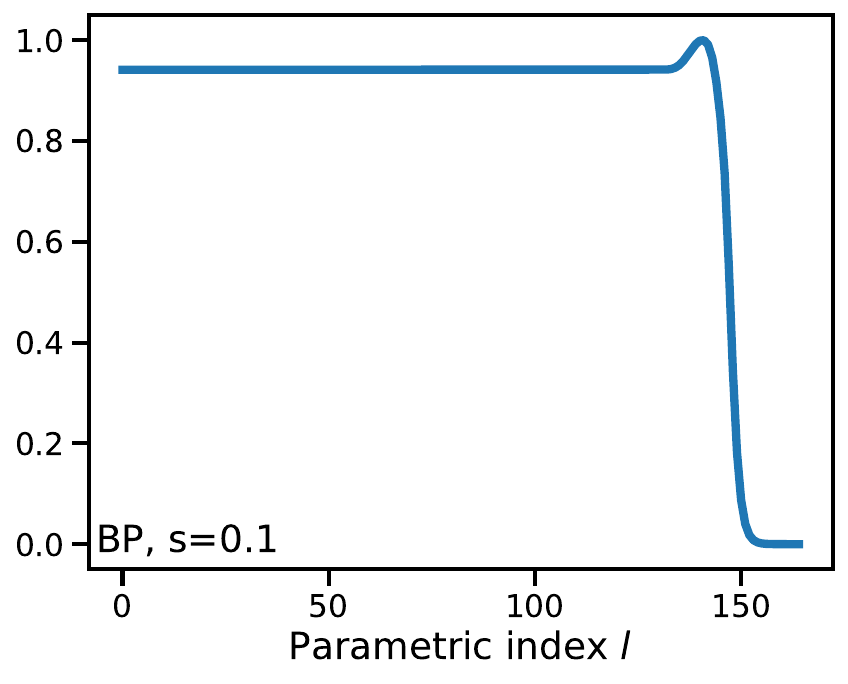}
        \includegraphics[width=0.32\textwidth]{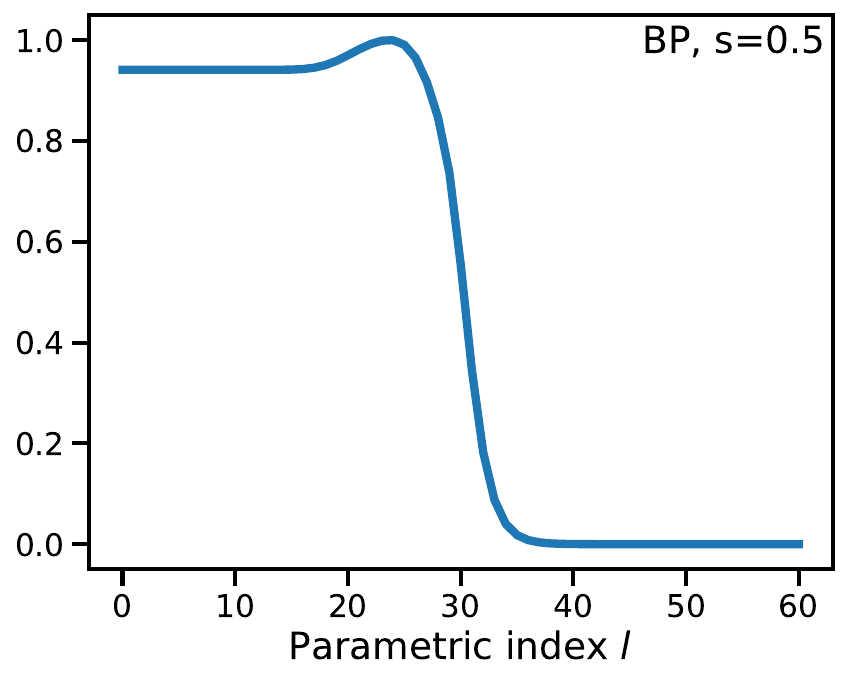}
        \includegraphics[width=0.32\textwidth]{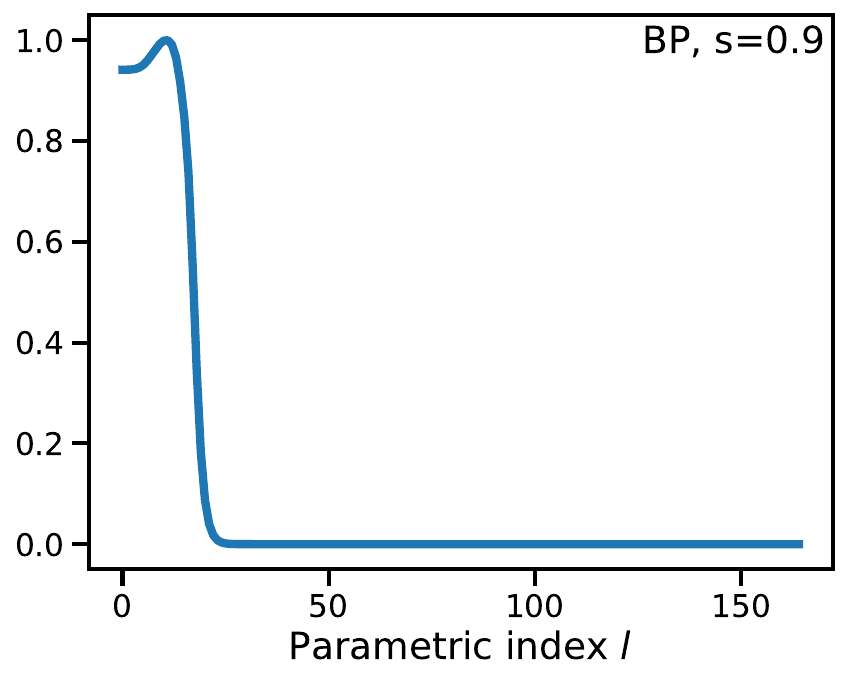}\\
        \includegraphics[width=0.32\textwidth]{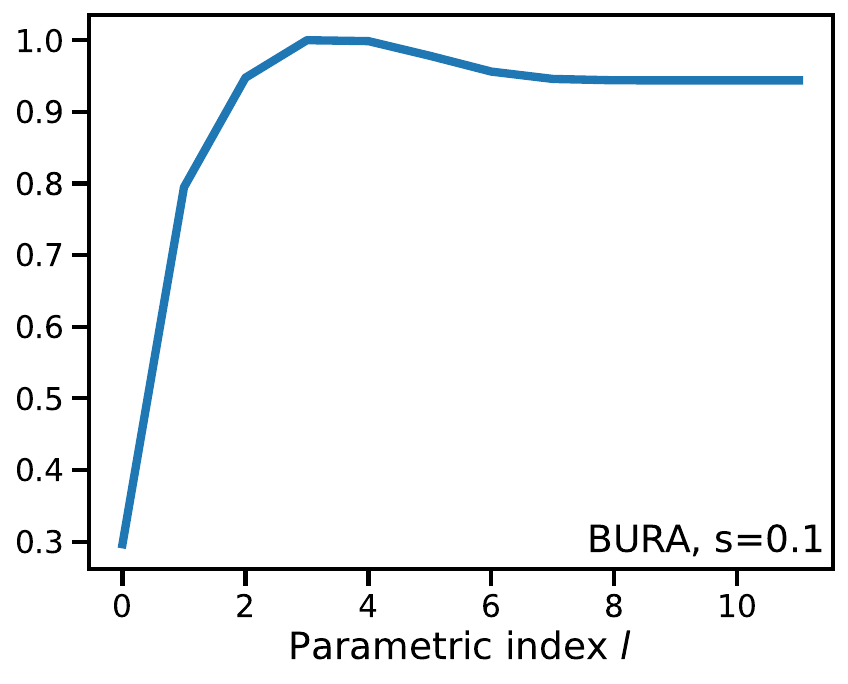}
        \includegraphics[width=0.32\textwidth]{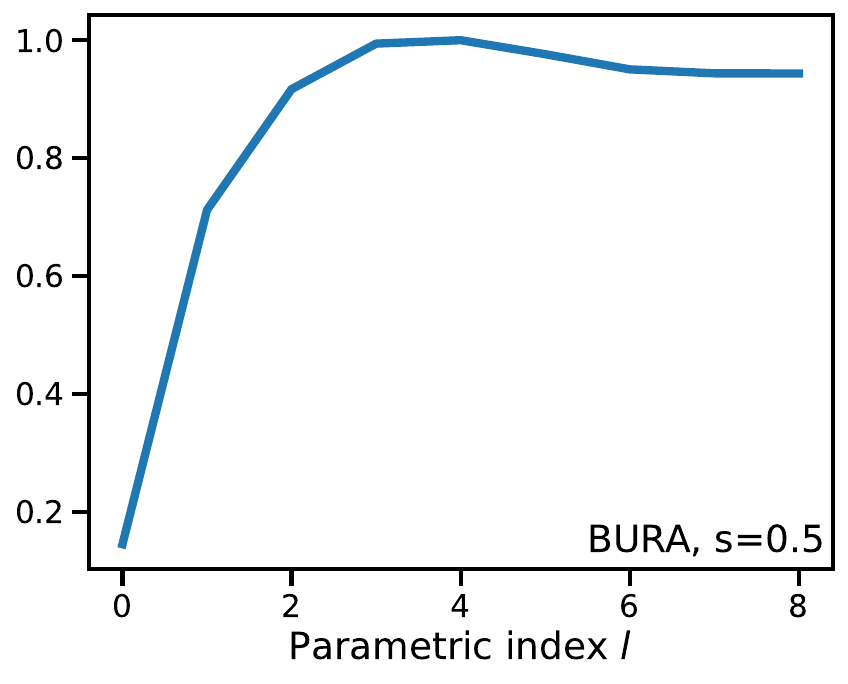}
        \includegraphics[width=0.32\textwidth]{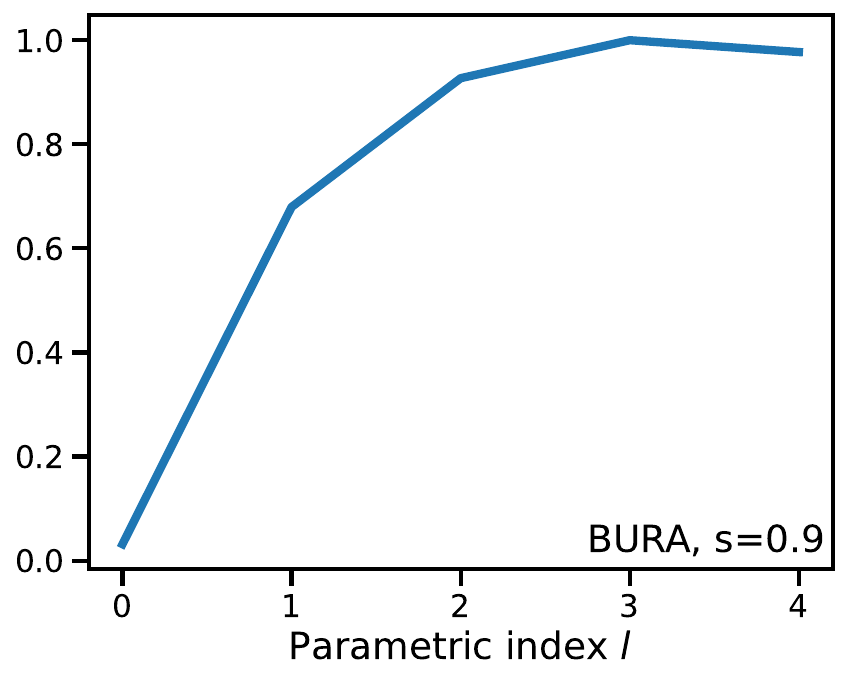}
    \end{center}
    \caption{\textbf{Two-dimensional product of sines text case:} variation of the normalized exact parametric errors $e_l/\max_l(e_l)$ with respect to the index $l \in \llbracket -M, N\rrbracket$ for BP and BURA and for three
    different fractional powers.}
    \label{fig:2d_sines_parametric}
\end{figure}

Since we know the analytical solutions to the parametric problems in this case,
it is possible to compute the exact parametric discretization errors $e_l := \norm{u_l
- u_{l,1}}_{L^2}$, for each $l \in \llbracket -M, N \rrbracket$.
This allows us to investigate the consequences of using the same mesh for all the
parametric problems.
In \cref{fig:2d_sines_parametric} we have plotted the exact parametric errors
after five steps of (uniform) refinement.
As we can notice, the same mesh leads to a wide range of parametric errors
values.
These errors are particularly low when the diffusion part of the parametric operator is dominant.
When the reaction part begin to domine the mesh seems to have a constant impact on the parametric errors (this is particularly striking for the BP scheme).

As expected these results suggest that the method can be optimized by using
different meshes depending on $l$.
In particular, coarser meshes would be sufficient when the diffusion coefficient is dominant.
These results are obtained for mesh uniform refinement, further investigations
deserve to be carried out for mesh adaptive refinement.

As we explained earlier, using a different hierarchy of meshes for each
parametric problem may be computationally advantageous, at the expense of ease
of implementation.
Several hierarchies of meshes would need to be stored and, in the case of
adaptive mesh refinement, interpolation between possibly non-nested
meshes would be required in order to compute the fractional solution $u$.
To avoid these complications when mesh adaptive refinement is used, we propose the
following:
\begin{enumerate}
    \item use the same hierarchy of meshes for all the parametric problems
        but not the same mesh.
        Some parametric problems might be solved on coarser meshes from the
        hierarchy and others on finer ones.
        This would allow to keep only one hierarchy of meshes stored in memory.
        Moreover, it would avoid the interpolation between non-nested meshes,
        since meshes from the same hierarchy are always nested.
    \item selectively refine the mesh hierarchy: estimate the error
        globally for each parametric problem (this can be done using the local
        parametric Bank--Weiser solutions) and mark the parametric problems for
        which a finer mesh is required, using e.g.\ a marking algorithm similar to
        Dörfler's marking strategy.
\end{enumerate}

\subsection{Two-dimensional checkerboard test case}\label{subsec:2d-checkerboard}
\begin{figure}
    \begin{center}
        \includegraphics[width=0.9\textwidth]{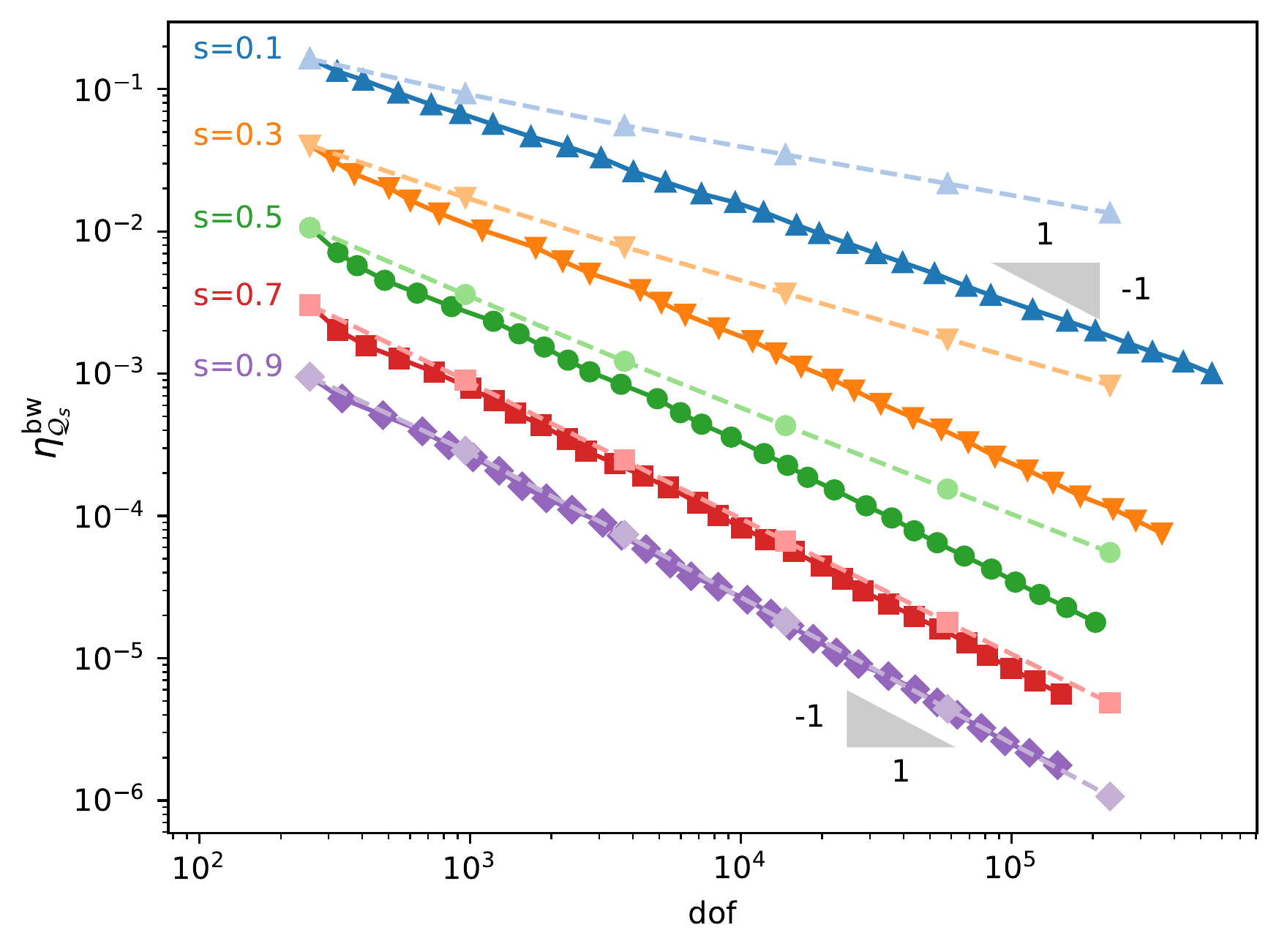}
    \end{center}
    \caption{\textbf{Two-dimensional checkerboard test case:} the dashed lines represent the Bank--Weiser estimator when uniform mesh refinement is performed and the solid lines represent the estimator when adaptive mesh refinement is performed.}
    \label{fig:2d-checkerboard_convergence}
\end{figure}

\begin{figure}
    \begin{center}
        \includegraphics[width=0.32\textwidth]{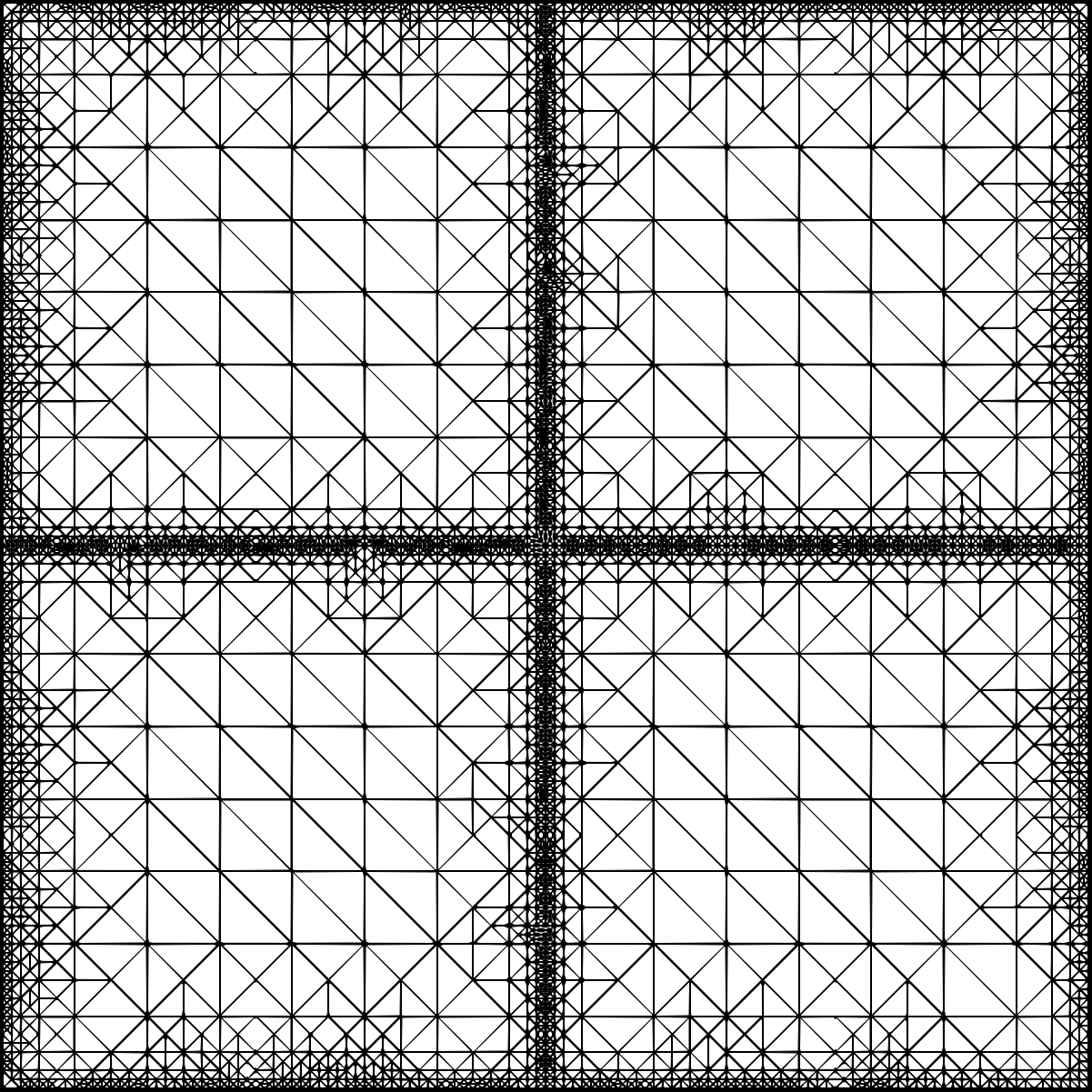}
        \includegraphics[width=0.32\textwidth]{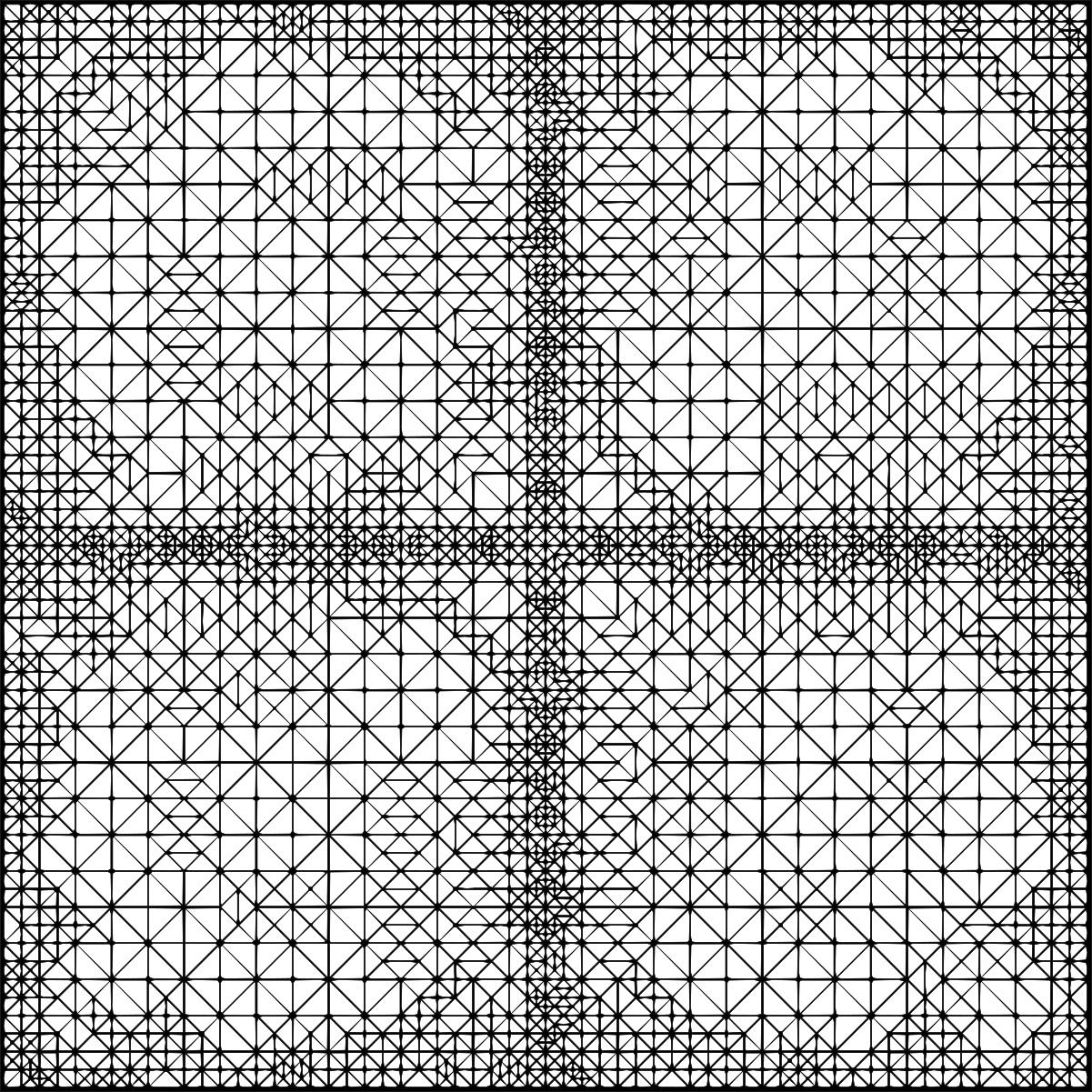}
        \includegraphics[width=0.32\textwidth]{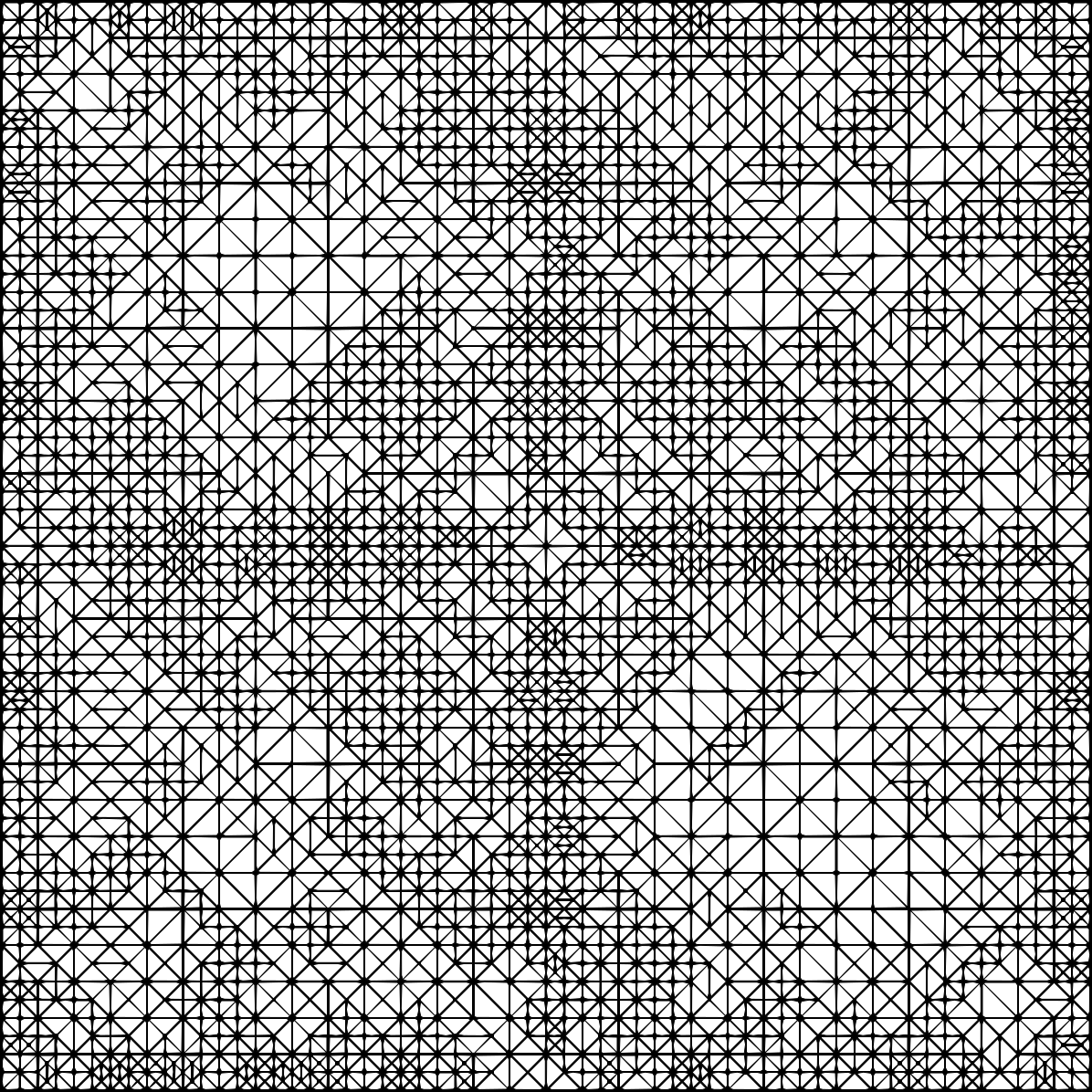}
    \end{center}
    \caption{\textbf{Two-dimensional checkerboard test case:} meshes
        obtained after 10 steps of mesh adaptive refinement steered by
    the Bank--Weiser estimator for $s=0.1,\ s=0.5$ and $s=0.9$ from left to right.}
    \label{fig:checkerboard_adapted_meshes}
\end{figure}

We solve the problem introduced in the numerical results of\
\cite{Bonito2013}.
We consider a unit square $\Omega = (0, 1)^2$ with data $f :
\Omega \to \R$ given for all $(x_1, x_2) \in \Omega$ by
\begin{equation}
f(x_1, x_2) = \begin{cases}
    1,& \text{if } (x_1 - 0.5)(x_2 - 0.5) > 0, \\
    -1,& \text{otherwise}.
\end{cases}
\end{equation}
The data $f$ belongs to $\mathbb H^{1/2 - \epsilon}(\Omega)$ for all $\epsilon > 0$.
So in Theorem 4.3 of \cite{Bonito2013} the index $\delta < 1/2$ and since
$\Omega$ is convex, again $\alpha$ can be chosen equal to $1$.
Then, the predicted convergence rate (for uniform refinement) is
$\ln(\sqrt{\mathrm{dof}})\mathrm{dof}^{-\beta}$ with
\begin{equation}
\beta = \begin{cases}
    1,& \text{if } s > \tfrac{3}{4}, \\
    s + \tfrac{1}{4},& \text{otherwise}.
\end{cases}
\label{eq:convergence_rate_beta}
\end{equation}
The predicted (if we omit the logarithmic term) and calculated convergence
rates for different choices of $s$ are given in
\cref{tab:2d-checkerboard_slopes}.
We recall that, to our knowledge, Theorem 4.3 is only available for BP and for BURA.
As we can see on this table, the convergence rates for the total
estimator is globally coherent with the predictions.
\Cref{tab:2d-checkerboard_slopes,fig:2d-checkerboard_convergence} show that mesh adaptive refinement improves the convergence rate for small fractional powers.
This is expected, the deterioration in the convergence rate is due to the
boundary layer behavior of the solution that is getting stronger as the fractional power decreases.
In the limit as the fractional power approaches $s \to 1$, the solution behaves like the solution to a non-fractional problem, i.e. there is no boundary layer and mesh adaptive refinement is no longer needed to improve the convergence rate. Nevertheless, it appears that our estimation strategy deals properly with this limit case in terms of recovering the expected convergence rates for all $s$.
This behavior can be seen on \cref{fig:checkerboard_adapted_meshes}, after 10 steps of
mesh adaptive refinement, the mesh associated to fractional power $s=0.9$ is almost
uniformly refined while the meshes associated to $s=0.5$ and $s=0.1$ show strongly localized refinement.
This explains why in \cref{fig:2d-checkerboard_convergence} we see the expected
behavior, i.e. no improvement in the convergence rate, when the mesh is
adaptively refined compared to uniformly refined when $s\geqslant 0.7$.

\begin{table}
    \begin{center}
    \begin{tabular}{llrrrrr}\toprule
        & Frac. power      & 0.1                        & 0.3                        & 0.5                        & 0.7                        & 0.9                       \\ \midrule
        BP   & Theory \cite{Bonito2013}      & -0.35                      & -0.55                      & -0.75                      & -0.95                      & -1.00                     \\ 
        & Unif. mesh ref.  & -0.35 & -0.56 & -0.77 & -0.94 & -1.00 \\ 
        & Adapt. mesh ref. & -0.66                      & -0.85                      & -0.95                      & -0.96                      & -1.02                     \\ \midrule
        BURA & Unif. mesh ref.  & -0.34                      & -0.59                      & -0.77                      & -0.94                      & -1.00                     \\
        & Adapt. mesh ref. & -0.54                      & -0.89                      & -0.95                      & -0.96                      & -1.02 \\ \bottomrule             
\end{tabular}
\end{center}
\caption{\textbf{Two-dimensional checkerboard test case:} convergence rates of the Bank--Weiser estimator for uniform mesh refinement and adaptive mesh refinement and for BP and BURA methods. In the case of the BP method, the convergence rates are compared with the values predicted by \cite{Bonito2013} for various fractional powers.}
    \label{tab:2d-checkerboard_slopes}
\end{table}

\subsection{Three-dimensional product of sines test case}\label{subsec:3d-sines}

\begin{figure}
    \begin{center}
        \includegraphics[width=0.45\textwidth]{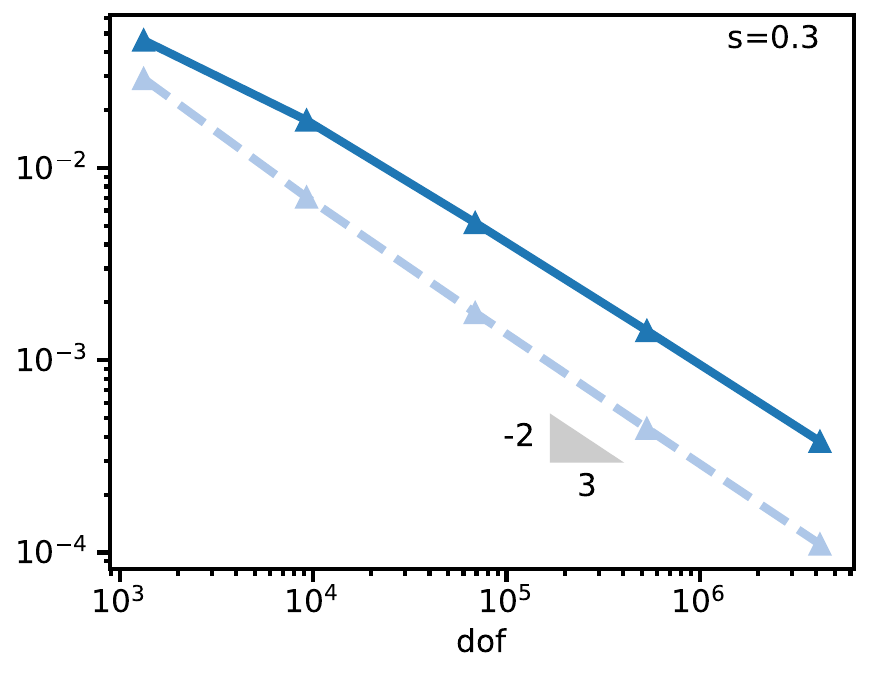}
        \includegraphics[width=0.45\textwidth]{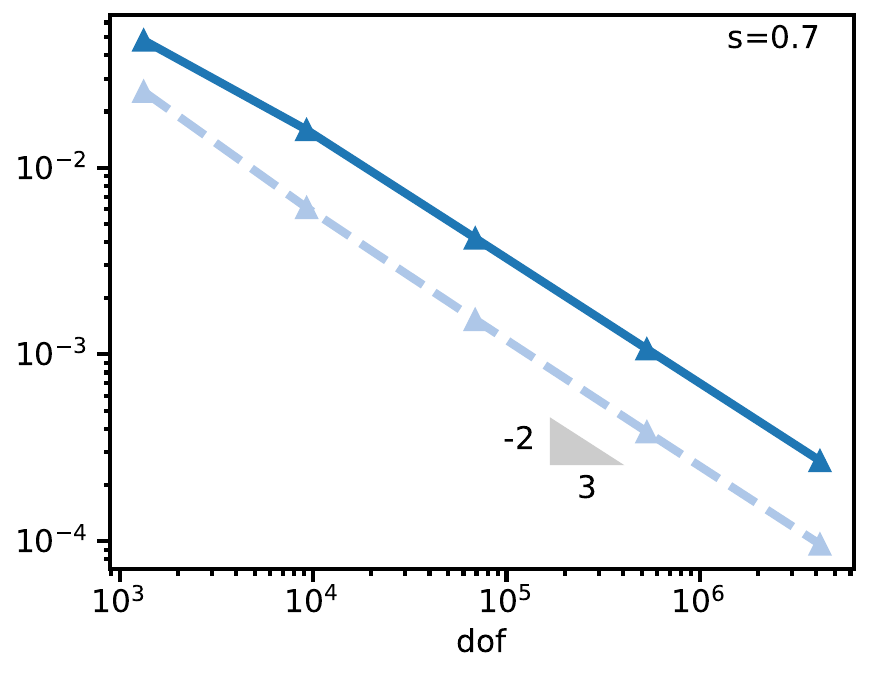}
    \end{center}
    \caption{\textbf{Three-dimensional product of sines test case:} the Bank--Weiser estimator
    $\eta^{\bw}_{\mathcal Q_s}$ in solid dark blue line is compared to the exact error in dashed light blue line for $s=0.3$ and $s=0.7$ when the BP rational scheme is used.}
\label{fig:3d_sines_convergence}
\end{figure}

This test case is the three-dimensional equivalent of the last test case.
We solve \cref{eq:fractional_strong_form} on the cube $\Omega = (0,\pi)^3$ with
data $f(x,y,z)=(2/\pi)^{3/2} \sin(x) \sin(y) \sin(z)$.
The analytical solution to this problem is given by $u(x,y,z) = 3^{-s}
(2/\pi)^{3/2} \sin(x) \sin(y) \sin(z)$.
The problem is solved on a hierarchy of uniformly refined Cartesian
(tetrahedral) meshes.
As for the two-dimensional case, the solution $u$ shows no boundary layer
behavior and mesh adaptive refinement is not required.
For the same reasons as for the two-dimensional case, Theorem 4.3 from
\cite{Bonito2013} predicts a convergence rate of $\mathrm{dof}^{-2/3}$ for the
finite element scheme.
\cref{fig:3d_sines_convergence} shows the values of the Bank--Weiser estimator
and of the exact error (computed from the knowledge of the analytical solution)
for $s=0.3$ and $s=0.7$.
As in the two-dimensional case, the efficiency indices are relatively robust
with respect to the fractional powers.
They are shown for various fractional powers in \cref{tab:3d_sines_efficiency}
and are computed by taking the average of the indices from the three last meshes
of the hierarchy.
As we can see, the Bank--Weiser estimator efficiency indices for this
three-dimensional case are not as good as in the two-dimensional case.
We have already observed this behavior for non-fractional problems
\cite{Bulle2020}.
We can notice that the convergence rates, given in \cref{tab:3d_sines_slopes},
are coherent with the predictions of Theorem 4.3 from \cite{Bonito2013}.
The convergence rates are computed from a linear regression on the values
computed from the three last meshes of the hierarchy.

\begin{table}
    \begin{center}
        \begin{tabular}{lrrrrr}
            Frac. power & 0.1 & 0.3 & 0.5 & 0.7 & 0.9\\
            \midrule
            Estimator & -0.56 & -0.60 & -0.63 & -0.65 & -0.66\\
            Exact error & -0.69 & -0.69 & -0.69 & -0.69 & -0.69
        \end{tabular}
    \end{center}
    \caption{\textbf{Three-dimensional product of sines test case:} convergence rates of the Bank--Weiser estimator and of the exact error for various fractional powers when the BP rational scheme is used.}
    \label{tab:3d_sines_slopes}
\end{table}

\begin{table}
    \begin{center}
        \begin{tabular}{r c c c c c}
            Frac. power & 0.1 & 0.3 & 0.5 & 0.7 & 0.9\\
            \midrule
            Est. eff. index & 2.12 & 3.20 & 3.08 & 2.77 & 2.45
        \end{tabular}
    \end{center}
    \caption{\textbf{Three-dimensional product of sines test case:} efficiency indices of
    the Bank--Weiser estimator for various fractional powers.}
    \label{tab:3d_sines_efficiency}
\end{table}

\subsection{Three-dimensional checkerboard test case}\label{subsec:3d-checkerboard}
\begin{figure}
    \begin{center}
        \includegraphics[width=0.5\textwidth]{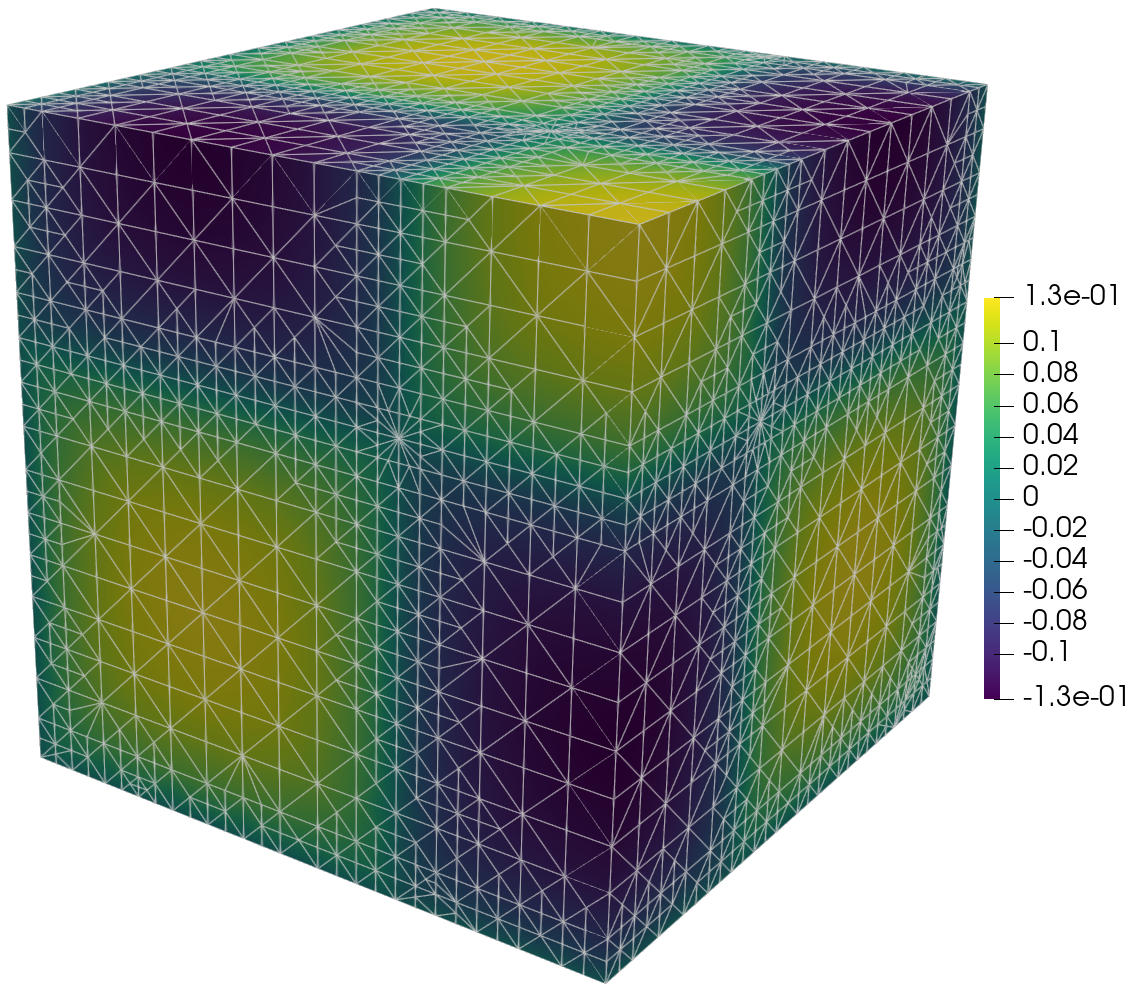}
    \end{center}
    \caption{\textbf{Three-dimensional checkerboard test case:} finite element
    solution and mesh after six steps of mesh adaptive refinement when $s=0.5$. The
    unit cube domain $(0,1)^3$ is truncated by the
    three planes passing through the point $(0.25, 0.25, 0.25)$ and
    orthogonal to the vectors $(1, 0, 0)$, $(0, 1, 0)$ and $(0, 0, 1)$
    respectively.}
    \label{fig:3d_checkerboard_solution_6}
\end{figure}

This test case is the three-dimensional version of the above checkerboard
problem.
We solve \cref{eq:fractional_strong_form} on the unit cube $\Omega = (0,1)^3$,
with data $f$ such that
\begin{equation}
    f(x_1, x_2, x_3) = \begin{cases}
        1, & \text{if } (x_1 - 0.5)(x_2 - 0.5) > 0 \text{ and } (x_3 - 0.5) < 0,\\
        1, & \text{if } (x_1 - 0.5)(x_2 - 0.5) < 0 \text{ and } (x_3 - 0.5) > 0,\\
        -1, & \text{otherwise.}
    \end{cases}
\end{equation}
The finite element solution $u_1$ and the corresponding mesh after six steps of
mesh adaptive refinement are shown in \cref{fig:3d_checkerboard_solution_6} for the
fractional power $s=0.5$.
As for the two-dimensional case, $f \in \mathbb H^{1/2 - \epsilon}(\Omega)$ for
all $\epsilon > 0$.
Consequently, once again Theorem 4.3 of \cite{Bonito2013} predicts a convergence
rate (for uniform refinement) equal to
$\ln\left(\mathrm{dof}^{1/3}\right)\mathrm{dof}^{-2\beta/3}$ with $\beta$ given by
\cref{eq:convergence_rate_beta}.

\begin{figure}
    \begin{center}
    \includegraphics[width=0.8\textwidth]{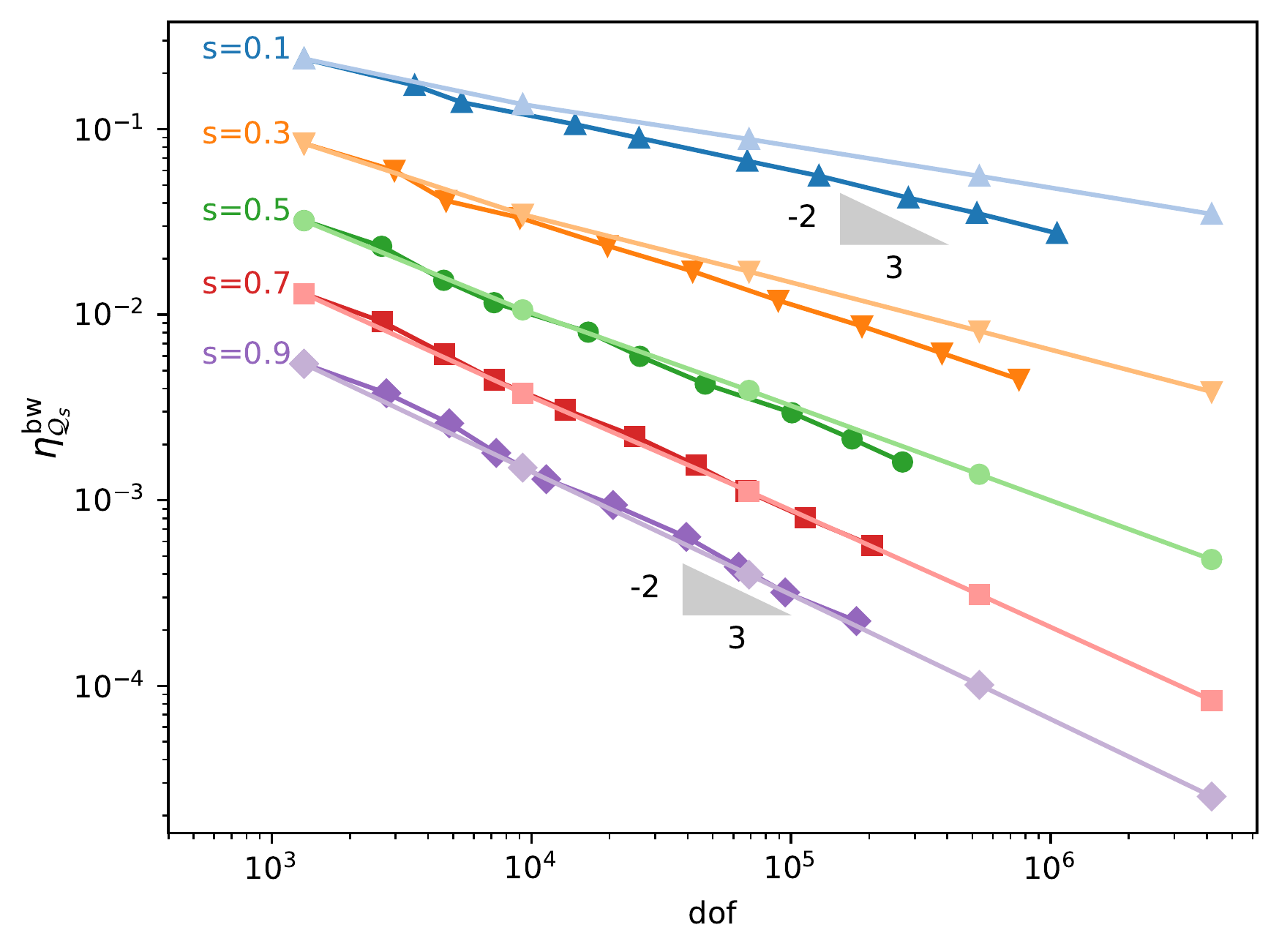}
    \end{center}\label{fig:3d_checkerboard_convergence}
    \caption{\textbf{Three-dimensional checkerboard test case:} for each fractional power we compare the
        values of the Bank--Weiser estimator $\eta^{\bw}_{\mathcal Q_s}$ when uniform
        refinement is performed (light lines) and when adaptive
        refinement is performed (darker lines).}
\end{figure}

Once again, if we omit the logarithmic term, the predicted and calculated
convergence rates are given in \cref{tab:3d_checkerboard_slopes}.
As in the two-dimensional case, the convergence rates of the Bank--Weiser
estimator are globally coherent with the predictions and the boundary layer
behavior becomes stronger as the fractional power decreases leading to poorer
convergence rates.
\cref{fig:3d_checkerboard_convergence} shows the values of the Bank--Weiser
estimator for mesh uniform and adaptive refinement and for several fractional powers.

\begin{table}
    \begin{center}
        \begin{tabular}{l r r r r r}
            Frac. power & 0.1 & 0.3 & 0.5 & 0.7 & 0.9\\
            \midrule
            Theory \cite{Bonito2013} & -0.23 & -0.37 & -0.50 & -0.63 & -0.67\\
            Unif. mesh ref. & -0.24 & -0.38 & -0.52 & -0.62 & -0.67\\
            Adapt. mesh ref. & -0.33 & -0.46 & -0.55 & -0.65 & -0.68\\
        \end{tabular}
    \end{center}
    \caption{\textbf{Three-dimensional checkerboard test case:} convergence rates of the
    Bank--Weiser estimator for mesh uniform refinement and for mesh adaptive refinement
    compared to the values predicted by \cite{Bonito2013} for various fractional powers when the BP rational scheme is used.}
    \label{tab:3d_checkerboard_slopes}
\end{table}

\subsection{Three--dimensional torus test case}\label{subsec:3d-torus}
We show now the application of our method on an
unstructured three-dimensional mesh. We generate a mesh of a filled torus with
unit inner radius and outer radius of $1.3$ using
gmsh~\cite{geuzaine_gmsh:2009}. The initial mesh has 12040 tetrahedral cells.
We set $s = 0.5$ and $f = 1$ for all $\Omega$. We estimate the lowest
eigenvalue of the standard Laplacian on this mesh $\lambda_0 \simeq 60$ and calculate of the data $\norm{f}_{L^2} \simeq 1.323$ giving a quadrature fineness parameter $\kappa = 0.48$ for the BP scheme (45 quadrature points).

\begin{figure}
\begin{center}
\includegraphics[width=0.8\textwidth]{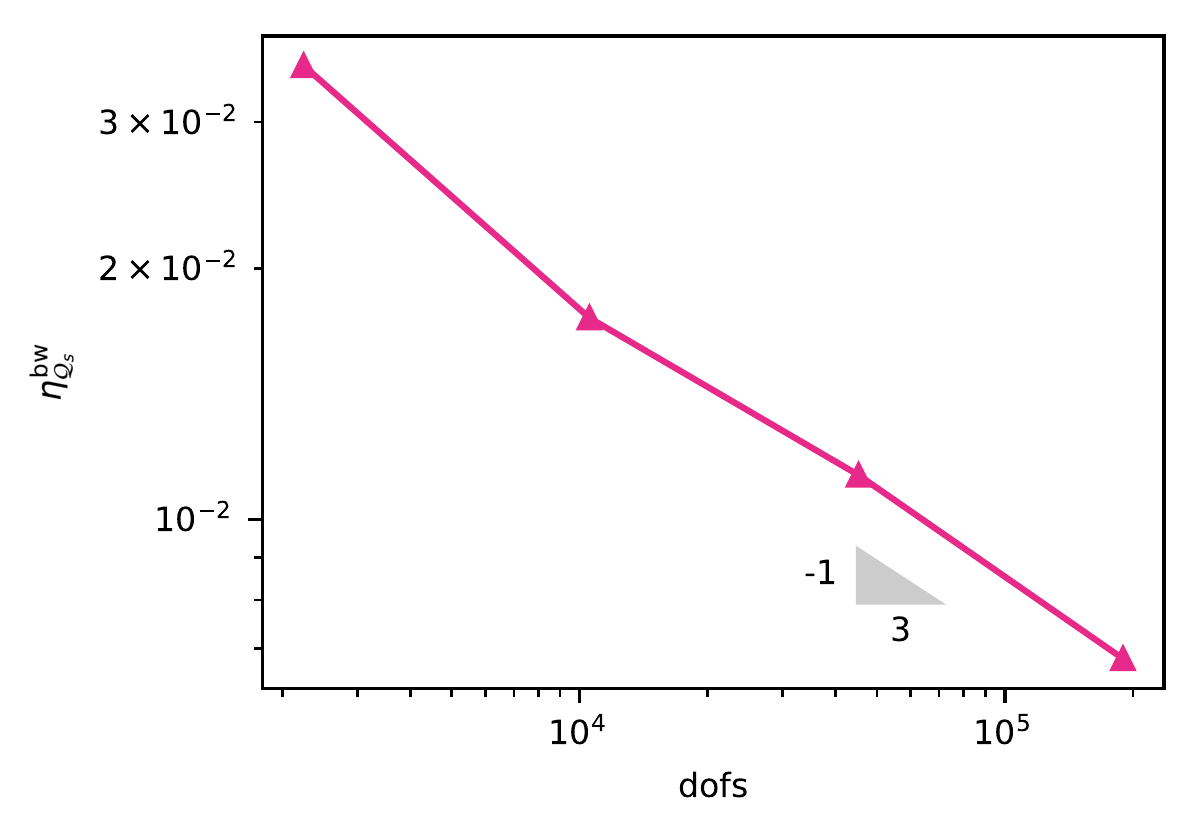}
	\caption{\textbf{Three--dimensional torus test case:} convergence of estimator under adaptive mesh refinement.}\label{fig:3d_torus_convergence}
\end{center}
\end{figure}

\begin{figure}
\begin{center}
	\includegraphics[width=0.8\textwidth]{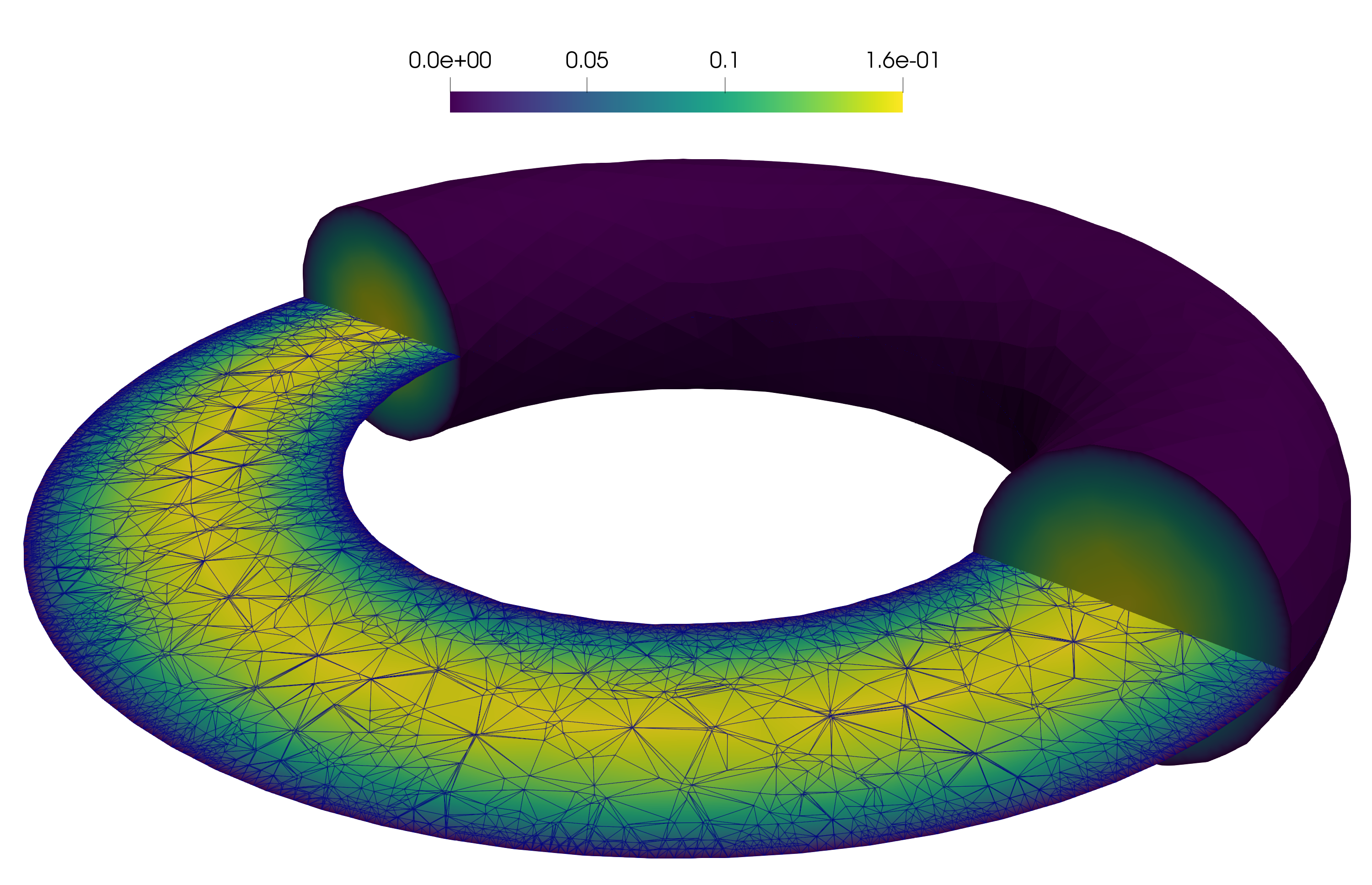}
	\caption{\textbf{Three--dimensional torus test case:} Plot of solution
	after the final mesh adaptive refinement. Clipped with plane passing through
	origin with normal $(1, 0, 0)$. Sliced with plane passing through
	origin with normal $(0, 0, 1)$. Slice shows mesh refinement with
	particularly strong refinement to resolve the boundary
	layer.}\label{fig:3d_torus_solution}
\end{center}
\end{figure}

We perform adaptive mesh refinement until a discretization error of less than $1 \times 10^{-2}$ has been reached. The evolution of the estimator is given in
\cref{fig:3d_torus_convergence}. A plot of the solution and an idea of the mesh
refinement on a cut at the third and final mesh is given in
\cref{fig:3d_torus_solution}. Clearly evident is the very strong refinement
near the boundary needed to capture the strong boundary layer in this problem.

\subsection{Adaptive rational scheme}\label{subsec:ra_adaptive_results}

In this section, we apply our method combining finite element mesh and rational scheme adaptive refinement (sketched in \cref{fig:fe_ra_refinement_algorithm}) to the two-dimensional test cases from \cref{subsec:2d-checkerboard,subsec:2d-sines}.
We provide a comparison between this method and the method where the rational schemes are fixed.
When an adaptive rational scheme is used the number of parametric problems is variable from one step to another.
Thus, in order to obtain a meaningful comparison, we use the total number of degrees of freedom (i.e.\ the number of degrees of freedom times the number of parametric problems) on the plots x--axis.
In addition, the rational approximation error is no longer negligible so our interest is on the total error estimator and the total discretization error (respectively defined in \cref{eq:total_estimator} and \cref{eq:total_discretization_error}) rather than the Bank--Weiser estimator and the finite element discretization error.

\subsubsection{Two-dimensional product of sines test case}\label{subsubsec:2d-pdt-sines}

\begin{figure}
    \begin{center}
        \includegraphics[width=0.45\textwidth]{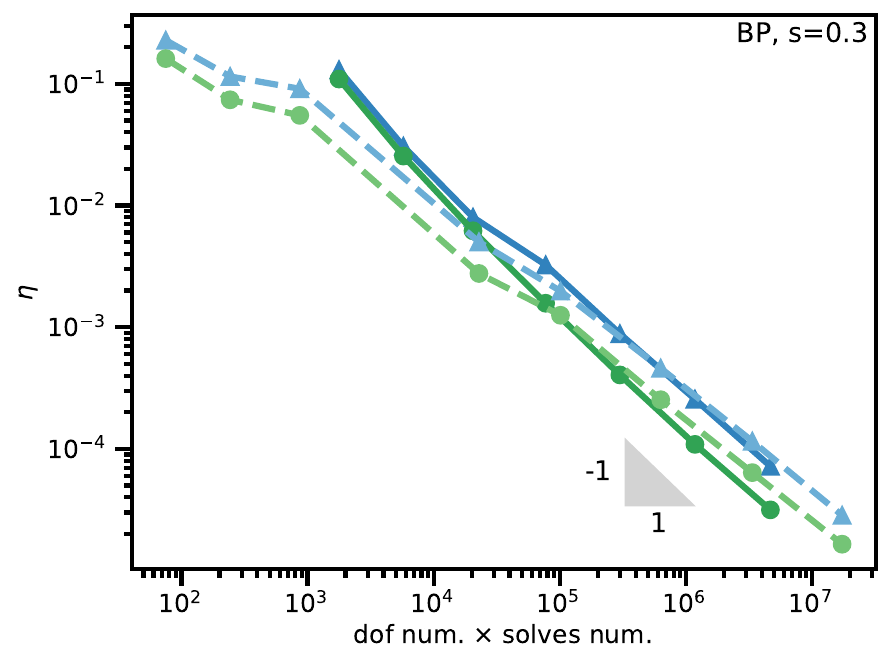}
        \includegraphics[width=0.45\textwidth]{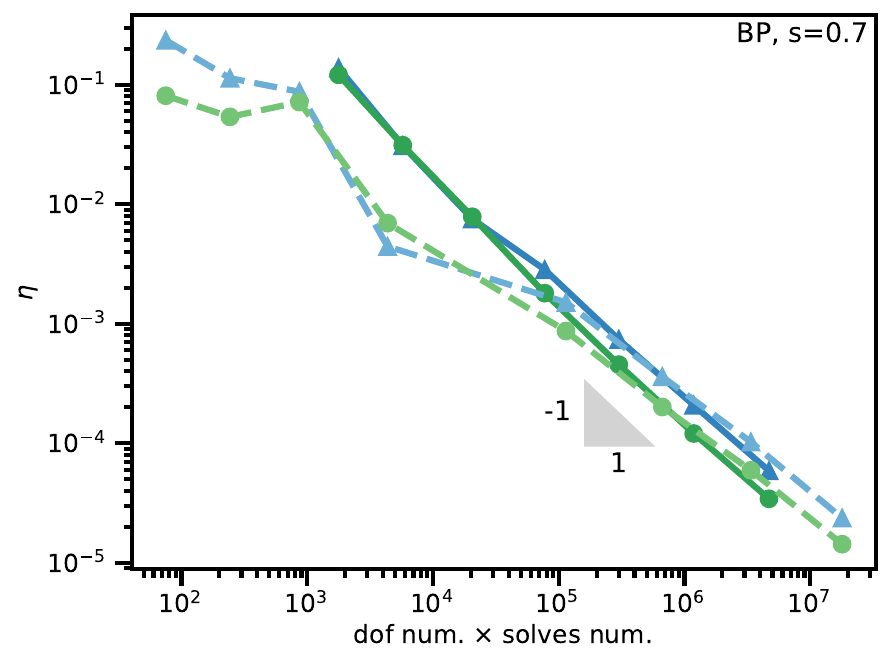}\\
        \includegraphics[width=0.45\textwidth]{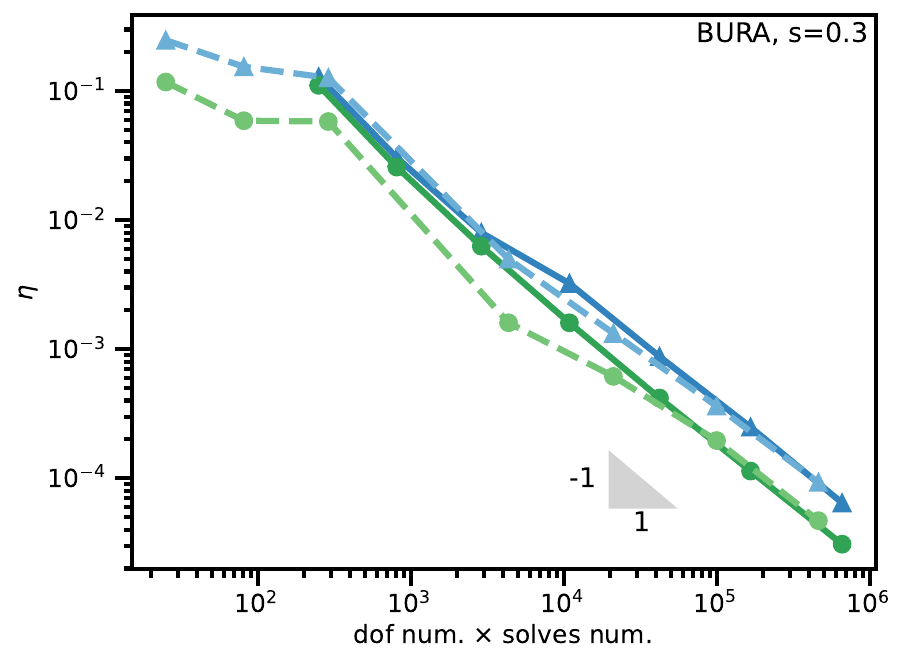}
        \includegraphics[width=0.45\textwidth]{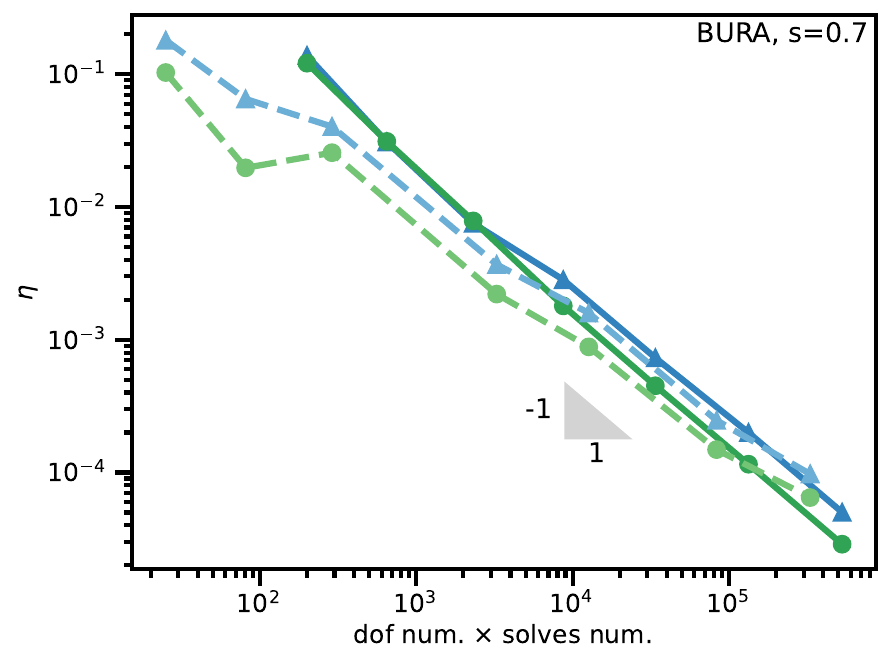}
    \end{center}
    \caption{\textbf{Two-dimensional product of sines test case:} the solid lines represent fixed BP and BURA schemes while the dashed lines represent adaptive BP and BURA schemes. Triangular markers represent the values of the total estimator and circular markers represent the values of the exact total error.}
    \label{fig:2d_sines_cumulative_convergence}
\end{figure}

In \cref{fig:2d_sines_cumulative_convergence} we compare the values of the total estimator $\eta$ and of the exact error with and without the use of an adaptive rational scheme.
As in \cref{subsec:2d-sines}, adaptive mesh refinement is not required on this test case so we only perform uniform mesh refinement.

In \cref{tab:2d_sines_ra_adaptive_param_problems} we compare the number of parametric problems solved with or without the use of an adaptive rational scheme.
The total numbers of parametric problems to solve in order to reach our tolerance is reduced by 30 to 43 \% for the BP method and by 40 to 48 \% for the BURA method.
However, and more surprisingly, the use of an adaptive rational scheme does not induce any gain in the overall computational cost --measured by the number of degrees of freedom times the number of parametric problems-- for this test case, as we can see in \cref{fig:2d_sines_cumulative_convergence}.
We notice that the introduction of the adaptive rational scheme slightly deteriorate the convergence rates, shown in \cref{tab:2d_sines_ra_adaptive_slopes}.
Our method consists in reducing the convergence rate of the rational approximation to the convergence rate of the Bank--Weiser estimator.
Thus, it assumes that the rational scheme does not influence the Bank--Weiser estimator rate of convergence.
However, the results in \cref{tab:2d_sines_ra_adaptive_slopes} suggest that the Bank--Weiser estimator is influenced by the rational approximation scheme.

\begin{table}[]
    \begin{center}
    \begin{tabular}{llrrrrr}
    \toprule
         & Frac. power          & 0.1  & 0.3 & 0.5 & 0.7 & 0.9  \\ \midrule
    BP   & Fixed ra. scheme     & 1155 & 497 & 427 & 497 & 1155 \\ 
         & Adaptive ra. scheme  & 504  & 209 & 178 & 199 & 358  \\ \midrule
    BURA & Fixed ra. scheme     & 96   & 77  & 63  & 49  & 35   \\
         & Adaptive ra. scheme  & 42   & 33  & 29  & 20  & 17   \\
    \bottomrule
    \end{tabular}
    \end{center}
    \caption{\textbf{Two-dimensional product of sines test case:} comparison of the total number of solves to reach the tolerance for each method and each fractional power with and without the adaptively refined rational scheme.}
    \label{tab:2d_sines_ra_adaptive_param_problems}
\end{table}

\begin{table}
    \begin{center}
        \begin{tabular}{lllrrrrr}
            \toprule
                 &              & Frac. power &   0.1 &   0.3 &   0.5 &   0.7 &   0.9 \\
            \midrule
            BP & Fixed scheme & Exact error & -1.03 & -1.03 & -1.04 & -1.04 & -1.04 \\
            &              & Estimator & -0.92 & -0.93 & -0.95 & -0.96 & -0.97 \\
               & Adaptive scheme & Exact error & -0.83 & -0.75 & -0.77 & -0.73 & -0.74 \\
            &              & Estimator & -0.78 & -0.74 & -0.83 & -0.74 & -0.77 \\
                 
            \midrule
            BURA & Fixed scheme & Exact error & -0.80 & -1.03 & -1.05 & -1.06 & -1.06 \\
            &              & Estimator & -0.78 & -0.94 & -0.96 & -0.98 & -0.99 \\
                & Adaptive scheme & Exact error & -0.67 & -0.83 & -0.84 & -0.77 & -0.77 \\
            &              & Estimator & -0.71 & -0.85 & -0.83 & -0.80 & -0.79 \\
                 
            \bottomrule
        \end{tabular}
    \end{center}
    \caption{\textbf{Two-dimensional product of sines test case:} convergence rates of the total estimator and of the exact total error for various fractional powers and for fixed or adaptive BP and BURA schemes.}
    \label{tab:2d_sines_ra_adaptive_slopes}
\end{table}

\begin{table}
    \begin{center}
        \begin{tabular}{llrrrrr}
            \toprule
                 & Frac. power      & 0.1  & 0.3   & 0.5   & 0.7   & 0.9 \\
            \midrule
            BP   & Fixed scheme     & 1.76 &  2.03 &  1.79 &  1.52 &  1.29 \\
                 & Adaptive scheme  & 1.69 &  1.75 &  1.73 &  1.51 &  1.85 \\
            BURA & Fixed scheme     & 1.02 &  1.93 &  1.81 &  1.53 &  1.30 \\
                 & Adaptive scheme  & 1.44 &  2.25 &  2.00 &  1.64 &  2.07 \\
            \bottomrule
        \end{tabular}
    \end{center}
    \caption{\textbf{Two-dimensional product of sines test case:} efficiency indices of the total estimator for various fractional powers and for fixed or adaptive BP and BURA schemes.}
    \label{tab:2d_sines_ra_adaptive_efficiency}
\end{table}

\subsubsection{Two-dimensional checkerboard test case}\label{subsubsec:2d-checkerboard}

\begin{figure}
    \begin{center}
        \includegraphics[width=0.45\textwidth]{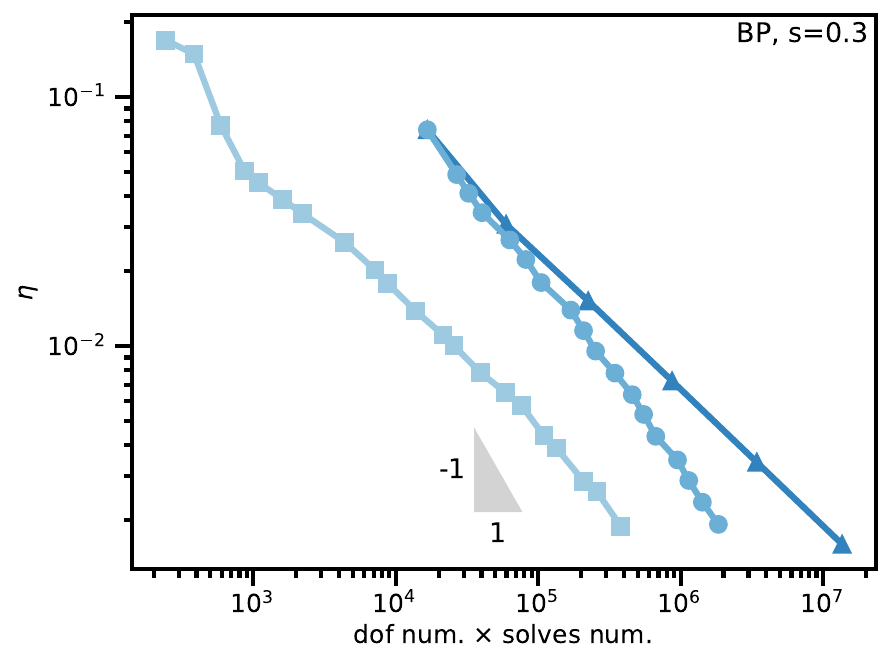}
        \includegraphics[width=0.45\textwidth]{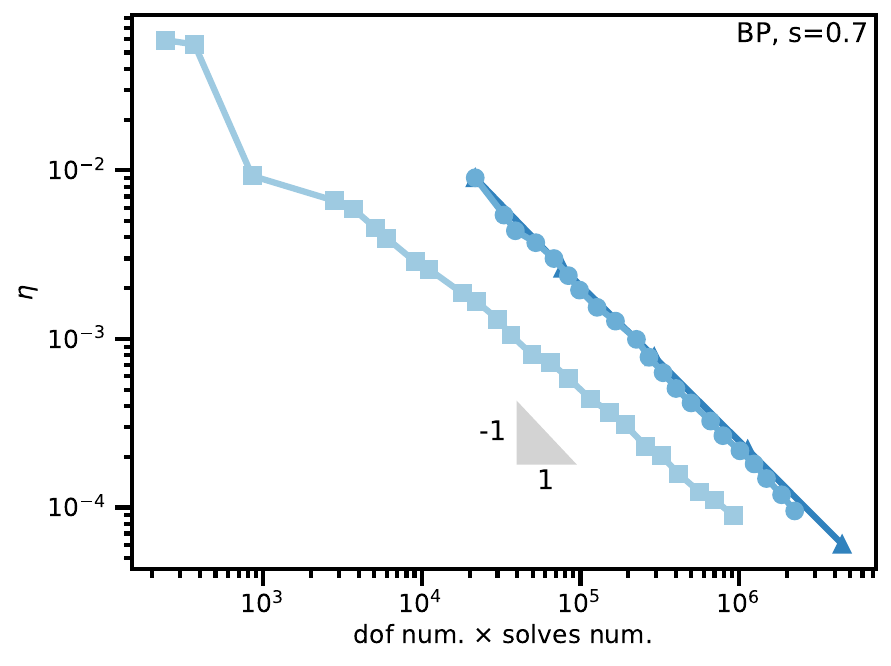}\\
        \includegraphics[width=0.45\textwidth]{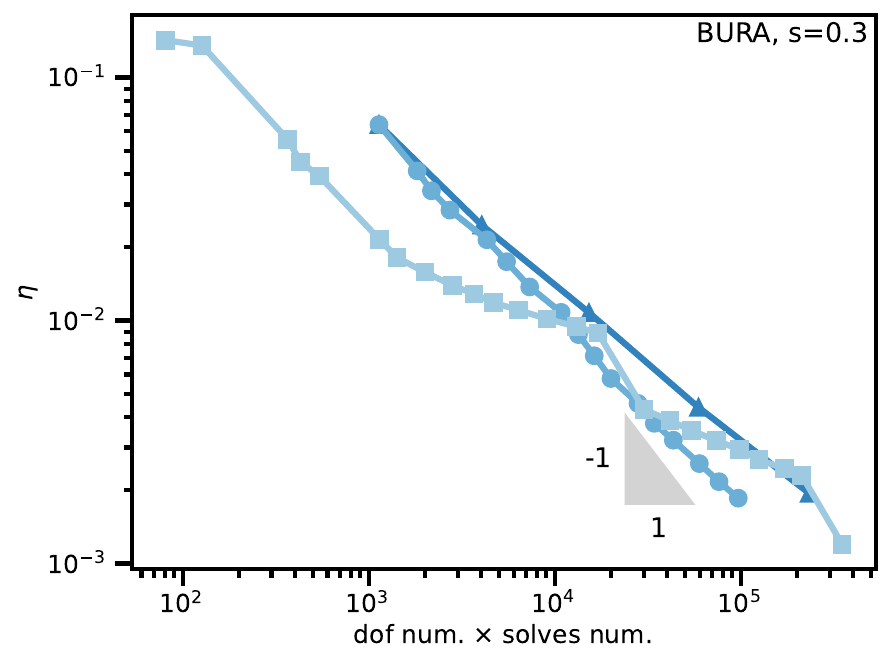}
        \includegraphics[width=0.45\textwidth]{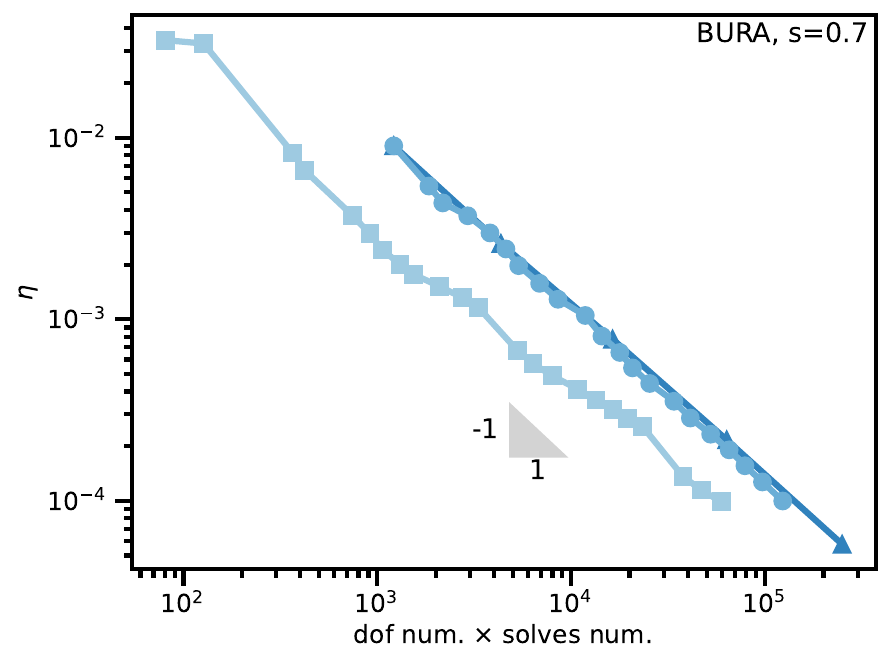}
    \end{center}
    \caption{\textbf{Two-dimensional checkerboard test case:} we compare the values of the total estimator for three refinement strategies, when the mesh is uniformly refined and a fixed rational scheme is used (triangular markers), when the mesh is adaptively refined and a fixed rational scheme is used (circular markers) and when the mesh is adaptively refined and an adaptive rational scheme is used (square markers).}
    \label{fig:2d-checkerboard_cumulative_convergence}
\end{figure}

In \cref{fig:2d-checkerboard_cumulative_convergence} we compare the values of the total estimator $\eta$ and of the exact error with and without the use of an adaptive rational scheme.
As in \cref{subsec:2d-checkerboard}, the boundary layer behavior, present for small values of $s$, requires adaptive mesh refinement but gets less and less prominent as $s$ tends to 1.

The gain in the number of parametric problems to solve is comparable to what we obtain in \cref{tab:2d_sines_ra_adaptive_param_problems}.
The use of an adaptive rational scheme clearly induces a gain in the precision with respect to the total number of degrees of freedom, except for 
small fractional powers when the BURA rational scheme is used.

As for the two--dimensional product of sines test case, we can notice on \cref{tab:2d-checkerboard_ra_adaptive_slopes} that the use of an adaptive rational scheme also deteriorates the convergence rates in most of the cases.
However, unlike the test case from \cref{subsubsec:2d-pdt-sines}, a gain in the computational cost clearly appears when an adaptive rational scheme is used, especially with the BP method.

The discrepancies between the BP and BURA methods might be due to the high sensitivity of the BURA scheme with respect to the parameter $N$: balancing the rational and finite element estimator is more difficult with the BURA scheme.

\begin{table}
    \begin{center}
        \begin{tabular}{llrrrrr}
            \toprule
                 & Frac. power &   0.1 &   0.3 &   0.5 &   0.7 &   0.9 \\
            \midrule
            BP  & Theory \cite{Bonito2013} & -0.35 & -0.55 & -0.75 & -0.95 & -1.00 \\
                & Unif. mesh ref. & -0.35 & -0.56 & -0.77 & -0.94 & -1.00 \\
                & Adapt. mesh ref. & -0.57 & -0.77 & -0.94 & -0.97 & -1.00 \\
                & Adapt. mesh ref. \& ra. sch. & -0.44 & -0.56 & -0.75 & -0.78 & -0.82 \\
            \midrule
            BURA & Unif. mesh ref. & -0.21 & -0.65 & -0.83 & -0.94 & -1.00 \\
                 & Adapt. mesh ref. & -0.36 & -0.79 & -0.97 & -0.96 & -1.00 \\
                 & Adapt. mesh ref. \& ra. sch. & -0.45 & -0.50 & -0.65 & -0.80 & -0.86 \\
            \bottomrule
        \end{tabular}
    \end{center}
    \caption{\textbf{Two-dimensional checkerboard test case:} convergence rates of the total estimator for uniform mesh refinement, adaptive mesh refinement and adaptive mesh refinement combined with rational scheme adaptation respectively for BP and BURA methods.
    The convergence rates are compared with the values predicted by \cite{Bonito2013} (in the case of uniform refinement) for various fractional powers.}
    \label{tab:2d-checkerboard_ra_adaptive_slopes}
\end{table}

\section{Concluding remarks}
In this work we presented a novel a posteriori error estimation method for the spectral fractional Laplacian.
This method benefits from the parallel character of both the
Bank--Weiser error estimator and the rational approximation methods, thus keeping
the appealing computational aspects of the underlying methodology in\ \cite{Bonito2013}.
Here are some important points we want to make to conclude this paper.
First, the Bank--Weiser estimator seems to be equivalent to the $L^2$ exact
error at least when structured meshes are used and when the solution $u$ is
smooth.
Second, adaptive mesh refinement methods drastically improves the convergence rate compared to uniform refinement for fractional powers close to $0$.
Third, the use of an adaptive rational scheme can reduce the number of total parametric problems to solve by $30$ to $48 \%$ depending on the scheme and the fractional power.
However, it does not necessarily induce a gain in the precision with respect to the total number of degrees of freedom.

Finally, we give some future directions that we think are worth considering.
More numerical tests could be performed, especially for higher order elements
and/or using variants of the Bank--Weiser error estimator as considered in
\cite{Bulle2021a}.

We would like also to study the derivation of an algorithm that allows to use
different meshes to discretize the parametric problems in order to save
computational time, as explained in \cref{subsubsec:parametric_problems}.
However, the implementation of this strategy requires to interpolate functions between different meshes, which is not currently available in the FEniCSx software and by consequence, is beyond the scope of this study.
The use of an adaptive rational scheme and in particular the dependence of the Bank--Weiser estimator with respect to the choice of the rational scheme deserve a deeper investigation in order to fully explain the results in \cref{subsubsec:2d-pdt-sines,subsubsec:2d-checkerboard}.
The a posteriori error estimation of the error in the ``natural'' norm of the
problem i.e.\ the spectral fractional norm defined in
\cref{eq:spectral_sobolev} is another extension of this work that is worth to
consider.
The replacement of the Bank--Weiser estimator by an anisotropic a posteriori
error estimator would improve the convergence rate even further in case of
boundary layers, see e.g.\ \cite{Banjai2019,Faustmann2020},
Another interesting extension would be to test our method on fractional powers
of other kinds of elliptic operators, following \cite{Bonito2013}, on another
definition of the fractional Laplacian operator \cite{Bonito2019a} and/or other
boundary conditions, following \cite{Antil2018}.

\section*{Supplementary material}
A minimal example of adaptive finite element method for the two--dimensional
spectral fractional Laplacian and the three--dimensional torus can be found in
the following FEniCSx--Error--Estimation repository
\url{https://github.com/jhale/fenicsx-error-estimation}. This minimal example
code (LGPLv3) is also archived at
\url{https://doi.org/10.6084/m9.figshare.19086695.v3}. A Docker image
\cite{hale_containers_2017} is provided in which this code can be executed.

\section*{Acknowledgements}
R.B.\ would like to acknowledge the support of the ASSIST research project of
the University of Luxembourg. This publication has been prepared in the
framework of the DRIVEN TWINNING project funded by the European Union's Horizon
2020 Research and Innovation programme under Grant Agreement No.\ 811099. F.C. is grateful of the Center for Mathematical Modeling grant FB20005. His
work is partially supported by the I-Site BFC project NAANoD and the EIPHI
Graduate School (contract ANR-17-EURE-0002).

\bibliographystyle{abbrvnat}
\bibliography{ms}

\begin{thebibliography}{10}

\bibitem{Aceto2019}
Lidia Aceto, Daniele Bertaccini, Fabio Durastante, and Paolo Novati.
\newblock {Rational Krylov methods for functions of matrices with applications
  to fractional partial differential equations}.
\newblock {\em J. Comput. Phys.}, 396:470--482, nov 2019.
\newblock \href {https://doi.org/10.1016/j.jcp.2019.07.009}
  {\path{doi:10.1016/j.jcp.2019.07.009}}.

\bibitem{Aceto2017}
Lidia Aceto and Paolo Novati.
\newblock {Rational Approximation to the Fractional Laplacian Operator in
  Reaction-Diffusion Problems}.
\newblock {\em SIAM J. Sci. Comput.}, 39(1):A214--A228, jan 2017.
\newblock \href {http://arxiv.org/abs/1607.04166} {\path{arXiv:1607.04166}},
  \href {https://doi.org/10.1137/16M1064714} {\path{doi:10.1137/16M1064714}}.

\bibitem{Aceto2018}
Lidia Aceto and Paolo Novati.
\newblock {Rational approximations to fractional powers of self-adjoint
  positive operators}.
\newblock {\em Numer. Math.}, 143(1):1--16, sep 2019.
\newblock \href {http://arxiv.org/abs/1807.10086} {\path{arXiv:1807.10086}},
  \href {https://doi.org/10.1007/s00211-019-01048-4}
  {\path{doi:10.1007/s00211-019-01048-4}}.

\bibitem{Acosta2017a}
Gabriel Acosta, Francisco~M Bersetche, and Juan~Pablo Borthagaray.
\newblock {A short FE implementation for a 2d homogeneous Dirichlet problem of
  a fractional Laplacian}.
\newblock {\em Comput. Math. with Appl.}, 74(4):784--816, aug 2017.
\newblock \href {http://arxiv.org/abs/1610.05558} {\path{arXiv:1610.05558}},
  \href {https://doi.org/10.1016/j.camwa.2017.05.026}
  {\path{doi:10.1016/j.camwa.2017.05.026}}.

\bibitem{Acosta2017}
Gabriel Acosta and Juan~Pablo Borthagaray.
\newblock {A Fractional Laplace Equation: Regularity of Solutions and Finite
  Element Approximations}.
\newblock {\em SIAM J. Numer. Anal.}, 55(2):472--495, jan 2017.
\newblock \href {http://arxiv.org/abs/1507.08970} {\path{arXiv:1507.08970}},
  \href {https://doi.org/10.1137/15M1033952} {\path{doi:10.1137/15M1033952}}.

\bibitem{Ainsworth2017a}
Mark Ainsworth and Christian Glusa.
\newblock {Aspects of an adaptive finite element method for the fractional
  Laplacian: A priori and a posteriori error estimates, efficient
  implementation and multigrid solver}.
\newblock {\em Comput. Methods Appl. Mech. Eng.}, 327:4--35, dec 2017.
\newblock \href {http://arxiv.org/abs/1708.03912} {\path{arXiv:1708.03912}},
  \href {https://doi.org/10.1016/j.cma.2017.08.019}
  {\path{doi:10.1016/j.cma.2017.08.019}}.

\bibitem{Ainsworth2018a}
Mark Ainsworth and Christian Glusa.
\newblock {Hybrid Finite Element--Spectral Method for the Fractional Laplacian:
  Approximation Theory and Efficient Solver}.
\newblock {\em SIAM J. Sci. Comput.}, 40(4):A2383--A2405, jan 2018.
\newblock URL: \url{https://epubs.siam.org/doi/10.1137/17M1144696}, \href
  {http://arxiv.org/abs/1709.01639} {\path{arXiv:1709.01639}}, \href
  {https://doi.org/10.1137/17M1144696} {\path{doi:10.1137/17M1144696}}.

\bibitem{Ainsworth2018}
Mark Ainsworth and Christian Glusa.
\newblock {Towards an Efficient Finite Element Method for the Integral
  Fractional Laplacian on Polygonal Domains}.
\newblock {\em Contemp. Comput. Math. - A Celebr. 80th Birthd. Ian Sloan},
  40:17--57, 2018.
\newblock \href {http://arxiv.org/abs/1708.01923} {\path{arXiv:1708.01923}},
  \href {https://doi.org/10.1007/978-3-319-72456-0_2}
  {\path{doi:10.1007/978-3-319-72456-0_2}}.

\bibitem{Ainsworth2017}
Mark Ainsworth and Zhiping Mao.
\newblock {Analysis and Approximation of a Fractional Cahn--Hilliard Equation}.
\newblock {\em SIAM J. Numer. Anal.}, 55(4):1689--1718, jan 2017.
\newblock \href {https://doi.org/10.1137/16M1075302}
  {\path{doi:10.1137/16M1075302}}.

\bibitem{Ainsworth2000}
Mark Ainsworth and J~Tinsley Oden.
\newblock {\em {A Posteriori Error Estimation in Finite Element Analysis}}.
\newblock Pure and Applied Mathematics (New York). John Wiley \& Sons, Inc.,
  Hoboken, NJ, USA, aug 2000.
\newblock \href {https://doi.org/10.1002/9781118032824}
  {\path{doi:10.1002/9781118032824}}.

\bibitem{Akagi2016}
Goro Akagi, Giulio Schimperna, and Antonio Segatti.
\newblock {Fractional Cahn–Hilliard, Allen–Cahn and porous medium
  equations}.
\newblock {\em J. Differ. Equ.}, 261(6):2935--2985, sep 2016.
\newblock \href {http://arxiv.org/abs/1502.06383} {\path{arXiv:1502.06383}},
  \href {https://doi.org/10.1016/j.jde.2016.05.016}
  {\path{doi:10.1016/j.jde.2016.05.016}}.

\bibitem{Allaire2007}
Grégoire Allaire.
\newblock {\em Numerical analysis and optimization: an introduction to
  mathematical modelling and numerical simulation}.
\newblock Numerical mathematics and scientific computation. Oxford University
  Press, Oxford ; New York, 2007.
\newblock OCLC: ocm82671667.

\bibitem{Alnaes2015}
Martin~S Aln{\ae}s, Jan Blechta, Johan Hake, August Johansson, Benjamin Kehlet,
  Anders Logg, Chris Richardson, Johannes Ring, Marie~E Rognes, and Garth~N
  Wells.
\newblock {The FEniCS Project Version 1.5}.
\newblock {\em Arch. Numer. Softw.}, 2015.
\newblock \href {https://doi.org/10.11588/ans.2015.100.20553}
  {\path{doi:10.11588/ans.2015.100.20553}}.

\bibitem{Antil2015}
Harbir Antil and Enrique Ot{\'{a}}rola.
\newblock {A FEM for an Optimal Control Problem of Fractional Powers of
  Elliptic Operators}.
\newblock {\em SIAM J. Control Optim.}, 53(6):3432--3456, jan 2015.
\newblock \href {http://arxiv.org/abs/1406.7460} {\path{arXiv:1406.7460}},
  \href {https://doi.org/10.1137/140975061} {\path{doi:10.1137/140975061}}.

\bibitem{Antil2018}
Harbir Antil, Johannes Pfefferer, and Sergejs Rogovs.
\newblock {Fractional operators with inhomogeneous boundary conditions:
  analysis, control, and discretization}.
\newblock {\em Commun. Math. Sci.}, 16(5):1395--1426, 2018.
\newblock \href {http://arxiv.org/abs/1703.05256} {\path{arXiv:1703.05256}},
  \href {https://doi.org/10.4310/CMS.2018.v16.n5.a11}
  {\path{doi:10.4310/CMS.2018.v16.n5.a11}}.

\bibitem{Atangana2018}
Abdon Atangana.
\newblock {Fractional Operators and Their Applications}.
\newblock In {\em Fract. Oper. with Constant Var. Order with Appl. to
  Geo-Hydrology}, pages 79--112. Elsevier, 2018.
\newblock \href {https://doi.org/10.1016/B978-0-12-809670-3.00005-9}
  {\path{doi:10.1016/B978-0-12-809670-3.00005-9}}.

\bibitem{Balakrishnan1959}
Alampallam~V Balakrishnan.
\newblock {Fractional powers of closed operators and the semigroups generated
  by them}.
\newblock {\em Pacific J. Math.}, 10(2):419--437, jun 1960.
\newblock \href {https://doi.org/10.2140/pjm.1960.10.419}
  {\path{doi:10.2140/pjm.1960.10.419}}.

\bibitem{balay_petsc_2016}
Satish Balay, Shrirang Abhyankar, Mark~F Adams, Jed Brown, Peter Brune, Kris
  Buschelman, Lisandro Dalcin, Victor Eijkhout, William~D Gropp, Dinesh
  Kaushik, Matthew~G Knepley, Lois~Curfman McInnes, Karl Rupp, Barry~F Smith,
  Stefano Zampini, Hong Zhang, and Hong Zhang.
\newblock {{PETSc} {Users} {Manual}}.
\newblock Technical Report ANL-95/11 - Revision 3.7, Argonne National
  Laboratory, 2016.
\newblock URL: \url{http://www.mcs.anl.gov/petsc}.

\bibitem{Banjai2019}
Lehel Banjai, Jens~M Melenk, Ricardo~H Nochetto, Enrique Ot{\'{a}}rola, Abner~J
  Salgado, and Christoph Schwab.
\newblock {Tensor FEM for Spectral Fractional Diffusion}.
\newblock {\em Found. Comput. Math.}, 19(4):901--962, aug 2019.
\newblock \href {http://arxiv.org/abs/1707.07367} {\path{arXiv:1707.07367}},
  \href {https://doi.org/10.1007/s10208-018-9402-3}
  {\path{doi:10.1007/s10208-018-9402-3}}.

\bibitem{Banjai2020}
Lehel Banjai, Jens~M Melenk, and Christoph Schwab.
\newblock {Exponential convergence of hp FEM for spectral fractional diffusion
  in polygons}.
\newblock {\em arXiv}, pages 1--37, 2020.
\newblock \href {http://arxiv.org/abs/2011.05701} {\path{arXiv:2011.05701}}.

\bibitem{Bank1985}
Randolph~E Bank and Alan Weiser.
\newblock {Some A Posteriori Error Estimators for Elliptic Partial Differential
  Equations}.
\newblock {\em Math. Comput.}, 44(170):283, apr 1985.
\newblock \href {https://doi.org/10.2307/2007953} {\path{doi:10.2307/2007953}}.

\bibitem{Barrera2020}
O~Barrera.
\newblock A unified modelling and simulation for coupled anomalous transport in
  porous media and its finite element implementation.
\newblock {\em Computational Mechanics}, pages 1--16, 2021.

\bibitem{Bolin2017}
David Bolin, Kristin Kirchner, and Mih{\'{a}}ly Kov{\'{a}}cs.
\newblock {Numerical solution of fractional elliptic stochastic PDEs with
  spatial white noise}.
\newblock {\em IMA J. Numer. Anal.}, 40(2):1051--1073, apr 2020.
\newblock \href {http://arxiv.org/abs/1705.06565} {\path{arXiv:1705.06565}},
  \href {https://doi.org/10.1093/imanum/dry091}
  {\path{doi:10.1093/imanum/dry091}}.

\bibitem{Bonito2018a}
Andrea Bonito, Juan~Pablo Borthagaray, Ricardo~H Nochetto, Enrique
  Ot{\'{a}}rola, and Abner~J Salgado.
\newblock {Numerical methods for fractional diffusion}.
\newblock {\em Comput. Vis. Sci.}, 19(5-6):19--46, dec 2018.
\newblock \href {http://arxiv.org/abs/1707.01566} {\path{arXiv:1707.01566}},
  \href {https://doi.org/10.1007/s00791-018-0289-y}
  {\path{doi:10.1007/s00791-018-0289-y}}.

\bibitem{Bonito2017}
Andrea Bonito, Wenyu Lei, and Joseph~E Pasciak.
\newblock {Numerical Approximation of Space-Time Fractional Parabolic
  Equations}.
\newblock {\em Comput. Methods Appl. Math.}, 17(4):679--705, oct 2017.
\newblock \href {http://arxiv.org/abs/1704.04254} {\path{arXiv:1704.04254}},
  \href {https://doi.org/10.1515/cmam-2017-0032}
  {\path{doi:10.1515/cmam-2017-0032}}.

\bibitem{Bonito2016}
Andrea Bonito, Wenyu Lei, and Joseph~E Pasciak.
\newblock {The approximation of parabolic equations involving fractional powers
  of elliptic operators}.
\newblock {\em J. Comput. Appl. Math.}, 315:32--48, may 2017.
\newblock \href {http://arxiv.org/abs/1607.07832} {\path{arXiv:1607.07832}},
  \href {https://doi.org/10.1016/j.cam.2016.10.016}
  {\path{doi:10.1016/j.cam.2016.10.016}}.

\bibitem{Bonito2019a}
Andrea Bonito, Wenyu Lei, and Joseph~E Pasciak.
\newblock {Numerical approximation of the integral fractional Laplacian}.
\newblock {\em Numer. Math.}, 142(2):235--278, jun 2019.
\newblock \href {http://arxiv.org/abs/1707.04290} {\path{arXiv:1707.04290}},
  \href {https://doi.org/10.1007/s00211-019-01025-x}
  {\path{doi:10.1007/s00211-019-01025-x}}.

\bibitem{Bonito2019}
Andrea Bonito, Wenyu Lei, and Joseph~E Pasciak.
\newblock {On sinc quadrature approximations of fractional powers of regularly
  accretive operators}.
\newblock {\em J. Numer. Math.}, 27(2):57--68, jun 2019.
\newblock \href {http://arxiv.org/abs/1709.06619} {\path{arXiv:1709.06619}},
  \href {https://doi.org/10.1515/jnma-2017-0116}
  {\path{doi:10.1515/jnma-2017-0116}}.

\bibitem{Bonito2021a}
Andrea Bonito and Murtazo Nazarov.
\newblock {Numerical Simulations of Surface Quasi-Geostrophic Flows on Periodic
  Domains}.
\newblock {\em SIAM J. Sci. Comput.}, 43(2):B405--B430, jan 2021.
\newblock \href {http://arxiv.org/abs/2006.01180} {\path{arXiv:2006.01180}},
  \href {https://doi.org/10.1137/20M1342616} {\path{doi:10.1137/20M1342616}}.

\bibitem{Bonito2013}
Andrea Bonito and Joseph~E Pasciak.
\newblock {Numerical approximation of fractional powers of elliptic operators}.
\newblock {\em Math. Comput.}, 84(295):2083--2110, mar 2015.
\newblock \href {http://arxiv.org/abs/1307.0888} {\path{arXiv:1307.0888}},
  \href {https://doi.org/10.1090/S0025-5718-2015-02937-8}
  {\path{doi:10.1090/S0025-5718-2015-02937-8}}.

\bibitem{Bonito2017a}
Andrea Bonito and Joseph~E Pasciak.
\newblock {Numerical approximation of fractional powers of regularly accretive
  operators}.
\newblock {\em IMA J. Numer. Anal.}, page drw042, aug 2016.
\newblock \href {http://arxiv.org/abs/1508.05869} {\path{arXiv:1508.05869}},
  \href {https://doi.org/10.1093/imanum/drw042}
  {\path{doi:10.1093/imanum/drw042}}.

\bibitem{Bonito2020a}
Andrea Bonito and Peng Wei.
\newblock {Electroconvection of thin liquid crystals: Model reduction and
  numerical simulations}.
\newblock {\em J. Comput. Phys.}, 405:109140, mar 2020.
\newblock \href {https://doi.org/10.1016/j.jcp.2019.109140}
  {\path{doi:10.1016/j.jcp.2019.109140}}.

\bibitem{Bordas2016}
St{\'{e}}phane~P.A. Bordas, Sundararajan Natarajan, and Alexander Menk.
\newblock {\em {Partition of unity methods}}.
\newblock Wiley-blackwell edition, 2016.

\bibitem{Borthagaray2019a}
Juan Borthagaray, Wenbo Li, and Ricardo~H Nochetto.
\newblock {Linear and nonlinear fractional elliptic problems}.
\newblock pages 69--92, jun 2020.
\newblock \href {http://arxiv.org/abs/1906.04230} {\path{arXiv:1906.04230}},
  \href {https://doi.org/10.1090/conm/754/15145}
  {\path{doi:10.1090/conm/754/15145}}.

\bibitem{Borthagaray2020}
Juan~Pablo Borthagaray, Dmitriy Leykekhman, and Ricardo~H Nochetto.
\newblock {Local Energy Estimates for the Fractional Laplacian}.
\newblock {\em SIAM J. Numer. Anal.}, 59(4):1918--1947, jan 2021.
\newblock \href {http://arxiv.org/abs/2005.03786} {\path{arXiv:2005.03786}},
  \href {https://doi.org/10.1137/20M1335509} {\path{doi:10.1137/20M1335509}}.

\bibitem{Bulle2021}
Rapha{\"{e}}l Bulle, Gioacchino Alotta, Gregorio Marchiori, Matteo Berni,
  Nicola~F. Lopomo, Stefano Zaffagnini, St{\'{e}}phane P.~A. Bordas, and Olga
  Barrera.
\newblock {The Human Meniscus Behaves as a Functionally Graded Fractional
  Porous Medium under Confined Compression Conditions}.
\newblock {\em Appl. Sci.}, 11(20):9405, oct 2021.
\newblock \href {https://doi.org/10.3390/app11209405}
  {\path{doi:10.3390/app11209405}}.

\bibitem{Bulle2020}
Rapha{\"{e}}l Bulle, Franz Chouly, Jack~S Hale, and Alexei Lozinski.
\newblock {Removing the saturation assumption in Bank–Weiser error estimator
  analysis in dimension three}.
\newblock {\em Appl. Math. Lett.}, 107:106429, sep 2020.
\newblock \href {https://doi.org/10.1016/j.aml.2020.106429}
  {\path{doi:10.1016/j.aml.2020.106429}}.

\bibitem{bulle_fenics-ee_2019}
Rapha{\"{e}}l Bulle and Jack~S Hale.
\newblock {{FEniCS} {Error} {Estimation} {(FEniCS-EE)}}, jan 2019.
\newblock \href {https://doi.org/10.6084/m9.figshare.10732421}
  {\path{doi:10.6084/m9.figshare.10732421}}.

\bibitem{Bulle2021a}
Raphaël Bulle, Jack~S. Hale, Alexei Lozinski, Stéphane~P.A. Bordas, and Franz
  Chouly.
\newblock Hierarchical a posteriori error estimation of {Bank}–{Weiser} type
  in the {FEniCS} {Project}.
\newblock {\em Computers \& Mathematics with Applications}, 131:103--123,
  February 2023.
\newblock URL:
  \url{https://linkinghub.elsevier.com/retrieve/pii/S0898122122004722}, \href
  {https://doi.org/10.1016/j.camwa.2022.11.009}
  {\path{doi:10.1016/j.camwa.2022.11.009}}.

\bibitem{Caffarelli2007}
Luis Caffarelli and Luis Silvestre.
\newblock {An Extension Problem Related to the Fractional Laplacian}.
\newblock {\em Commun. Partial Differ. Equations}, 32(8):1245--1260, aug 2007.
\newblock \href {http://arxiv.org/abs/0608640} {\path{arXiv:0608640}}, \href
  {https://doi.org/10.1080/03605300600987306}
  {\path{doi:10.1080/03605300600987306}}.

\bibitem{Cances2017}
Eric Canc{\`{e}}s, Genevi{\`{e}}ve Dusson, Yvon Maday, Benjamin Stamm, and
  Martin Vohral{\'{i}}k.
\newblock {Guaranteed and Robust a Posteriori Bounds for Laplace Eigenvalues
  and Eigenvectors: Conforming Approximations}.
\newblock {\em SIAM J. Numer. Anal.}, 55(5):2228--2254, jan 2017.
\newblock \href {https://doi.org/10.1137/15M1038633}
  {\path{doi:10.1137/15M1038633}}.

\bibitem{Caputo1967}
Michele Caputo.
\newblock {Linear Models of Dissipation whose Q is almost Frequency
  Independent--II}.
\newblock {\em Geophys. J. Int.}, 13(5):529--539, nov 1967.
\newblock \href {https://doi.org/10.1111/j.1365-246X.1967.tb02303.x}
  {\path{doi:10.1111/j.1365-246X.1967.tb02303.x}}.

\bibitem{Caputo2000}
Michele Caputo.
\newblock {Models of flux in porous media with memory}.
\newblock {\em Water Resour. Res.}, 36(3):693--705, mar 2000.
\newblock \href {https://doi.org/10.1029/1999WR900299}
  {\path{doi:10.1029/1999WR900299}}.

\bibitem{Carlson2019}
Max Carlson, Robert~M Kirby, and Hari Sundar.
\newblock {A scalable framework for solving fractional diffusion equations}.
\newblock {\em Proc. 34th ACM Int. Conf. Supercomput.}, pages 1--11, jun 2020.
\newblock \href {http://arxiv.org/abs/1911.11906} {\path{arXiv:1911.11906}},
  \href {https://doi.org/10.1145/3392717.3392769}
  {\path{doi:10.1145/3392717.3392769}}.

\bibitem{Carstensen2014}
Carsten Carstensen and Joscha Gedicke.
\newblock {Guaranteed lower bounds for eigenvalues}.
\newblock {\em Math. Comput.}, 83(290):2605--2629, apr 2014.
\newblock \href {https://doi.org/10.1090/S0025-5718-2014-02833-0}
  {\path{doi:10.1090/S0025-5718-2014-02833-0}}.

\bibitem{Carstensen2010}
Carsten Carstensen and Christian Merdon.
\newblock {Estimator Competition for poisson Problems}.
\newblock {\em J. Comput. Math.}, 28(3):309--330, 2010.
\newblock \href {https://doi.org/10.4208/jcm.2009.10-m1010}
  {\path{doi:10.4208/jcm.2009.10-m1010}}.

\bibitem{Chen2018}
Huyuan Chen.
\newblock {The Dirichlet elliptic problem involving regional fractional
  Laplacian}.
\newblock {\em J. Math. Phys.}, 59(7):071504, jul 2018.
\newblock \href {http://arxiv.org/abs/1509.05838} {\path{arXiv:1509.05838}},
  \href {https://doi.org/10.1063/1.5046685} {\path{doi:10.1063/1.5046685}}.

\bibitem{Chen2015}
Long Chen, Ricardo~H Nochetto, Enrique Ot{\'{a}}rola, and Abner~J Salgado.
\newblock {A PDE approach to fractional diffusion: A posteriori error
  analysis}.
\newblock {\em J. Comput. Phys.}, 293:339--358, jul 2015.
\newblock \href {https://doi.org/10.1016/j.jcp.2015.01.001}
  {\path{doi:10.1016/j.jcp.2015.01.001}}.

\bibitem{Collier2015}
Nathan Collier, Abdul-Lateef Haji-Ali, Fabio Nobile, Erik von Schwerin, and
  Raúl Tempone.
\newblock A continuation multilevel {Monte} {Carlo} algorithm.
\newblock {\em BIT Numerical Mathematics}, 55(2):399--432, June 2015.
\newblock URL: \url{http://link.springer.com/10.1007/s10543-014-0511-3}, \href
  {https://doi.org/10.1007/s10543-014-0511-3}
  {\path{doi:10.1007/s10543-014-0511-3}}.

\bibitem{Cusimano2020}
Nicole Cusimano, F{\'{e}}lix del Teso, and Luca Gerardo-Giorda.
\newblock {Numerical approximations for fractional elliptic equations via the
  method of semigroups}.
\newblock {\em ESAIM Math. Model. Numer. Anal.}, 54(3):751--774, may 2020.
\newblock \href {https://doi.org/10.1051/m2an/2019076}
  {\path{doi:10.1051/m2an/2019076}}.

\bibitem{Cusimano2018}
Nicole Cusimano, F{\'{e}}lix del Teso, Luca Gerardo-Giorda, and Gianni Pagnini.
\newblock {Discretizations of the Spectral Fractional Laplacian on General
  Domains with Dirichlet, Neumann, and Robin Boundary Conditions}.
\newblock {\em SIAM J. Numer. Anal.}, pages 1243--1272, jan 2018.
\newblock \href {http://arxiv.org/abs/1708.03602} {\path{arXiv:1708.03602}},
  \href {https://doi.org/10.1137/17M1128010} {\path{doi:10.1137/17M1128010}}.

\bibitem{Danczul2019}
Tobias Danczul and Joachim Sch{\"{o}}berl.
\newblock {A Reduced Basis Method For Fractional Diffusion Operators I}.
\newblock pages 1--19, 2019.
\newblock URL: \url{http://arxiv.org/abs/1904.05599}, \href
  {http://arxiv.org/abs/1904.05599} {\path{arXiv:1904.05599}}.

\bibitem{Danczul2020}
Tobias Danczul and Joachim Sch{\"{o}}berl.
\newblock {A reduced basis method for fractional diffusion operators II}.
\newblock {\em J. Numer. Math.}, 29(4):269--287, oct 2021.
\newblock \href {http://arxiv.org/abs/2005.03574} {\path{arXiv:2005.03574}},
  \href {https://doi.org/10.1515/jnma-2020-0042}
  {\path{doi:10.1515/jnma-2020-0042}}.

\bibitem{Defterli2015}
Ozlem Defterli, Marta D'Elia, Qiang Du, Max Gunzburger, Rich Lehoucq, and
  Mark~M Meerschaert.
\newblock {Fractional Diffusion on Bounded Domains}.
\newblock {\em Fract. Calc. Appl. Anal.}, pages 1689--1699, jan 2015.
\newblock \href {http://arxiv.org/abs/arXiv:1011.1669v3}
  {\path{arXiv:arXiv:1011.1669v3}}, \href
  {https://doi.org/10.1515/fca-2015-0023} {\path{doi:10.1515/fca-2015-0023}}.

\bibitem{DElia2013}
Marta D'Elia and Max Gunzburger.
\newblock {The fractional Laplacian operator on bounded domains as a special
  case of the nonlocal diffusion operator}.
\newblock {\em Comput. Math. with Appl.}, pages 1245--1260, oct 2013.
\newblock \href {http://arxiv.org/abs/1303.6934} {\path{arXiv:1303.6934}},
  \href {https://doi.org/10.1016/j.camwa.2013.07.022}
  {\path{doi:10.1016/j.camwa.2013.07.022}}.

\bibitem{Dinh2021}
Huy Dinh, Harbir Antil, Yanlai Chen, Elena Cherkaev, and Akil Narayan.
\newblock {Model reduction for fractional elliptic problems using Kato's
  formula}.
\newblock {\em Math. Control Relat. Fields}, 2021.
\newblock \href {https://doi.org/10.3934/mcrf.2021004}
  {\path{doi:10.3934/mcrf.2021004}}.

\bibitem{dorfler_convergent_1996}
Willy D{\"{o}}rfler.
\newblock {A Convergent Adaptive Algorithm for Poisson's Equation}.
\newblock {\em SIAM J. Numer. Anal.}, 33(3):1106--1124, jun 1996.
\newblock \href {https://doi.org/10.1137/0733054} {\path{doi:10.1137/0733054}}.

\bibitem{Du2019}
Qiang Du, Jiang Yang, and Zhi Zhou.
\newblock {Time-Fractional Allen–Cahn Equations: Analysis and Numerical
  Methods}.
\newblock {\em J. Sci. Comput.}, page~42, nov 2020.
\newblock \href {http://arxiv.org/abs/1906.06584} {\path{arXiv:1906.06584}},
  \href {https://doi.org/10.1007/s10915-020-01351-5}
  {\path{doi:10.1007/s10915-020-01351-5}}.

\bibitem{Duo2019}
Siwei Duo, Hong Wang, and Yanzhi Zhang.
\newblock {A comparative study on nonlocal diffusion operators related to the
  fractional Laplacian}.
\newblock {\em Discret. Contin. Dyn. Syst. - B}, pages 231--256, 2019.
\newblock \href {http://arxiv.org/abs/1711.06916} {\path{arXiv:1711.06916}},
  \href {https://doi.org/10.3934/dcdsb.2018110}
  {\path{doi:10.3934/dcdsb.2018110}}.

\bibitem{falgout_hypre_2002}
Robert~D Falgout and Ulrike~Meier Yang.
\newblock {hypre: A Library of High Performance Preconditioners}.
\newblock In Peter M~A Sloot, Alfons~G Hoekstra, C~J~Kenneth Tan, and Jack~J
  Dongarra, editors, {\em Comput. {Science} — {ICCS} 2002}, number 2331 in
  Lecture {Notes} in {Computer} {Science}, pages 632--641. Springer Berlin
  Heidelberg, apr 2002.
\newblock \href {https://doi.org/10.1007/3-540-47789-6_66}
  {\path{doi:10.1007/3-540-47789-6_66}}.

\bibitem{Fall2020}
Mouhamed~M Fall.
\newblock {Regional fractional Laplacians: Boundary regularity}.
\newblock {\em arXiv}, jul 2020.
\newblock URL: \url{http://arxiv.org/abs/2007.04808}, \href
  {http://arxiv.org/abs/2007.04808} {\path{arXiv:2007.04808}}.

\bibitem{Faustmann2020}
Markus Faustmann, Michael Karkulik, and Jens~M Melenk.
\newblock {Local convergence of the FEM for the integral fractional Laplacian}.
\newblock {\em arXiv}, pages 1--20, may 2020.
\newblock URL: \url{http://arxiv.org/abs/2005.14109}, \href
  {http://arxiv.org/abs/2005.14109} {\path{arXiv:2005.14109}}.

\bibitem{Faustmann2019}
Markus Faustmann, Jens~M Melenk, and Dirk Praetorius.
\newblock {Quasi-optimal convergence rate for an adaptive method for the
  integral fractional Laplacian}.
\newblock {\em Math. Comput.}, 90(330):1557--1587, apr 2021.
\newblock \href {http://arxiv.org/abs/1903.10409} {\path{arXiv:1903.10409}},
  \href {https://doi.org/10.1090/mcom/3603} {\path{doi:10.1090/mcom/3603}}.

\bibitem{Gavrilyuk2003}
Ivan~P Gavrilyuk, Wolfgang Hackbusch, and Boris~N Khoromskij.
\newblock {Data-sparse approximation to the operator-valued functions of
  elliptic operator}.
\newblock {\em Math. Comput.}, pages 1297--1325, jul 2003.
\newblock \href {https://doi.org/10.1090/S0025-5718-03-01590-4}
  {\path{doi:10.1090/S0025-5718-03-01590-4}}.

\bibitem{Gavrilyuk2004}
Ivan~P Gavrilyuk, Wolfgang Hackbusch, and Boris~N Khoromskij.
\newblock {Data-sparse approximation to a class of operator-valued functions}.
\newblock {\em Math. Comput.}, pages 681--709, aug 2004.
\newblock \href {https://doi.org/10.1090/S0025-5718-04-01703-X}
  {\path{doi:10.1090/S0025-5718-04-01703-X}}.

\bibitem{geuzaine_gmsh:2009}
Christophe Geuzaine and Jean-François Remacle.
\newblock Gmsh: {A} 3-{D} finite element mesh generator with built-in pre- and
  post-processing facilities.
\newblock {\em International Journal for Numerical Methods in Engineering},
  79(11):1309--1331, September 2009.
\newblock URL:
  \url{http://onlinelibrary.wiley.com/doi/10.1002/nme.2579/abstract}, \href
  {https://doi.org/10.1002/nme.2579} {\path{doi:10.1002/nme.2579}}.

\bibitem{Gimperlein2019}
Heiko Gimperlein and Jakub Stocek.
\newblock {Space–time adaptive finite elements for nonlocal parabolic
  variational inequalities}.
\newblock {\em Comput. Methods Appl. Mech. Eng.}, pages 137--171, aug 2019.
\newblock \href {http://arxiv.org/abs/1810.06888} {\path{arXiv:1810.06888}},
  \href {https://doi.org/10.1016/j.cma.2019.04.019}
  {\path{doi:10.1016/j.cma.2019.04.019}}.

\bibitem{Grubb2015}
Gerd Grubb.
\newblock {Fractional Laplacians on domains, a development of H{\"{o}}rmander's
  theory of $\mu$-transmission pseudodifferential operators}.
\newblock {\em Adv. Math. (N. Y).}, 268:478--528, jan 2015.
\newblock \href {https://doi.org/10.1016/j.aim.2014.09.018}
  {\path{doi:10.1016/j.aim.2014.09.018}}.

\bibitem{Habera2020}
Michal Habera, Jack~S Hale, Chris Richardson, Johannes Ring, Marie Rognes, Nate
  Sime, and Garth~N Wells.
\newblock {FEniCSX: A sustainable future for the FEniCS Project}.
\newblock 2020.
\newblock \href {https://doi.org/10.6084/m9.figshare.11866101.v1}
  {\path{doi:10.6084/m9.figshare.11866101.v1}}.

\bibitem{hale_containers_2017}
Jack~S. Hale, Lizao Li, Christopher~N. Richardson, and Garth~N. Wells.
\newblock Containers for {Portable}, {Productive}, and {Performant}
  {Scientific} {Computing}.
\newblock {\em Computing in Science \& Engineering}, 19(6):40--50, November
  2017.
\newblock \href {https://doi.org/10.1109/MCSE.2017.2421459}
  {\path{doi:10.1109/MCSE.2017.2421459}}.

\bibitem{Harizanov2016}
Stanislav Harizanov, Raytcho Lazarov, Svetozar Margenov, Pencho Marinov, and
  Yavor Vutov.
\newblock {Optimal solvers for linear systems with fractional powers of sparse
  SPD matrices}.
\newblock {\em Numer. Linear Algebr. with Appl.}, 25(5):e2167, oct 2018.
\newblock \href {http://arxiv.org/abs/1612.04846} {\path{arXiv:1612.04846}},
  \href {https://doi.org/10.1002/nla.2167} {\path{doi:10.1002/nla.2167}}.

\bibitem{Harizanov2020}
Stanislav Harizanov, Svetozar Margenov, and Nedyu Popivanov.
\newblock {Spectral Fractional Laplacian with Inhomogeneous Dirichlet Data:
  Questions, Problems, Solutions}.
\newblock pages 123--138, 2021.
\newblock \href {http://arxiv.org/abs/arXiv:2010.01383v1}
  {\path{arXiv:arXiv:2010.01383v1}}, \href
  {https://doi.org/10.1007/978-3-030-71616-5_13}
  {\path{doi:10.1007/978-3-030-71616-5_13}}.

\bibitem{Higham2013}
Nicholas~J. Higham and Lijing Lin.
\newblock {An Improved Schur--Pad{\'{e}} Algorithm for Fractional Powers of a
  Matrix and Their Fr{\'{e}}chet Derivatives}.
\newblock {\em SIAM J. Matrix Anal. Appl.}, 34(3):1341--1360, jan 2013.
\newblock \href {https://doi.org/10.1137/130906118}
  {\path{doi:10.1137/130906118}}.

\bibitem{Hofreither2020}
Clemens Hofreither.
\newblock {A unified view of some numerical methods for fractional diffusion}.
\newblock {\em Comput. Math. with Appl.}, 80(2):332--350, jul 2020.
\newblock \href {https://doi.org/10.1016/j.camwa.2019.07.025}
  {\path{doi:10.1016/j.camwa.2019.07.025}}.

\bibitem{Hofreither2021}
Clemens Hofreither.
\newblock {An algorithm for best rational approximation based on barycentric
  rational interpolation}.
\newblock {\em Numer. Algorithms}, 88(1):365--388, sep 2021.
\newblock \href {https://doi.org/10.1007/s11075-020-01042-0}
  {\path{doi:10.1007/s11075-020-01042-0}}.

\bibitem{Kwasnicki2017}
Mateusz Kwa{\'{s}}nicki.
\newblock {Ten equivalent definitions of the fractional laplace operator}.
\newblock {\em Fract. Calc. Appl. Anal.}, 20(1):7--51, jan 2017.
\newblock \href {http://arxiv.org/abs/1507.07356} {\path{arXiv:1507.07356}},
  \href {https://doi.org/10.1515/fca-2017-0002}
  {\path{doi:10.1515/fca-2017-0002}}.

\bibitem{Lindgren2011a}
Finn Lindgren, H{\aa}vard Rue, and Johan Lindstr{\"{o}}m.
\newblock {An explicit link between Gaussian fields and Gaussian Markov random
  fields: the stochastic partial differential equation approach}.
\newblock {\em J. R. Stat. Soc. Ser. B (Statistical Methodol.)},
  73(4):423--498, sep 2011.
\newblock \href {https://doi.org/10.1111/j.1467-9868.2011.00777.x}
  {\path{doi:10.1111/j.1467-9868.2011.00777.x}}.

\bibitem{Lischke2018}
Anna Lischke, Guofei Pang, Mamikon Gulian, Fangying Song, Christian Glusa,
  Xiaoning Zheng, Zhiping Mao, Wei Cai, Mark~M Meerschaert, Mark Ainsworth, and
  George~Em Karniadakis.
\newblock {What is the fractional Laplacian? A comparative review with new
  results}.
\newblock {\em J. Comput. Phys.}, 404:109009, mar 2020.
\newblock \href {http://arxiv.org/abs/1801.09767} {\path{arXiv:1801.09767}},
  \href {https://doi.org/10.1016/j.jcp.2019.109009}
  {\path{doi:10.1016/j.jcp.2019.109009}}.

\bibitem{Meidner2018}
Dominik Meidner, Johannes Pfefferer, Klemens Sch{\"{u}}rholz, and Boris Vexler.
\newblock {$hp$-Finite Elements for Fractional Diffusion}.
\newblock {\em SIAM J. Numer. Anal.}, 56(4):2345--2374, jan 2018.
\newblock \href {http://arxiv.org/abs/1706.04066} {\path{arXiv:1706.04066}},
  \href {https://doi.org/10.1137/17M1135517} {\path{doi:10.1137/17M1135517}}.

\bibitem{Mou2015}
Chenchen Mou and Yingfei Yi.
\newblock {Interior Regularity for Regional Fractional Laplacian}.
\newblock {\em Commun. Math. Phys.}, 340(1):233--251, nov 2015.
\newblock \href {https://doi.org/10.1007/s00220-015-2445-2}
  {\path{doi:10.1007/s00220-015-2445-2}}.

\bibitem{Nochetto2015}
Ricardo~H Nochetto, Enrique Ot{\'{a}}rola, and Abner~J Salgado.
\newblock {A PDE Approach to Fractional Diffusion in General Domains: A Priori
  Error Analysis}.
\newblock {\em Found. Comput. Math.}, 15(3):733--791, jun 2015.
\newblock \href {http://arxiv.org/abs/1302.0698} {\path{arXiv:1302.0698}},
  \href {https://doi.org/10.1007/s10208-014-9208-x}
  {\path{doi:10.1007/s10208-014-9208-x}}.

\bibitem{Nochetto2010}
Ricardo~H Nochetto, Tobias von Petersdorff, and Chen-Song Zhang.
\newblock {A posteriori error analysis for a class of integral equations and
  variational inequalities}.
\newblock {\em Numer. Math.}, 116(3):519--552, sep 2010.
\newblock \href {https://doi.org/10.1007/s00211-010-0310-y}
  {\path{doi:10.1007/s00211-010-0310-y}}.

\bibitem{plaza_local_2000}
Angel Plaza and Graham~F Carey.
\newblock {Local refinement of simplicial grids based on the skeleton}.
\newblock {\em Appl. Numer. Math.}, 32(2):195--218, feb 2000.
\newblock \href {https://doi.org/10.1016/S0168-9274(99)00022-7}
  {\path{doi:10.1016/S0168-9274(99)00022-7}}.

\bibitem{Podlubny2009}
Igor Podlubny, Aleksei Chechkin, Tomas Skovranek, YangQuan Chen, and Blas~M.
  {Vinagre Jara}.
\newblock {Matrix approach to discrete fractional calculus II: Partial
  fractional differential equations}.
\newblock {\em J. Comput. Phys.}, 228(8):3137--3153, may 2009.
\newblock \href {https://doi.org/10.1016/j.jcp.2009.01.014}
  {\path{doi:10.1016/j.jcp.2009.01.014}}.

\bibitem{Stahl2003}
Herbert~R. Stahl.
\newblock {Best uniform rational approximation of $x^{\alpha}$ on $[0, 1]$}.
\newblock {\em Acta Mathematica}, 190(2):241 -- 306, 2003.
\newblock \href {https://doi.org/10.1007/BF02392691}
  {\path{doi:10.1007/BF02392691}}.

\bibitem{Stinga2019}
Pablo~Ra{\'{u}}l Stinga.
\newblock {User's guide to the fractional Laplacian and the method of
  semigroups}.
\newblock In {\em Fract. Differ. Equations}, pages 235--266. De Gruyter, feb
  2019.
\newblock \href {http://arxiv.org/abs/1808.05159} {\path{arXiv:1808.05159}},
  \href {https://doi.org/10.1515/9783110571660-012}
  {\path{doi:10.1515/9783110571660-012}}.

\bibitem{Stinga2010}
Pablo~Ra{\'{u}}l Stinga and Jos{\'{e}}~Luis Torrea.
\newblock {Extension Problem and Harnack's Inequality for Some Fractional
  Operators}.
\newblock {\em Commun. Partial Differ. Equations}, 35(11):2092--2122, oct 2010.
\newblock \href {http://arxiv.org/abs/0910.2569} {\path{arXiv:0910.2569}},
  \href {https://doi.org/10.1080/03605301003735680}
  {\path{doi:10.1080/03605301003735680}}.

\bibitem{Vabishchevich2019}
Petr~N Vabishchevich.
\newblock {Approximation of a fractional power of an elliptic operator}.
\newblock {\em Numer. Linear Algebr. with Appl.}, 27(3):1--17, may 2020.
\newblock \href {http://arxiv.org/abs/1905.10838} {\path{arXiv:1905.10838}},
  \href {https://doi.org/10.1002/nla.2287} {\path{doi:10.1002/nla.2287}}.

\bibitem{Verfurth1998}
R{\"{u}}diger Verf{\"{u}}rth.
\newblock {Robust a posteriori error estimators for a singularly perturbed
  reaction-diffusion equation}.
\newblock {\em Numer. Math.}, 78(3):479--493, jan 1998.
\newblock \href {https://doi.org/10.1007/s002110050322}
  {\path{doi:10.1007/s002110050322}}.

\bibitem{Weiss2019}
Christian~J Weiss, Bart~G {van Bloemen Waanders}, and Habir Antil.
\newblock {Fractional Operators Applied to Geophysical Electromagnetics}.
\newblock {\em Geophys. J. Int.}, pages 1242--1259, nov 2019.
\newblock \href {http://arxiv.org/abs/1902.05096} {\path{arXiv:1902.05096}},
  \href {https://doi.org/10.1093/gji/ggz516} {\path{doi:10.1093/gji/ggz516}}.

\bibitem{Zhao2017}
Xuan Zhao, Xiaozhe Hu, Wei Cai, and George~Em Karniadakis.
\newblock {Adaptive finite element method for fractional differential equations
  using hierarchical matrices}.
\newblock {\em Comput. Methods Appl. Mech. Eng.}, 325:56--76, oct 2017.
\newblock \href {http://arxiv.org/abs/1603.01358} {\path{arXiv:1603.01358}},
  \href {https://doi.org/10.1016/j.cma.2017.06.017}
  {\path{doi:10.1016/j.cma.2017.06.017}}.

\end{thebibliography}
\end{document}